\documentclass{article}
\usepackage{temp}

\usepackage[utf8]{inputenc} % allow utf-8 input
\usepackage[T1]{fontenc}    % use 8-bit T1 fonts
\usepackage{hyperref}       % hyperlinks
\usepackage{url}            % simple URL typesetting
\usepackage{booktabs}  % professional-quality tables
\usepackage{amsmath}
\usepackage{amsthm}
\usepackage{amsfonts}       % blackboard math symbols
\usepackage{nicefrac}       % compact symbols for 1/2, etc.
\usepackage{microtype}      % microtypography
\usepackage{cleveref}       % smart cross-referencing
\usepackage{lipsum}         % Can be removed after putting your text content
\usepackage{graphicx}
\usepackage{float}
\usepackage{subcaption}
\usepackage{doi}
\usepackage{bm}
\usepackage{algorithm}
\usepackage{algorithmic}
\usepackage{makecell}
\usepackage{changepage}
\theoremstyle{plain}
\newtheorem{theorem}{Theorem}[section]
\newtheorem{lemma}[theorem]{Lemma}

\theoremstyle{definition}

\theoremstyle{remark}

\usepackage{mathtools}
\title{Global Energy Minimization for Simplex Mesh Optimization: A Radius
Ratio Approach to Sliver Elimination}

\date{}

\author{{Dong Wang} \\
    School of Mathematics and Computational Science \\
	XiangTan University \\
    China, Xiangtan 411105 \\
    wangdongs@smail.xtu.edu.cn 
	\And
    {Chunyu Chen} \\
    School of Mathematics and Computational Science \\
	XiangTan University \\
    China, Xiangtan 411105 \\
    cbtxs@smail.xtu.edu.cn 
	\And
    {Huayi Wei*} \\
    School of Mathematics and Computational Science, Xiangtan University \\
	National Center of Applied Mathematics in Hunan \\
    Hunan Key Laboratory for Computation and Simulation 
    in Science and Engineering \\
    China, Xiangtan 411105 \\
    weihuayi@xtu.edu.cn 
}

\usepackage[square,sort,comma,numbers]{natbib}
\begin{document}
\maketitle

\begin{abstract}
This paper constructs an energy function for simplex mesh based on the radius 
ratio and develops a corresponding mesh optimization method. The method 
combines vertex relocation and connectivity improvement, and can effectively 
remove slivers and improve the overall mesh quality. Based on the structure of 
the gradient of the energy function, we design a preconditioner, which reduces 
the number of iterations and improves the efficiency of the optimization 
algorithm. Numerical experiments show that the proposed method is effective 
in both sliver removal and mesh quality improvement.
\end{abstract}

% keywords can be removed

\section{Introduction}
In numerical computation, mesh plays an important role. Mainstream 
simulation methods, including the Finite Element Method (FEM) and the Finite 
Volume Method (FVM), are based on mesh generation~\cite{13ref,14ref}. 
As basic geometric entities, simplices are widely used in many fields, and the 
study of simplex mesh generation has continued for several decades~\cite{15ref}. 
Current methods primarily fall into three categories: the 
Advancing-Front technique~\cite{31ref}, spatial decomposition
algorithms~\cite{32ref,33ref}, and Delaunay-based methods~\cite{2ref,34ref}.

Delaunay algorithms are widely used in mesh generation due to its algorithmic 
simplicity and broad applicability. However, in three-dimensional cases, it 
tends to produce ill-shaped tetrahedron that require optimization. 
Common improvement strategies include Laplacian smoothing,  local transformations, 
and hybrid approaches~\cite{5ref,7ref},  However, these methods cannot completely 
remove a type of nearly degenerate element called sliver.The existence of
slivers was first pointed out by Cavendish~\cite{1ref}. A sliver has 
faces that are well-shaped triangles, but its volume is close to zero, and its 
four vertices are almost coplanar. Despite this, it still satisfies the Delaunay 
property in a topological sense.  Such elements can seriously affect the 
numerical simulation~\cite{43ref}. 
To address the sliver problem, many methods have been proposed. For example, Chew~\cite{16ref} 
introduced a randomized point insertion method, Li~\cite{42sliver} proposed a 
vertex perturbation method, and Cheng~\cite{38ref} developed the sliver exudation 
method.
\begin{figure}[htbp]
\centering
\subfloat[Delaunay Tetrahedron and its circumsphere]{
\includegraphics[width=0.4\linewidth]{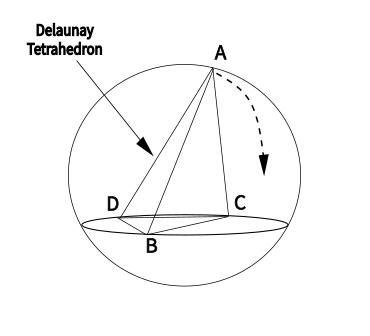}}
\hspace{0.01\linewidth}
\subfloat[Sliver Element is a valid Delaunay Element]{
\includegraphics[width=0.4\linewidth]{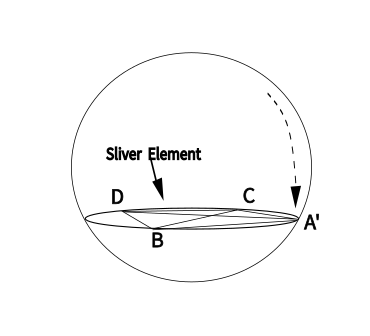}}
\caption{Sliver element}
\label{fig:sliver}
\end{figure}

In the past two decades, researchers have also studied mesh optimization from 
the perspective of optimization theory. Examples include variational mesh 
optimization methods such as Optimal Delaunay Triangulation(ODT)~\cite{19ref,40ref}
and Centroidal Voronoi Tessellation(CVT)~\cite{11ref,12ref} algorithms, as well
as algebraic quality-based methods such as the Target-Matrix Optimization
Paradigm(TMOP)~\cite{36ref,37ref}. Although these methods differ in their design 
ideas, they all define an energy function and formulate mesh optimization as 
the minimization of an objective function. A high-quality mesh is then obtained 
by applying numerical optimization methods. However, these approaches are not 
designed specifically for slivers. As a result, the optimized meshes 
may still contain slivers. Some researchers have designed energy functions that 
directly target slivers. For example, Saifeng Ni~\cite{38ref} proposed a 
shape-matching strategy to construct an energy function. This energy is 
sensitive to the height of tetrahedral elements and can effectively suppress slivers.

In addition to element height, the radius ratio is also an important quality 
measure that can effectively identify slivers. In fact, when Cavendish
first pointed out the existence of sliver, they introduced a quality metric 
based on the ratio between the inradius and the circumradius, and used it to 
determine whether an element is a sliver. This measure can capture the 
geometric feature that the element volume approaches zero while the edge 
lengths remain normal. Since then, it has been widely used for the detection 
and evaluation of slivers, and has become a classical quality metric in
this area. In this work, we introduce the radius ratio 
explicitly into the construction of the energy function, in order to build an 
optimization model that is highly sensitive to slivers. In addition, 
based on the structure of the gradient of the energy function, we construct a 
symmetric positive definite (SPD) and diagonally dominant matrix as a 
preconditioner, which leads to an efficient global optimization algorithm. The 
proposed method can be used as a standalone mesh optimization tool, and it can 
also be applied as a post-processing step for existing optimization methods.
\section{Radius Ratio Energy Function}\label{sec:rro}
The convex hull of a set of points 
$
\{\bm{x}_i\}_{i=0}^d\subset\mathbb{R}^d
$ 
that do not lie in one hyperplane
$$
\tau=\{\bm x =\sum_{i=0}^d\lambda_i\bm x_i|0\le\lambda_i\le 1,\sum_{i=0}^d\lambda_i=1\}
$$
is called a geometric d-simplex generated by $\{\bm x_i\}_{i=0}^d$. For example, 
an interval is a 1-simplex, a triangle is a 2-simplex, and a tetrahedron is a 
3-simplex.

For a positively-oriented simplex $\tau$ in $\mathbb{R}^d$, $d=2,3$, the radius 
ratio metric is defined as follows~\cite{1ref}:
$$
\mu=\frac{R}{dr}
$$
where $R$ is the circumradius and $r$ is the inradius of $\tau$. $\mu\in[1,+\infty]$ 
With $\mu=1$ if and only if $\tau$ is equilateral simplex.

Radius ratio metric $\mu$ is a function of the coordinates of the simplex vertices.
Given a triangulation $\mathcal T$ with $N_v$ vertices $\{\bm{x}_i\}_{i=0}^{N_v-1}$ and 
$N_c$ simplex $\{\tau_n\}_{n=0}^{N_c-1}$. We can define an energy function of radius 
ratio metric on the $\mathcal{T}$
$$
F=\frac{1}{N_c}\sum_{n=0}^{N_c-1}\mu_n
$$

Here we call $F$ is global radius ratio energy function. By minimizing $F$, we 
can develop efficient global mesh smoothing algorithm by moving the 
points in $\mathcal T$. The key problem is how to compute the gradient of $F$ 
about mesh point coordinates. In the following, we derive the gradient of
$\mu_n$ with respect to each vertex of a d-simplex $\tau_n$, and then obtain the 
global gradient.

\subsection{2D Case}
Given a 2D triangulation $\mathcal{T}$ with $N_v$ vertices
$\{\bm{x}_i\}_{i=0}^{N_v-1}\subset\mathbb{R}^2$ 
and $N_c$ triangles $\{\tau_n\}_{n=0}^{N_c-1}$. For $\forall\tau_n=(\bm{x}_0, \bm{x}_1,
\bm{x}_2)$, let $e_0:=(\bm{x}_1,\bm{x}_2), e_1:=(\bm{x}_2,\bm{x}_0)$ and 
$e_2:=(\bm{x}_0,\bm{x}_1)$ be the three edges of $\tau_n$ with length $l_0,l_1$ 
and $l_2$ respectively; $|\tau_n|$ be the area of $\tau_n$.

The circumradius of $\tau_n$ is:
$$
R_n=\frac{q}{4|\tau_n|}
$$
with $q=l_0l_1l_2$. The inradius of $\tau_n$ is:
$$
r_n=\frac{2|\tau_n|}{p}
$$
with $p=l_0+l_1+l_2$. The radius ratio metric of $\tau_n$ is:
$$
\mu_n = \frac{R_n}{2r_n} =\frac{pq}{16|\tau_n|^2}
$$
Next we compute the gradient $\mu_n$ about $\bm{x}_0$.
$$
\nabla_{\bm x_0} \mu_n 
= \frac{pq}{16|\tau_n|^2}(\frac{\nabla_{\bm x_0}p}{p}+
\frac{\nabla_{\bm x_0}q}{q}-\frac{2}{|\tau_n|}\nabla_{\bm x_0}|\tau_n|)
$$
where
$$
\begin{aligned}
\nabla_{\bm x_0} p &=\frac{1}{l_1}(\bm x_0-\bm x_2) + \frac{1}{l_2}(\bm x_0-\bm
x_1) \\
\nabla_{\bm x_0}q &= \frac{q}{l_1^2}(\bm x_0-\bm x_2)+\frac{q}{l_2^2}(\bm
x_0-\bm x_1) \\
\nabla_{\bm x_0}|\tau_n| &=\frac{1}{2}\bm W(\bm x_2-\bm x_1)
\end{aligned}
$$
with
$$
\bm W =
\begin{pmatrix}
0 & -1 \\ 1 & 0
\end{pmatrix}
$$
We have
$$
\nabla_{\bm x_0}\mu_n = 
\mu_n\left\{\left[\frac{1}{pl_1}+\frac{1}{l_1^2}\right](\bm x_0 - \bm x_2)+
\left[\frac{1}{pl_2}+\frac{1}{l_2^2}\right](\bm x_0 - \bm x_1)+ 
\frac{1}{|\tau_n|}\bm W(\bm x_1- \bm x_2)\right\}
$$
For simplicity, we denote
$$
c=\frac{1}{|\tau_n|}, c_0 = \frac{1}{pl_0}+\frac{1}{l_0^2}, c_1
= \frac{1}{pl_1}+\frac{1}{l_1^2}, c_2=\frac{1}{pl_2}+\frac{1}{l_2^2},
$$
$$
c_{01}=c_0+c_1, c_{02}=c_0+c_2, c_{12}=c_1+c_2
$$
So
$$
\nabla_{\bm{x}_0}\mu_n=
\mu_n\left[
c_{12}\bm x_0-
(c_2\bm I-c\bm W)\bm x_1-
(c_1\bm I+c\bm W)\bm x_2
\right]
$$
Simiarly
$$
\begin{aligned}
\nabla_{\bm x_1}\mu_n &=
\mu_n\left[
c_{20}\bm x_1-
(c_0\bm I-c\bm W)\bm x_2-
(c_2\bm I+c\bm W)\bm x_0
\right] \\
\nabla_{\bm x_2}\mu_n &=
\mu_n\left[
c_{01}\bm x_2-
(c_1\bm I-c\bm W)\bm x_0-(c_0\bm I+c\bm W)\bm x_1
\right]
\end{aligned}
$$
We denote
$$
\bm{V}_n = [x_0,x_1,x_2,y_0,y_1,y_2]^T
$$
Then we can write it in matrix form
$$
\nabla_{\bm{V}_n}\mu_n = 
\begin{bmatrix}
\bm{A}_n & \bm{B}_n\\
-\bm{B}_n & \bm{A}_n
\end{bmatrix} \bm{V}_n
$$
with
$$
\bm{A}_n=\mu_n 
\begin{pmatrix}
c_{12} & -c_{2} & -c_{1} \\
-c_{2} & c_{20} & -c_{0} \\
-c_{1} & -c_{0} & c_{01}
\end{pmatrix},
\bm{B}_n = \mu_n
\begin{pmatrix}
0 & -c & c \\
c & 0 & -c \\
-c & c & 0
\end{pmatrix}
$$
Matrix $\bm{A}_n$ is Laplacian and $\bm{B}_n$ is antisymmetric.

Then, we can assemble the gradient matrix of each element into a global 
gradient matrix. Therefore, we get the following theorem
\begin{theorem}
\label{thm:2dgrad}
For the energy function
$$
F=\frac{1}{N_c}\sum_{n=0}^{N_c-1}\mu_n
$$
of the triangular mesh, we have
$$
\nabla_{\bm{V}}F=
\begin{bmatrix}
\bm{A} & \bm{B} \\
-\bm{B} & \bm{A} \\
\end{bmatrix}
\begin{bmatrix}
\bm{X} \\
\bm{Y}
\end{bmatrix}
=
\bm{G}_F\bm{V}
$$ 
where $\bm{A}$ is $N_v\times N_v$ Laplacian matrix, and $\bm{B}$ is
$N_v\times N_v$ antisymmetric matrix,
$$
\bm{V} = [\bm{X}^T,\bm{Y}^T]^T,
\bm{X}=[x_0,x_1,\cdots,x_{N_v-1}]^T,\bm{Y}=[y_0,y_1,\cdots,y_{N_v-1}]^T
$$
\end{theorem}

\subsection{3D Case}
Given a 3D Tetrahedralization $\mathcal{T}$, consisting of $N_v$ vertices
$\{\bm{x}_i\}_{i=0}^{N_v-1}\subset\mathbb{R}^3$ and $N_c$ tetrahedrons
$\{\tau_n\}_{n=0}^{N_c-1}$. For $\forall\tau_n=(\bm{x}_0,\bm{x}_1,\bm{x}_2,\bm{x}_3)$,
let $s_0,s_1,s_2,s_3$ be the area of triangle opposite to vertices 
$\bm{x}_0,\bm{x}_1,\bm{x}_2,\bm{x}_3$; $|\tau_n|$ be the volume of $\tau_n$.

The radius of circumsphere of $\tau_n$ is:
$$
R_n = \frac{|\bm d_0|}{12|\tau_n|}
$$
where 
$$
\bm{d}_0 =\bm{v}_{30}^2\bm {v}_{10}\times\bm{v}_{20}+ 
\bm{v}_{10}^2\bm{v}_{20}\times\bm{v}_{30}+ \bm{v}_{20}^2\bm{v}_{30}\times\bm{v}_{10}
$$
with 
$$
\bm v_{10}=\bm x_0-\bm x_1, \bm v_{20}=\bm x_0-\bm x_2, \bm v_{30}=\bm x_0-\bm
x_3, \bm{v}^2=\bm{v}\cdot\bm{v}
$$
The radius of inside sphere of $\tau_n$ is
$$
r_n = \frac{3|\tau_n|}{s}
$$
with $s=s_0+s_1+s_2+s_3$. The radius ratio metric of $\tau_n$ is
$$
\mu_n =\frac{R_n}{3r_n}=\frac{s|\bm d_0|}{108|\tau_n|^2}
$$

The gradient of $\mu_n$ at vertex $\bm{x}_0$ is 
$$
\nabla_{\bm x_0} \mu_n = \mu_n\left\{\frac{1}{|\bm{d}_0|}\nabla_{\bm x_0}|\bm{d}_0|+
\frac{1}{s}\nabla_{\bm x_0}s-\frac{2}{|\tau_n|}\nabla_{\bm x_0}|\tau_n|\right\}
$$
where
$$
\begin{aligned}
\nabla_{\bm x_0}|\bm d_0|=&
\frac{1}{|\bm d_0|}\{2[\bm d_0\cdot(\bm v_{20}\times\bm v_{30})\bm v_{10}+
\bm d_0\cdot(\bm v_{30}\times\bm v_{10})\bm{v}_{20}+
\bm{d}_0\cdot(\bm{v}_{10}\times\bm v_{20})\bm{v}_{30}] \\
&+\bm{d}_0\times[(\bm v_{30}^2-\bm{v}_{20}^2)\bm{v}_{10}+
(\bm v_{10}^2-\bm v_{30}^2)\bm v_{20}+(\bm v_{20}^2-\bm v_{10}^2)\bm{v}_{30}]\}
\end{aligned}
$$
$$
\nabla_{\bm x_0} s =
(\frac{\bm{v}_{31}\cdot\bm{v}_{30}}{4s_2}+\frac{\bm{v}_{21}\cdot\bm{v}_{20}}{4s_3})\bm{v}_{10}+
(\frac{\bm{v}_{32}\cdot\bm{v}_{30}}{4s_1}+\frac{\bm{v}_{12}\cdot\bm{v}_{10}}{4s_3})\bm{v}_{20}+ 
(\frac{\bm{v}_{23}\cdot\bm{v}_{20}}{4s_1}+\frac{\bm{v}_{13}\cdot\bm{v}_{10}}{4s_2})\bm{v}_{30}
$$
$$
\nabla_{\bm x_0}|\tau_n|= -\frac{1}{18}\left\{[\bm x_2 - 2\bm x_3 + \bm x_1 ]_\times \bm{v}_{10} +
[\bm x_3 - 2\bm x_1 + \bm x_2 ]_\times \bm{v}_{20} + 
[\bm x_1  - 2\bm x_2 + \bm x_3]_\times \bm{v}_{30}\right\}
$$
We denote
$$
\bm{V}_n=
[x_0,x_1,x_2,x_3,y_0,y_1,y_2,y_3,z_0,z_1,z_2,z_3]^T
$$
and calculate the matrix form of 
$\nabla_{\bm{x}_i}|\bm{d}_0|,\nabla_{\bm{x}_i}s,\nabla_{\bm{x}_i}|\tau_n|~(i=0,1,2,3)$.

For $\nabla_{\bm{x}_i}|\bm{d}_0|$, we denote
$$
c_{23}=\bm{d}_0\cdot(\bm{v}_{20}\times\bm{v}_{30}),
c_{31}=\bm{d}_0\cdot(\bm{v}_{30}\times\bm{v}_{10}),
c_{12}=\bm{d}_0\cdot(\bm{v}_{10}\times\bm{v}_{20})
$$
$$
k_{23}=\bm{v}_{30}^2-\bm{v}_{20}^2,
k_{31}=\bm{v}_{10}^2-\bm{v}_{30}^2,
k_{12}=\bm{v}_{20}^2-\bm{v}_{10}^2
$$
$$
[\bm{d}_0]_{\times}
=\bm{D}_0=
\begin{bmatrix}
0 & -d_2 & d_1 \\
d_2 & 0 & -d_0 \\
-d_1 & d_0 & 0
\end{bmatrix}
$$
so we have 
$$
\nabla_{\bm{x}_0}|\bm{d}_0|=
\frac{1}{12|\tau_n|R_n}
\{[2(c_{23}+c_{31}+c_{12})+(k_{23}+k_{31}+k_{12})\bm{D}_0]\bm{x}_0-
(2c_{23}+k_{23}\bm{D}_0)\bm{x}_1-
(2c_{31}+k_{31}\bm{D}_0)\bm{x}_2-
(2c_{12}+k_{12}\bm{D}_0)\bm{x}_3\}
$$

$$
\begin{aligned}
\nabla_{\bm{x}_1}|\bm{d}_0| &=
\frac{1}{12|\tau_n|R_n}\{(-2c_{23}+k_{23}\bm{D}_0)\bm{x}_0+2c_{23}\bm{x}_1-
\bm{v}_{30}^2\bm{D}_0\bm{x}_2+\bm{v}_{20}^2\bm{D}_0\bm{x}_3\}
\\
\nabla_{\bm{x}_2}|\bm{d}_0| &=
\frac{1}{12|\tau_n|R_n}\{(-2c_{31}+k_{31}\bm{D}_0)\bm{x}_0+\bm{v}_{30}^2\bm{D}_0\bm{x}_1+
2c_{31}\bm{x}_2-\bm{v}_{10}^2\bm{D}_0\bm{x}_3\} \\
\nabla_{\bm{x}_3}|\bm{d}_0| &= 
\frac{1}{12|\tau_n|R_n}\{(-2c_{12}+k_{12}\bm{D}_0)\bm{x}_0-\bm{v}_{20}^2\bm{D}_0\bm{x}_1+
\bm{v}_{10}^2\bm{D}_0\bm{x}_2+2c_{12}\bm{x}_3\}
\end{aligned}
$$
Then we can write $\nabla_{\bm{V}_n}|\bm{d}_0|$ in matrix form
$$
\nabla_{\bm{V}_n}|\bm{d}_0| = 
\frac{1}{12|\tau_n|R_n}
\begin{bmatrix}
\bm{L} & -d_2\bm{K} & d_1\bm{K} \\
d_2\bm{K} & \bm{L} & -d_0\bm{K} \\
-d_1\bm{K} & d_0\bm{K} & \bm{L}
\end{bmatrix} \bm{V}_n
$$
with
$$
\bm{L}=
\begin{bmatrix}
2c & -2c_{23} & -2c_{31} & -2c_{12} \\
-2c_{23} & 2c_{23} & 0 & 0 \\
-2c_{31} & 0 & 2c_{31} & 0 \\
-2c_{12} & 0 & 0 & 2c_{12}
\end{bmatrix}, 
\bm{K}=
\begin{bmatrix}
0 & -k_{23} & -k_{31} & -k_{12} \\
k_{23} & 0 & -\bm{v}_{30}^2 & \bm{v}_{20}^2 \\
k_{31} & \bm{v}_{30}^2 & 0 & -\bm{v}_{10}^2 \\
k_{12} & -\bm{v}_{20}^2 & \bm{v}_{10}^2 & 0
\end{bmatrix}
$$
$c=c_{23}+c_{31}+c_{12}$, matrix $\bm{L}$ is symmetric, $\bm{K}$ is antisymmetric.

For $\nabla_{\bm{x}_i}s$, we denote
$$
\begin{aligned}
p_0 &= \frac{\bm{v}_{31}^2}{4s_2}+\frac{\bm{v}_{21}^2}{4s_3}+\frac{\bm{v}_{32}^2}{4s_1},\quad 
p_1 = \frac{\bm{v}_{32}^2}{4s_0}+\frac{\bm{v}_{02}^2}{4s_3}+\frac{\bm{v}_{30}^2}{4s_2},\quad\\ 
p_2 &= \frac{\bm{v}_{30}^2}{4s_1}+\frac{\bm{v}_{10}^2}{4s_3}+\frac{\bm{v}_{31}^2}{4s_0},\quad
p_3 = \frac{\bm{v}_{10}^2}{4s_2}+\frac{\bm{v}_{20}^2}{4s_1}+\frac{\bm{v}_{21}^2}{4s_0},\quad
\end{aligned}
$$

$$
\begin{aligned}
q_{01} &=-(\frac{\bm{v}_{31}\cdot\bm{v}_{30}}{4s_2}+\frac{\bm{v}_{21}\cdot\bm{v}_{20}}{4s_3}),
q_{02} =-(\frac{\bm{v}_{32}\cdot\bm{v}_{30}}{4s_1}+\frac{\bm{v}_{12}\cdot\bm{v}_{10}}{4s_3}),
q_{03} =-(\frac{\bm{v}_{23}\cdot\bm{v}_{20}}{4s_1}+\frac{\bm{v}_{13}\cdot\bm{v}_{10}}{4s_2}),
\\
q_{12} &=-(\frac{\bm{v}_{32}\cdot\bm{v}_{31}}{4s_0}+\frac{\bm{v}_{02}\cdot\bm{v}_{01}}{4s_3}),
q_{13} =-(\frac{\bm{v}_{03}\cdot\bm{v}_{01}}{4s_2}+\frac{\bm{v}_{23}\cdot\bm{v}_{21}}{4s_0}),
q_{23} =-(\frac{\bm{v}_{13}\cdot\bm{v}_{12}}{4s_0}+\frac{\bm{v}_{03}\cdot\bm{v}_{02}}{4s_1})
\end{aligned}
$$
We have
$$
\begin{aligned}
\nabla_{\bm{x_0}}s
&=p_0\bm{x}_0+q_{01}\bm{x}_1+q_{02}\bm{x}_{2}+q_{03}\bm{x}_{3},
\nabla_{\bm{x}_1}s
=q_{01}\bm{x}_{0}+p_1\bm{x}_1+q_{12}\bm{x}_2+q_{23}\bm{x}_{3} \\
\nabla_{\bm{x}_2}s
&=q_{02}\bm{x}_0+q_{12}\bm{x}_{1}+p_2\bm{x}_2+q_{23}\bm{x}_{3},
\nabla_{\bm{x}_3}s
=q_{03}\bm{x}_0+q_{13}\bm{x}_{1}+q_{23}\bm{x}_{2}+p_3\bm{x}_3
\end{aligned}
$$
Then we can write $\nabla_{\bm{V}_n}s$ in matrix form
$$
\nabla_{\bm{V}_n}\mu_n = 
\begin{bmatrix}
\bm{S} &  &  \\
 & \bm{S} &  \\
 &  & \bm{S}
\end{bmatrix} \bm{V}_n
$$
with 
$$
\bm{S} = 
\begin{bmatrix}
p_0 & q_{01} & q_{02} & q_{03}\\
q_{01} & p_1 & q_{12} & q_{13}\\
q_{02} & q_{12} & p_2 & q_{23}\\
q_{03} & q_{13} & q_{23} & p_3
\end{bmatrix}
$$
Matrix $\bm{S}$ is symmetric.

For $\nabla_{\bm{x}_i}|\tau_n|$, we have
$$
\begin{aligned}
\nabla_{\bm{x}_0}|\tau_n| 
&= \frac{1}{6}(\bm{x}_2\times\bm{x}_1+\bm{x}_3\times\bm{x}_2+\bm{x}_1\times\bm{x}_3),
\nabla_{\bm{x}_1}|\tau_n|
=\frac{1}{6}(\bm{x}_3\times\bm{x}_0+\bm{x}_0\times\bm{x}_2+\bm{x}_2\times\bm{x}_3),
\\
\nabla_{\bm{x}_2}|\tau_n| 
&=\frac{1}{6}(\bm{x}_1\times\bm{x}_0+\bm{x}_3\times\bm{x}_1+\bm{x}_0\times\bm{x}_3),
\nabla_{\bm{x}_3}|\tau_n|
=\frac{1}{6}(\bm{x}_2\times\bm{x}_0+\bm{x}_0\times\bm{x}_1+\bm{x}_1\times\bm{x}_2)
\end{aligned}
$$
We denote
$$
\bm{C}_0=
\begin{bmatrix}
0 & x_2 & x_3 & x_1 \\
x_3 & 0 & x_0 & x_2 \\
x_1 & x_3 & 0 & x_0 \\
x_2 & x_0 & x_1 & 0
\end{bmatrix},
\bm{C}_1=
\begin{bmatrix}
0 & y_2 & y_3 & y_1 \\
y_3 & 0 & y_0 & y_2 \\
y_1 & y_3 & 0 & y_0 \\
y_2 & y_0 & y_1 & 0
\end{bmatrix},
\bm{C}_2=
\begin{bmatrix}
0 & z_2 & z_3 & z_1 \\
z_3 & 0 & z_0 & z_2 \\
z_1 & z_3 & 0 & z_i \\
z_2 & z_0 & z_1 & 0
\end{bmatrix}
$$
Then we have
$$
\nabla_{\bm{V}_n}|\tau_n|=
\frac{1}{6}
\begin{bmatrix}
    0 & -\bm{C}_2 & \bm{C}_1 \\
    \bm{C}_2 & 0 & -\bm{C}_0 \\
    -\bm{C}_1 & \bm{C}_0 & 0
\end{bmatrix}\bm{V}_n = \bm{C}\bm{V}_n
$$
Where $|\tau_n|$ is a cubic form with respect to $\bm{V}_n$, we have the 
following lemma for cubic forms.
\begin{lemma}
If the function $f$ is a cubic form with respect to $\bm{X}=(x_0,x_1,\cdots,x_n)$, 
then there exist matrix $\bm{C}$ and $\bm{E}=\frac{1}{2}(\bm{C}+\bm{C}^T)$ such 
that
$$
\nabla f = \bm{C}\bm{X}=\bm{E}\bm{X}
$$
\begin{proof}
If the function $f$ is a cubic form with respect to $\bm{X}=(x_0,x_1,\cdots,x_n)$, 
then there exists a third-order matrix $\bm{U}$ such that
$$
f = U_{ijk}x_ix_jx_k
$$
Let $W_{ik}^j=\frac{1}{2}(U_{ijk}+U_{kji})$, for any fixed $j$, $\bm{W}^j$ 
is a symmetric matrix, and we have
$$
f = W_{ik}^j x_i x_jx_k
$$
So we have $\nabla_{x_k}f=U_{ijk}x_ix_j=W_{ik}^jx_ix_j$, let
$C_{ki} = U_{ijk}x_j, E_{ki} = W_{ik}^jx_j$, so
$\bm{E} = \frac{1}{2}(\bm{C}+\bm{C}^T)$, then we can express 
$\nabla f$ as
$$
\nabla f=\bm{C}\bm{X}=\bm{E}\bm{X}
$$
\end{proof}
\end{lemma}
From the above theorem, we can get
$\nabla_{\bm{V}_n}|\tau_n|=\frac{1}{2}(\bm{C}+\bm{C}^T)\bm{V}_n$.
$$
\frac{1}{2}(\bm{C}+\bm{C}^T)=
\frac{1}{12}
\begin{bmatrix}
0 & \bm{C}_2^T-\bm{C}_2 & \bm{C}_1^T-\bm{C}_1 \\
\bm{C}_2-\bm{C}_2^T & 0 & \bm{C}_0^T-\bm{C}_0 \\
\bm{C}_1^T-\bm{C}_1 & \bm{C}_0^T-\bm{C}_0 & 0
\end{bmatrix}
$$
It is evident that each block matrix in $\frac{1}{2}(\bm{C}+\bm{C}^T)$ is an 
antisymmetric matrix.

From the equation
$$
\nabla_{\bm x_0} \mu_n = \mu_n\left\{\frac{1}{|\bm{d}_0|}\nabla_{\bm
x_0}|\bm{d}_0|+\frac{1}{s}\nabla_{\bm x_0}s-\frac{2}{|\tau_n|}\nabla_{\bm
x_0}|\tau_n|\right\}
$$
and the gradients $\nabla_{\bm{V}_n}|\bm{d}_0|$, $\nabla_{\bm{V}_n}s$, and 
$\nabla_{\bm{V}_n}|\tau_n|$, we can obtain
$$
\nabla_{\bm{V}_n}\mu_n=
\mu_n\begin{bmatrix}
    \bm{A}_n & \bm{B}_{0n} & \bm{B}_{1n} \\
    -\bm{B}_{2n} & \bm{A}_n & \bm{B}_{0n} \\
    -\bm{B}_{1n} & -\bm{B}_{0n} & \bm{A}_n
\end{bmatrix}
\bm{V}_n
$$
with
$$
\bm{A}_n=\frac{1}{12|\tau_n|R_n}\bm{L}+\frac{1}{s}\bm{S}, 
\bm{B}_{0n}=\frac{-d_0}{12|\tau_n|R_n}\bm{K}-\frac{1}{6|\tau_n|}(\bm{C}_0^T-\bm{C}_0),
$$
$$
\bm{B}_{1n}=\frac{d_1}{12|\tau_n|R_n}\bm{K}-\frac{1}{6|\tau_n|}(\bm{C}_1^T-\bm{C}_1),
\bm{B}_{2n}=\frac{-d_2}{12|\tau_n|R_n}\bm{K}-\frac{1}{6|\tau_n|}(\bm{C}_2^T-\bm{C}_2)
$$
Obviously, $\bm{A}_n$ is symmetric matrix, and $\bm{B}_{0n}$, $\bm{B}_{1n}$, 
$\bm{B}_{2n}$ are antisymmetric matrices.

Then, we can assemble the gradient matrix of each element into a global 
gradient matrix. Therefore, we get the following theorem
\begin{theorem}
\label{thm:3dgrad}
For the energy function
$$
F=\frac{1}{N_c}\sum_{n=0}^{N_c-1}\mu_n
$$
of the tetrahedral mesh,we have
$$
\nabla_{\bm{V}}F=
\begin{bmatrix}
\bm{A} & \bm{B}_2 & \bm{B}_1 \\
-\bm{B}_2 & \bm{A} & \bm{B}_0 \\
-\bm{B}_1 & -\bm{B}_0 & \bm{A}
\end{bmatrix}
\begin{bmatrix}
\bm{X} \\
\bm{Y} \\
\bm{Z} 
\end{bmatrix}
=\bm{G}_F\bm{V} 
$$
where $\bm{A}$ is $N_v\times N_v$ symmetric matrix, and 
$\bm{B}_0$, $\bm{B}_1$, $\bm{B}_2$ are $N_v\times N_v$ antisymmetric matrices.
$$
\bm{V} = [\bm{X}^T,\bm{Y}^T,\bm{Z}^T]^T,
\bm{X}=[x_0,x_1,\cdots,x_{N_v-1}]^T,\bm{Y}=[y_0,y_1,\cdots,y_{N_v-1}]^T,
\bm{Z}=[z_0,z_1,\cdots,z_{N_v-1}]^T
$$.
\end{theorem}

\subsection{Anisotropic case}
The radius ratio energy function can be extended to anisotropic mesh. For a 
triangle element $\tau_n = (\bm x_0, \bm x_1, \bm x_2)$ in a two-dimensional 
triangulation $\mathcal{T}$, the radius ratio is  
$$
\mu_n = \frac{R_n}{2r_n} =\frac{pq}{16|\tau_n|^2}
$$
Let $\bm{M}(\bm{x})$ be the anisotropic metric at point $\bm{x}$. To use a 
consistent metric within each element, we define the element metric as the 
average of the vertex metrics:
$$
\bm{M}=\frac{\bm{M}(\bm{x}_0)+\bm{M}(\bm{x}_1)+\bm{M}(\bm{x}_2)}{3}=\begin{bmatrix}
m_0& m_2 \\
m_2 & m_1
\end{bmatrix}
$$
Under this metric, the lengths of the three edges of the triangle can be 
written as
$$
l_0 = \sqrt{<\bm{u}_1-\bm{u}_2,\bm{M}(\bm{u}_1-\bm{u}_2)>},
l_1 = \sqrt{<\bm{u}_0-\bm{u}_2,\bm{M}(\bm{u}_0-\bm{u}_2)>},
l_2 = \sqrt{<\bm{u}_0-\bm{u}_1,\bm{M}(\bm{u}_0-\bm{u}_1)>}
$$
The area of the element under the anisotropic metric satisfies
$$
|\tau_n| = \sqrt{\det{\bm{M}}}|\tau_n|_{\Omega}
$$
where $|\tau_n|_{\Omega}$ is the area under the Euclidean metric.

Based on these definitions, the gradients of the quantities in the radius ratio
with respect to the vertex positions can be derived. The derivation is similar
to the two-dimensional isotropic case. Finally, we obtain
$$
\nabla_{\bm{V}_n}\mu_n = 
\begin{bmatrix}
m_0\bm{A}_n & m_2\bm{A}_n+\bm{B}_n\\
m_2\bm{A}_n+\bm{B}^T_n & m_1\bm{A}_n
\end{bmatrix} \bm{V}_n
$$
with
$$
\bm{A}_n= 
\mu_n
\begin{pmatrix}
c_{12} & -c_{2} & -c_{1} \\
-c_{2} & c_{20} & -c_{0} \\
-c_{1} & -c_{0} & c_{01}
\end{pmatrix},
\bm{B}_n = \mu_n
\begin{pmatrix}
0 & -c & c \\
c & 0 & -c \\
-c & c & 0
\end{pmatrix}
$$
$$
\bm{V}_n=
\begin{bmatrix}
 x_0,x_1,x_2,y_0,y_1,y_2
\end{bmatrix}^T
$$
$$
c=\frac{\sqrt{\det\bm{M}}}{|\tau_n|}, 
c_0 = \frac{1}{pl_0}+\frac{1}{l_0^2}, 
c_1 = \frac{1}{pl_1}+\frac{1}{l_1^2}, 
c_2=\frac{1}{pl_2}+\frac{1}{l_2^2}
$$
$$
c_{01}=c_0+c_1, c_{02}=c_0+c_2, c_{12}=c_1+c_2
$$

For a tetrahedral element $\tau_n = (\bm x_0, \bm x_1, \bm x_2,\bm{x}_3)$ in a 
three-dimensional tetrahedralization $\mathcal{T}$, the radius ratio is  
$$
\mu_n =\frac{R_n}{3r_n}=\frac{s|\bm d_0|}{108|\tau|^2}
$$
We define the anisotropic metric of the element as the average of the vertex 
metrics:
$$
\bm{M}=\frac{\bm{M}(\bm{x}_0)+\bm{M}(\bm{x}_1)+\bm{M}(\bm{u}_2)+\bm{M}(\bm{x}_3)}{4}= 
\begin{bmatrix}
 m_{00} & m_{01} & m_{02} \\
m_{01} & m_{11} & m_{12} \\
m_{02} & m_{12} & m_{22}
\end{bmatrix}
$$
Based on this metric, the gradient of the radius ratio of the anisotropic
tetrahedral element with respect to the vertex positions can be derived. It can
be written in a matrix–vector form. Let
$$
\bm{V}_n=
\begin{bmatrix}
x_0, x_1, x_2, x_3, y_0, y_1, y_2, y_3,
z_0, z_1, z_2, z_3 
\end{bmatrix}^T
$$
then the element gradient can be expressed as
$$
\nabla_{\bm{V}_n}\mu_n=
\mu_n\begin{bmatrix}
\bm{A}_0 & \bm{B}_{2} & \bm{B}_{1} \\
\bm{B}_2^T & \bm{A}_1 & \bm{B}_{0} \\
\bm{B}_1^T & \bm{B}_0^T & \bm{A}_2
\end{bmatrix}
\bm{V}_n
$$
where $\bm{A}_i,\bm{B}_i~(i=0,1,2)$ are both symmetric matrices. Since the 
derivation is similar to the isotropic case, it is omitted here. 

Under an anisotropic metric, both two-dimensional triangular elements and
three-dimensional tetrahedral elements admit a radius ratio $\mu$ and its
gradient with respect to the vertex coordinates.Consider a simplex mesh
$\mathcal{T}$ with $N_v$ vertices
$\{\bm{x}_i\}_{i=0}^{N_v-1}\subset\mathbb{R}^3$ and  $N_c$ elements
$\{\tau_n\}_{n=0}^{N_c-1}$. We define the global anisotropic radius ratio
energy as
$$
F=\frac{1}{N_c}\sum_{n=0}^{N_c-1}\mu_n
$$
This energy sums the anisotropic radius ratios of all elements and reflects the 
overall shape quality of the mesh under the given anisotropic metric.

Similar to the isotropic case, the gradient of the radius ratio for each 
element can be written in a matrix–vector form. Therefore, the gradient of the 
global anisotropic energy can be obtained by assembling the element gradient 
matrices into a global matrix, in the same way as in Theorem \ref{thm:2dgrad} 
and Theorem \ref{thm:3dgrad}, and then computing the gradient of the energy function.

\section{Preconditioner design}
In this work, one key step in mesh optimization is to compute the minimizer of 
the energy function and obtain updated vertex positions. Since the gradient of 
the energy function can be computed, we can apply first-order numerical 
optimization methods, such as gradient descent, quasi-Newton methods, and
nonlinear conjugate gradient (NLCG) methods. From the derivation in
Section~\ref{sec:rro}, the gradient of the energy can be written in a 
matrix-vector form. The matrix contains partial Hessian information. Based on 
this structure, we design a preconditioning operator. This idea is mainly 
inspired by fast algorithms developed for Optimal Delaunay Triangulation (ODT) 
~\cite{10ref}and Centroidal Voronoi Tessellation (CVT)~\cite{11ref}.

From Theorem~\ref{thm:2dgrad}, the gradient of the radius ratio energy for a
two-dimensional simplex mesh can be written as
$$
\nabla_{\bm{V}}F=
\begin{bmatrix}
\bm{A} & \bm{B} \\
-\bm{B} & \bm{A} \\
\end{bmatrix}
\begin{bmatrix}
\bm{X} \\
\bm{Y}
\end{bmatrix}
=\bm{G}_F\bm{V}
$$
Based on this, the Hessian matrix can be written in the form
$$
\bm{H}=\nabla_{\bm{V}}^2F=\bm{G}_F+\frac{\partial\bm{G}_F}{\partial\bm{V}}\bm{V}
$$
The term $\frac{\partial\bm{G}_F}{\partial\bm{V}}\bm{V}$ describes how 
$\bm{G}_F$ changes with the vertex coordinates $\bm{V}$. This term is 
difficult to compute. In practice, we treat $\bm{G}_F$ as a constant matrix at 
the current iteration. Then $\bm{G}_F$ is used as the main approximation of the 
Hessian.

In $\bm{G}_F$, the matrix $\bm{B}$ is antisymmetric and does not contribute to 
the quadratic form. The quadratic form is fully determined by $\bm{A}$. 
Therefore, when constructing the preconditioner, we can ignore $\bm{B}$ and 
keep only the symmetric dominant part. The preconditioner is chosen as
$$
\bm{P}=\begin{bmatrix}
\bm{A} & \\
& \bm{A}
\end{bmatrix}
$$
From Theorem~\ref{thm:3dgrad}, the gradient of the radius ratio energy for a 
three-dimensional simplex mesh can be written as
$$
\nabla_{\bm{V}}F=
\begin{bmatrix}
\bm{A} & \bm{B}_2 & \bm{B}_1 \\
-\bm{B}_2 & \bm{A} & \bm{B}_0 \\
-\bm{B}_1 & -\bm{B}_0 & \bm{A}
\end{bmatrix}
\begin{bmatrix}
\bm{X} \\
\bm{Y} \\
\bm{Z}
\end{bmatrix}=
\bm{G}_F\bm{V}
$$
In the two-dimensional case, diagonal blocks can be directly used as the
preconditioner. In three dimensions, the block matrix $\bm{A}$ has a more
complex structure. It is assembled from the element-wise matrices $\bm{A}_n$,
where
$$
\bm{A}_n = \frac{1}{12|\tau_n|R_n}\bm{L}+\frac{1}{s}\bm{S}
$$
$\bm{S}$ is Laplacian matrix, and the entries of $\bm{L}$ are given by
$$
c_{ij}=\bm{d}_0\cdot(\bm{v}_{i0}\times\bm{v}_{j0})
$$
which are related to face normals. When the element quality is poor, $c_{ij}$ 
may become negative.This makes $\bm{L}$ indefinite and affects the stability of 
$\bm{A}_n$.

To address this issue, we modify $\bm{L}$ by taking absolute values:
$$
\bm{L}_{\text{abs}}=\begin{bmatrix}
2(|c_{23}|+|c_{31}|+c_{12}|) & -2|c_{23}| & -2|c_{31}| & -2|c_{12}| \\
-2|c_{23}| & 2|c_{23}| & 0 & 0 \\
-2|c_{31}| & 0 & 2|c_{31}| & 0 \\
-2|c_{12}| & 0 & 0 & 2|c_{12}|
\end{bmatrix}
$$
and define
$$
\bm{A}_{n,\text{abs}}=\frac{1}{12|\tau_n|R_n}\bm{L}_{\text{abs}}+\frac{1}{s}\bm{S}
$$
After assembling $\bm{A}_{n,abs}$ over all elements, we obtain a global matrix 
$\bm{A}_{\text{abs}}$, which is Laplacian. Based on this, we 
choose the following block diagonal form as the preconditioner:
$$
\bm{P}=\begin{bmatrix}
\bm{A}_{abs} & & \\
& \bm{A}_{abs} & \\
& & \bm{A}_{abs}
\end{bmatrix}
$$
In both two and three dimensions, the preconditioner $\bm{P}$ is a sparse and 
symmetric positive semi-definite Laplacian matrix. To make it symmetric 
positive definite in practice, we can fix the boundary nodes. If the boundary 
nodes are not fixed, we can apply a projection operator to restrict their 
motion to the geometric constraint manifold.

For each mesh node, we define a local $3\times3$ projection matrix $\bm{T}_i$: 
$$
\begin{cases}
\bm{T}_i = \bm{I}_3\qquad&\text{if} ~\bm{x}_i~\text{is an interior node} \\
\bm{T}_i = \bm{I}_3-\bm{n}_i\bm{n}_i^T\qquad&\text{if}~\bm{x}_i~\text{is on a surface boundary} \\
\bm{T}_i = \bm{t}_i\bm{t}_i^T \qquad&\text{if}~\bm{x}_i~\text{is on a curve boundary}
\end{cases}
$$
where $\bm{n}_i$ is the unit normal vector at a surface boundary node, and 
$\bm{t}_i$ is the unit tangent vector at a curve boundary node. These matrices 
can be assembled into a global projection matrix $\bm\Pi$.

By either fixing boundary nodes or applying the projection matrix  $\bm\Pi$, 
the preconditioner becomes symmetric positive definite. The resulting linear 
system can then be solved efficiently by the conjugate gradient
method~\cite{30ref} or algebraic multigrid methodcite~\cite{28ref,29ref}, which 
makes the approach suitable for large-scale mesh optimization problems.

\section{Energy Optimization}\label{sec:energyopt}
We call our mesh optimization algorithm the Radius Ratio Energy (RRE) mesh 
optimization algorithm. The RRE mesh optimization algorithm consists of two main 
parts:
\begin{enumerate}
\item Mesh vertex relocation based on the radius ratio energy
\item Connectivity update
\end{enumerate}
In the vertex relocation step, we define an energy function based on the radius 
ratio of mesh elements and compute its gradient with respect to the vertex 
coordinates. Then, we apply a first-order optimization method to minimize the 
energy. This gives new vertex positions and updates the mesh.

The BFGS method is one of the most widely used quasi-Newton methods. In this 
work, we mainly use its preconditioned limited-memory version (PLBFGS). Below 
we present one iteration of the PLBFGS method. More details can be found 
in~\cite{22ref,24ref}.
\begin{algorithm}[!h]
\caption{PLBFGS Algorithm}
\begin{algorithmic}[1]
    \STATE At current point $\bm{x}_k$, compute the gradient $\bm{g}_k$
    \STATE Set 
    $$\delta=
        \begin{cases}
            0,\quad\text{if}~k\le m \\
            k-m,\quad\text{if}~k>m
        \end{cases};
        L=
        \begin{cases}
        k,\quad\text{if}~k\le m \\
        m,\quad\text{if}~k>m
        \end{cases};
        \bm{q}=\bm{g}_k
    $$
    \STATE Backward loop:
    \FOR{$i=L-1,L-2,\cdots,1,0$:}
    \STATE $j=i+\delta$
    \STATE $\alpha_i=\rho_j\bm{s}_j^T\bm{q}$
    \STATE $\bm{q}_i=\bm{q}-\alpha_i\bm{y}_j$
    \ENDFOR
    \STATE If no preconditioner is used, set $\bm{r}=\bm{H}_k\bm{q}$, if a 
    preconditioner is used, solve $\bm{P}_k\bm{r}=\bm{q}$
    \STATE Forward Loop:
    \FOR{$i=0,1,2,\cdots,L-2,L-1$:}
    \STATE $j=i+\delta$
    \STATE $\beta_j=\rho_j\bm{y}_j^T\bm{r}$
    \STATE $\bm{r} = \bm{r}+(\alpha_i-\beta_i)\bm{s}_j$
    \ENDFOR
    \STATE Set the search direction $\bm{d}_k=-\bm{r}$. Use a line search to 
    find the step size $\alpha_k$, and update 
    $\bm{x}_{k+1}=\bm{x}_k+\alpha_k\bm{d}_k$
    \STATE Compute $\bm{s}_k=\bm{x}_{k+1}-\bm{x}_k$,$\bm{y}_k=\bm{g}_{k+1}-\bm{g}_k$
    and store the most recent $m$ pairs of vectors.
\end{algorithmic}
\end{algorithm}

Here, $\bm{P}_k$ is the chosen preconditioner, and
$$
\bm{H}_k=\frac{\bm{s}_{k-1}^T\bm{y}_{k-1}}{\bm{y}_{k-1}^T\bm{y}_{k-1}}
$$
The step size $\alpha_k$ is obtained by a line search that satisfies the strong 
Wolfe conditions, with initial guess $\alpha_k=1$.

In addition to the PLBFGS method, we also use the preconditioned nonlinear 
conjugate gradient (PNLCG) method~\cite{22ref,23ref}. We first describe one 
iteration from step $k-1$ to step $k$ ($k=1,2,\cdots$) without a preconditioner.
\begin{algorithm}[!h]
\caption{NLCG Algorithm}
\begin{algorithmic}[1]
    \STATE At the current point $\bm{x}_k$, compute the gradient $\bm{g}_k$
    \STATE Compute the scalar coefficient $\beta_k$, and update the search 
    direction $\bm{d}_k=-\bm{g}_k+\beta_k\bm{d}_{k-1}$
    \STATE Use a line search to compute the step size $\alpha_k$, and update
    $\bm{x}_k=\bm{x}_{k-1}+\alpha_k\bm{d}_k$
\end{algorithmic}
\end{algorithm}

Here, $\alpha_k$ is obtained by a line search that satisfies the strong Wolfe 
conditions, with initial guess $\alpha_k=1$. The choice of $\beta_k$ has a 
strong effect on performance~\cite{23ref}. Common choices include the 
Polak–Ribière($\beta_k^{PR}$)~\cite{25ref}, Hestenes–Stiefel($\beta_k^{HS}$)
~\cite{26ref}, Fletcher–Reeves($\beta_k^{FR}$)~\cite{27ref}, and 
Dai-Yuan($\beta_k^{DY}$)~\cite{39ref} formulas. In this work, we use the 
Polak–Ribière formula:
$$
\beta_{k}^{PR} = \frac{\bm{g}_{k}^T(\bm{g}_{k}-\bm{g}_{k-1})}{\bm{g}_{k-1}^T\bm{g}_{k-1}}
$$
Next, we give the iteration with a preconditioner. Let the preconditioned 
gradient $\bm{z}_k$ be defined by solving $\bm{P}_k\bm{z}_k=\bm{g}_k$ and set 
the initial search direction $\bm{d}_0=-\bm{z}_0$.  Then one iteration from 
step $k-1$ to step $k$ is:
\begin{algorithm}[!h]
\caption{PNLCG Algorithm}
\begin{algorithmic}[1]
    \STATE At the current point $\bm{x}_k$, compute the preconditioned gradient
    $\bm{z}_k$.
    \STATE Compute the scalar coefficient $\beta_k$, and update the search 
    direction $\bm{d}_k=-\bm{z}_k+\beta_k\bm{d}_{k-1}$
    \STATE Use a line search to compute $\alpha_k$, and update $\bm{x}_{k+1}=
    \bm{x}_k+\alpha_k\bm{d}_k$
\end{algorithmic}
\end{algorithm}

When a preconditioner is used, $\beta_k$ is also based on the preconditioned 
gradient. The PR-type formula used in this work is
$$
\beta_{k}^{PR} = \frac{\bm{z}_{k}^T(\bm{g}_{k}-\bm{g}_{k-1})}{\bm{z}_{k-1}^T\bm{g}_{k-1}}
$$

The main stopping criterion for vertex relocation is $\|F_{k+1}-F_k\|<\varepsilon$, 
where $\varepsilon$ is usually set to $10^{-6}$. After convergence of the vertex 
relocation step, we update the mesh connectivity. Given an initial connectivity, 
we improve it by flip operations. These operations are applied to tetrahedral
mesh, including 2-3 flips and 3-2 flips~\cite{44ref}. A flip is accepted if it reduces the 
average element energy. We examine all edges and faces that satisfy the flip 
conditions, and repeat until no further energy reduction is possible.

After updating the connectivity, we perform vertex relocation again. These two 
steps are applied in an alternating way. The process stops when no flip 
operation can further reduce the energy after the relocation step.

\section{Experiments}
We implement the proposed algorithm in Python based on the open-source CAX core 
library FEALPy~\cite{20ref}. All experiments are carried out on a laptop running 
Ubuntu 22.04, with an AMD Ryzen 7 7735H CPU (3.2 GHz, 16 MB cache) and 16 GB 
DDR5 RAM (5600 MHz).

Section 5.1 presents examples of two-dimensional triangular mesh, Section 5.2 
presents examples of three-dimensional tetrahedral mesh, and Section 5.3 
presents examples of anisotropic mesh. In the following, mesh quality refers 
to the distribution of the element radius ratio. During optimization, the energy 
function is constructed using the radius ratio $\mu=\frac{R}{dr}$. For 
convenience in statistics and visualization, we report results using 
$\frac{dr}{R}\in(0,1]$. Unless otherwise stated, we use the L-BFGS method and 
its preconditioned version to solve the energy minimization problem in all 
examples.
\subsection{2D TriangleMesh}
Two-dimensional Delaunay triangulations do not produce slivers. Our main 
focus is on tetrahedral mesh. Therefore, for triangular mesh, we only 
perform vertex relocation and use these examples to show the effectiveness of 
the method in improving mesh quality.
\begin{figure}[!h] 
\centering
\subfloat[initial mesh]{
\includegraphics[width=0.3\linewidth]{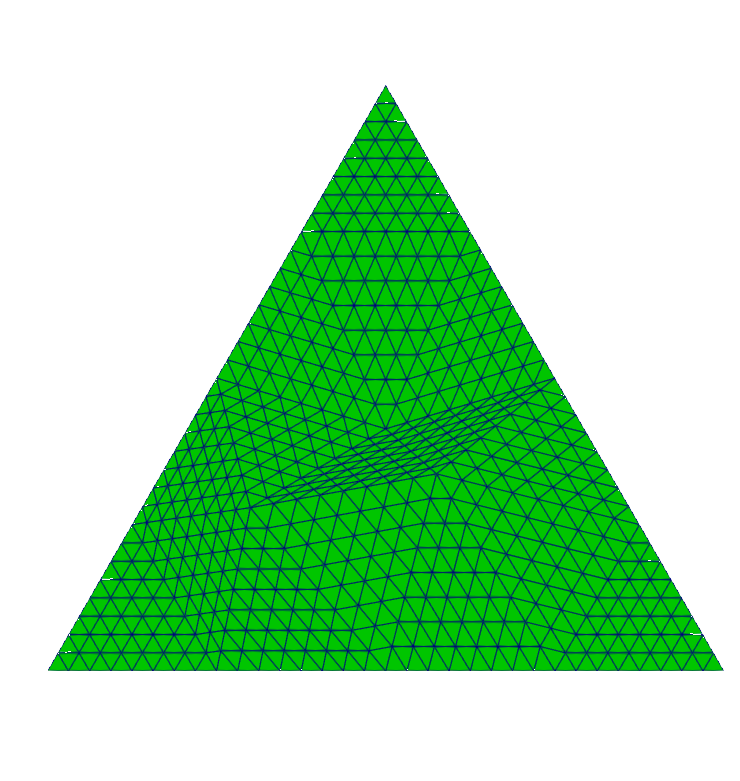}}
\hspace{0.01\linewidth}
\subfloat[ODT local mesh smoothing]{
\includegraphics[width=0.3\linewidth]{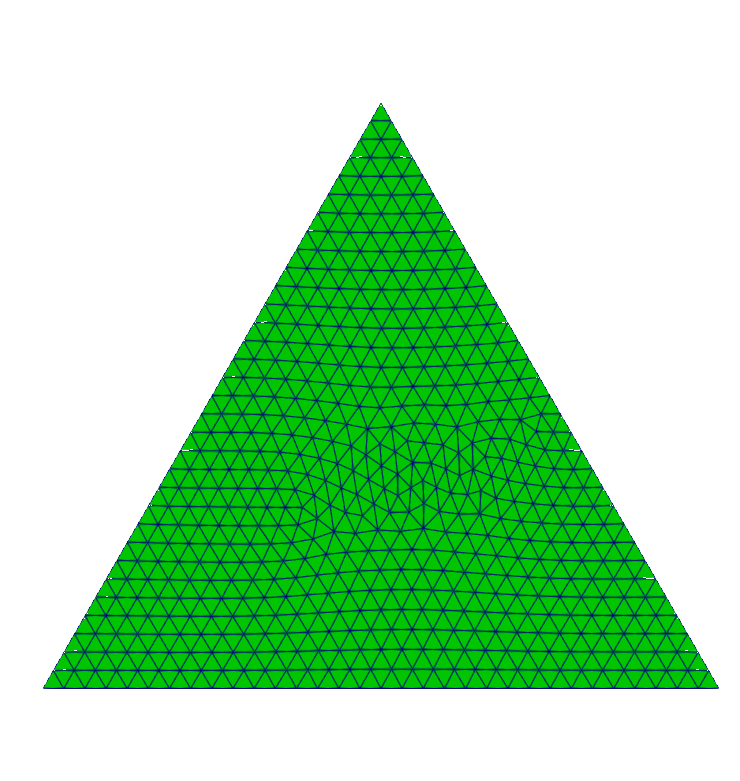}}
\hspace{0.01\linewidth}
\subfloat[RRE optimization]{
\includegraphics[width=0.3\linewidth]{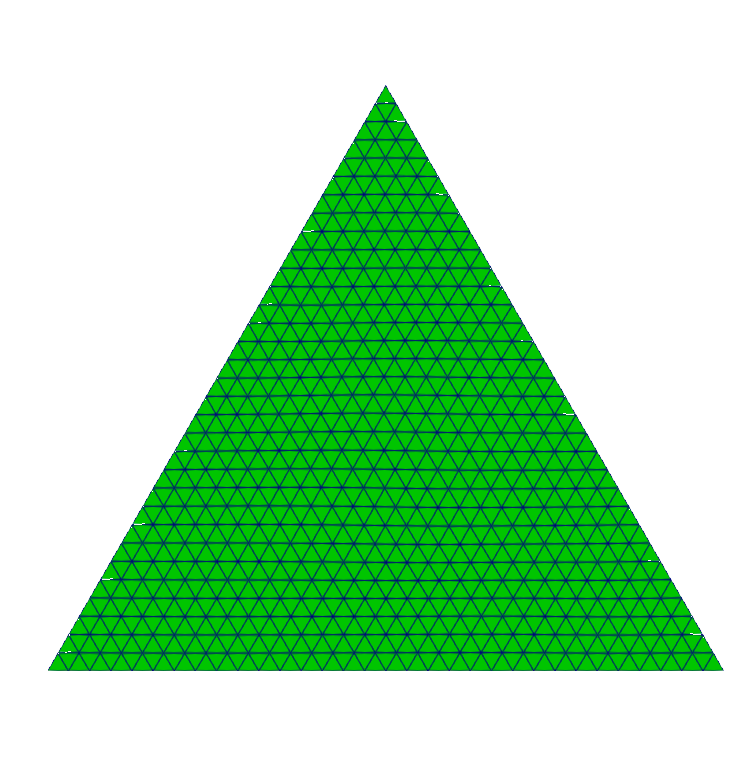}}
\caption{Comparison of mesh obtained by local and global methods}
\label{fig:tri}
\end{figure}

Following the setup in Chen and Holst~\cite{10ref}, we move three interior 
nodes to perturb a uniform mesh of equilateral triangles, and then apply global 
refinement to obtain Fig.~\ref{fig:tri}~(a). We then optimize the mesh using a 
local ODT method and the radius ratio method. Fig.~\ref{fig:tri}~(b) shows the 
result of the local ODT method (a simple Python implementation). This method
relies only on local mesh information and performs smoothing by minimizing a
local energy associated with each vertex or element. Although each iteration
guarantees a decrease of the local energy, this decrease is restricted to a
small neighborhood and does not necessarily lead to a reduction of the global
energy over the entire mesh. As a result, the optimization process may converge
to a configuration that is locally optimal but not globally optimal. Once such a
locally stable state is reached, the algorithm may terminate even though further
global improvement is still possible. As observed in Fig.~\ref{fig:tri}~(b), the 
overall mesh quality is improved compared with the initial mesh, but some 
triangles still deviate from the equilateral shape. Fig.~\ref{fig:tri}~(c) 
shows the result of the radius ratio method. This method defines a global 
energy based on the radius ratio and minimizes it using the L-BFGS method. All 
free nodes are updated together, so global energy decreases over the whole mesh. 
For this example, the mesh topology is already optimal, so the method converges 
to a mesh of equilateral triangles. This example shows that, compared with 
methods based only on local information, a global energy method has a stronger 
ability to improve the mesh and can achieve higher quality.
\begin{figure}[htbp]
\centering
\subfloat[initial mesh]{
\includegraphics[width=0.3\linewidth]{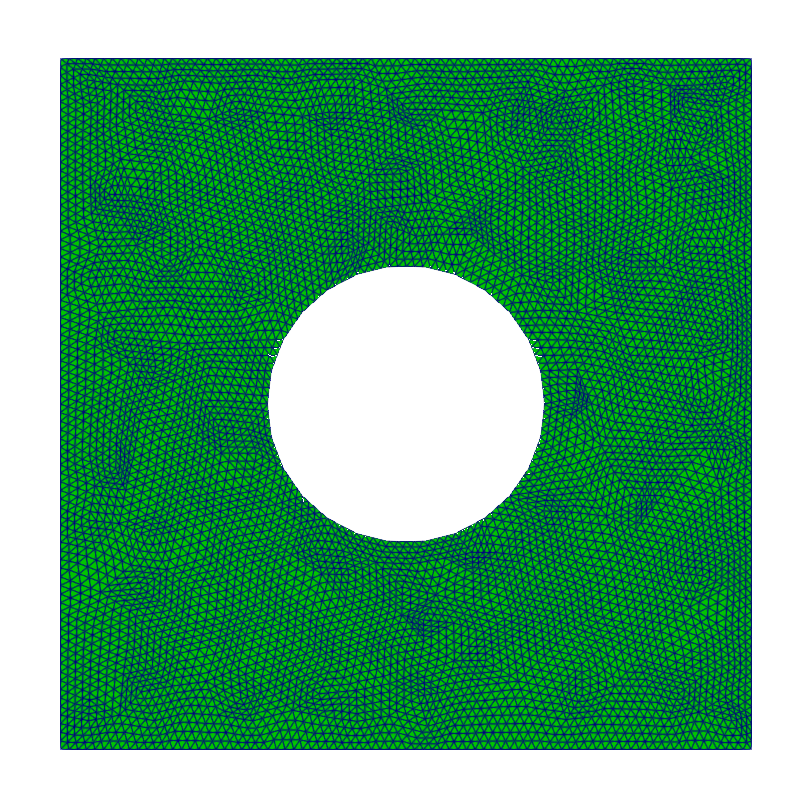}}
\hspace{0.02\linewidth}
\subfloat[RRE]{
\includegraphics[width=0.3\linewidth]{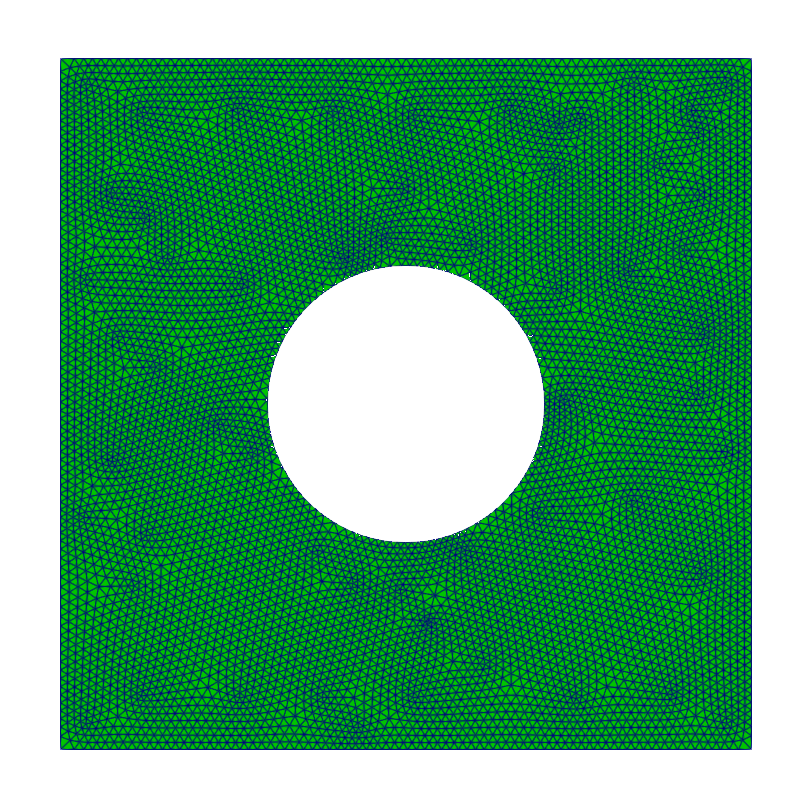}}
\hspace{0.02\linewidth}
\subfloat[Precondition RRE]{
\includegraphics[width=0.3\linewidth]{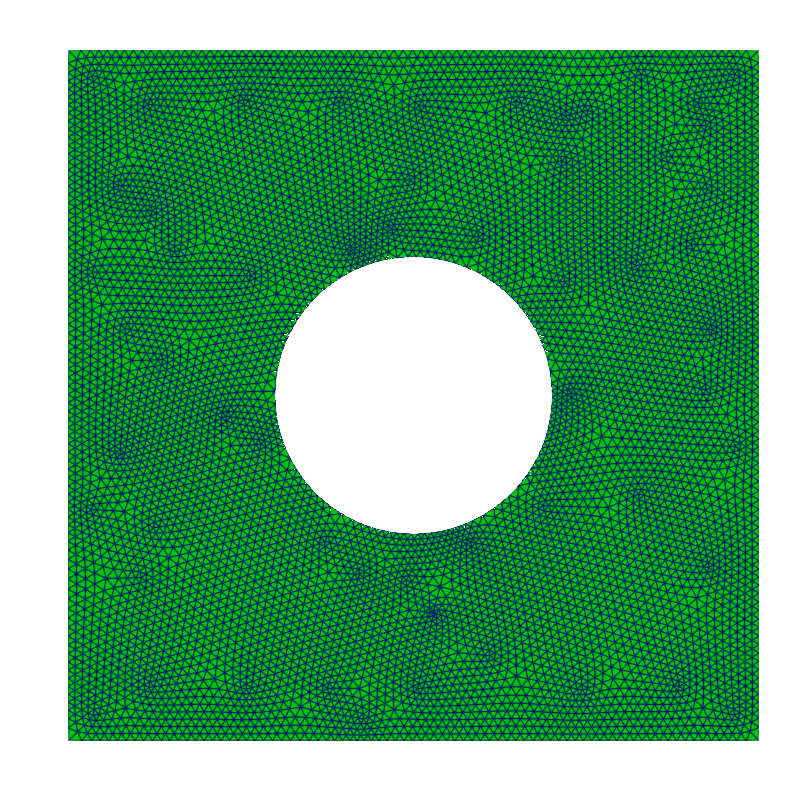}}
\hspace{0.02\linewidth}
\vfill
\subfloat[Quality of initial mesh]{
\includegraphics[width=0.3\linewidth]{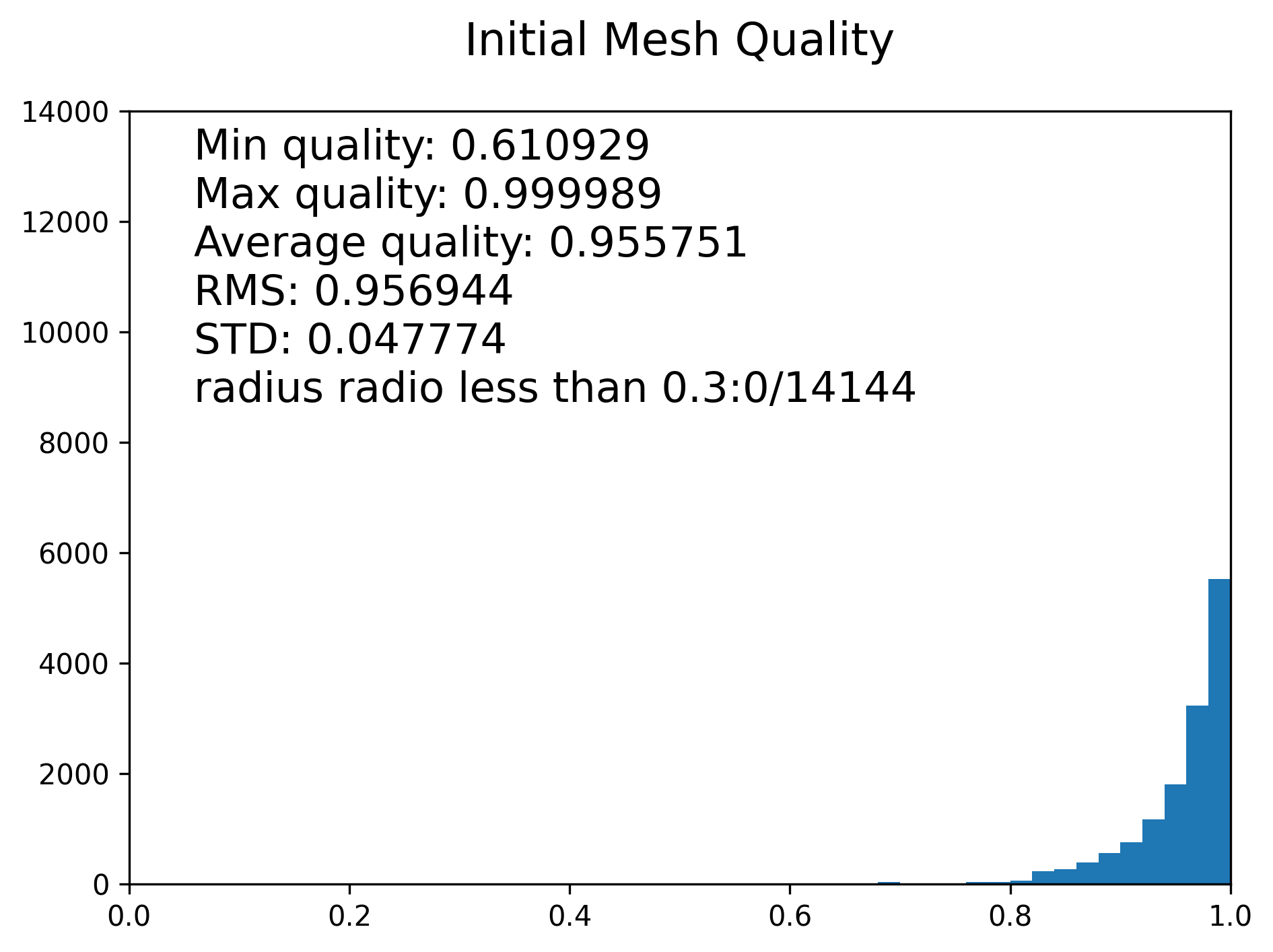}}
\hspace{0.02\linewidth}
\subfloat[Quality of RRE]{
\includegraphics[width=0.3\linewidth]{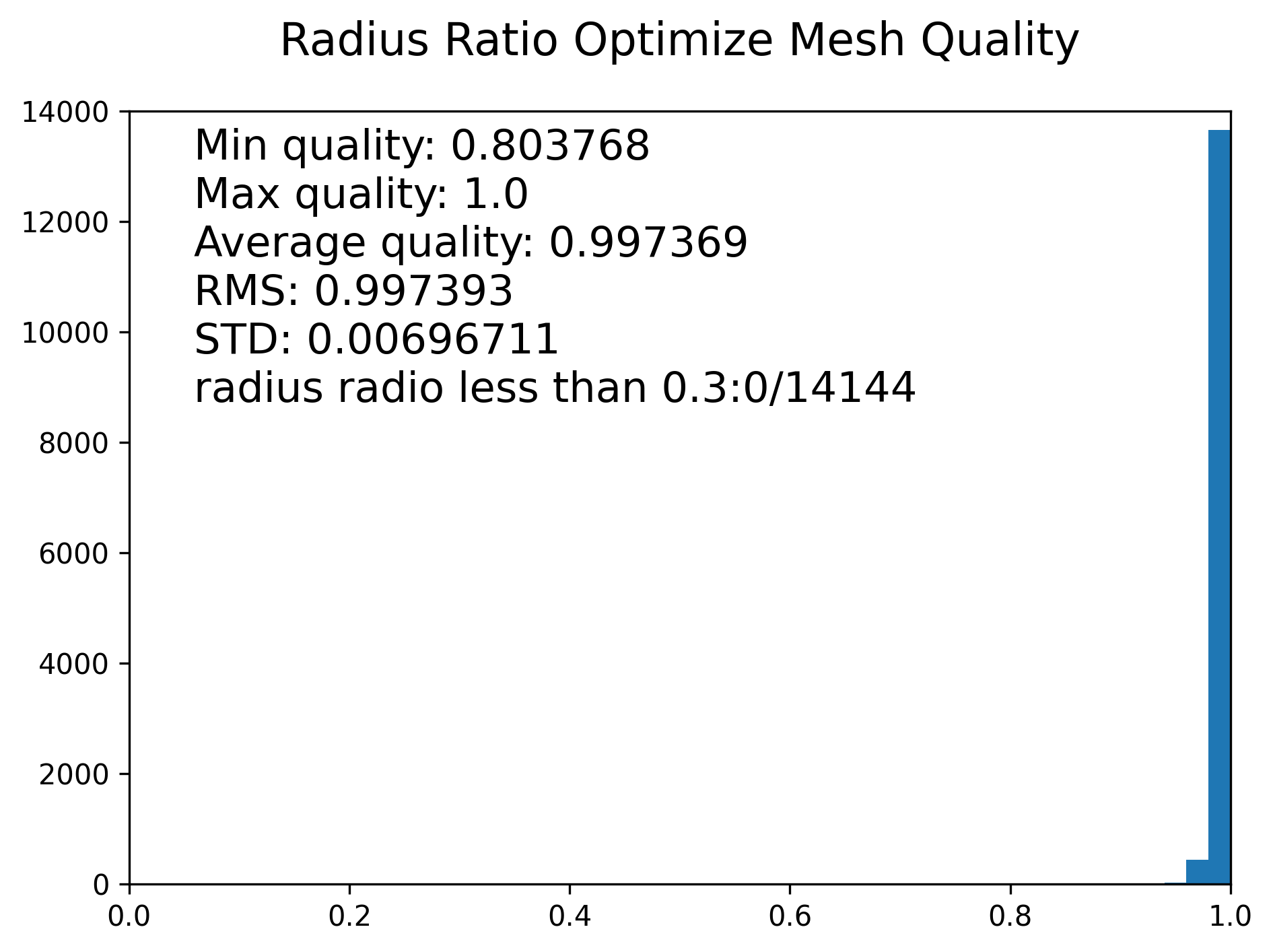}}
\hspace{0.02\linewidth}
\subfloat[Quality of Precondition RRE]{
\includegraphics[width=0.3\linewidth]{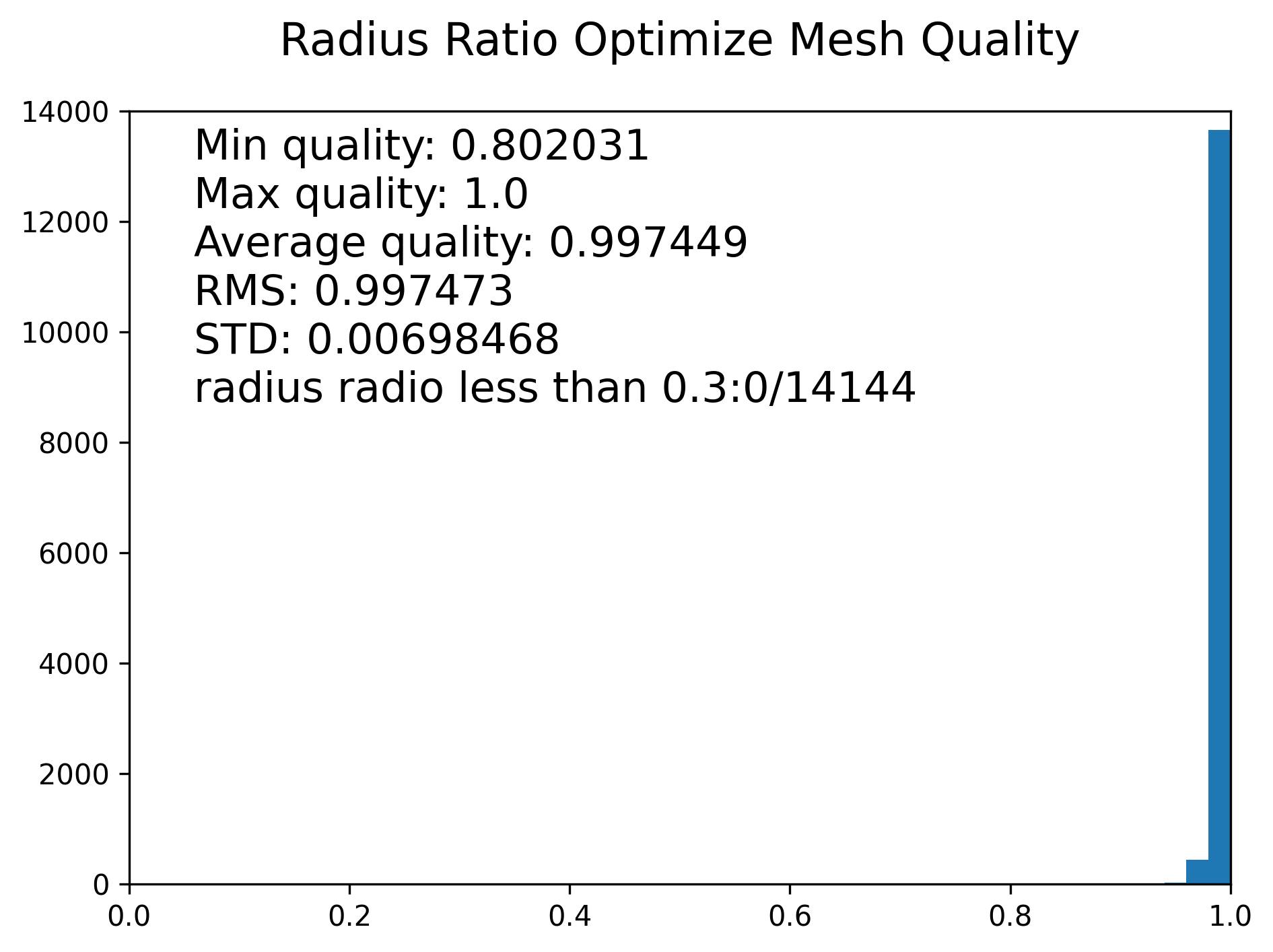}}
\caption{Square domain with a circular hole}
\label{fig:squareh}
\end{figure}

Fig.~\ref{fig:squareh} shows a square domain with a circular hole. The initial 
mesh is generated by Gmsh~\cite{21ref}. To test the ability of the method to improve poor 
elements, we add random perturbations and apply global refinement. As a result, 
the initial mesh contains a certain number of low-quality elements. The results 
in Fig.~\ref{fig:squareh}~(b) and (e) show that the number of low-quality 
elements is greatly reduced, and the minimum radius ratio is clearly improved. 
The optimized mesh is also more uniform, and the element shapes are closer to 
regular ones. Fig.~\ref{fig:squareh}~(c) and (f) show the results with the 
preconditioner. The results are almost the same as those without the 
preconditioner.
\begin{table}[htbp]
    \caption{Optimization results of two-dimensional examples}
    \label{table2d}
    \centering
    %\begin{adjustwidth}{-1.5cm}{0cm}
    \begin{tabular}{|c|c|c|c|c|c|c|}
        \hline
        Model & Method & Min radius ratio & Mean radius ratio &Iteration count
        &Times(sec.) \\
        \hline
        Triangle Domain & Init & 0.490 & 0.896 & / & / \\ 
        \hline
            & RRE & 1.0 & 1.0 & 52 & 0.039    \\
        \hline
        & Precondition RRE & 1.0 & 1.0 & 10 & 0.046 \\
        \hline
        \makecell{Square domain with\\ a circular hole} & Init & 0.611 & 0.956 & / & / \\
        \hline
        & RRE & 0.804 & 0.997 & 193 & 1.436 \\
        \hline
        & Precondition RRE & 0.802 & 0.997 & 72 & 2.119 \\
        \hline
    \end{tabular}
%\end{adjustwidth}
\end{table}

Table~\ref{table2d} lists the results for the two-dimensional examples. In terms of 
iteration counts, the use of the preconditioner reduces the number of iterations. 
Here, the iteration count is not the number of L-BFGS iterations. Instead, it 
is the number of calls to the radius ratio energy function during the Wolfe line 
search. In one L-BFGS iteration, the line search may call the energy function 
several times. Therefore, this count better reflects the efficiency of the 
algorithm. However, since the problem size is small, the cost of building the 
preconditioner is relatively high, and there is no clear advantage in total 
computation time. In the tetrahedral mesh examples, we will present more tests 
to show the benefit of the preconditioner in improving efficiency.
\subsection{3D TetrahedronMesh}
For tetrahedral mesh optimization, the workflow follows Chapter~\ref{sec:energyopt} 
and consists of two steps: vertex relocation and connectivity update.
\begin{figure}[htbp]
\centering
\subfloat[Model]{
\includegraphics[width=0.2\linewidth]{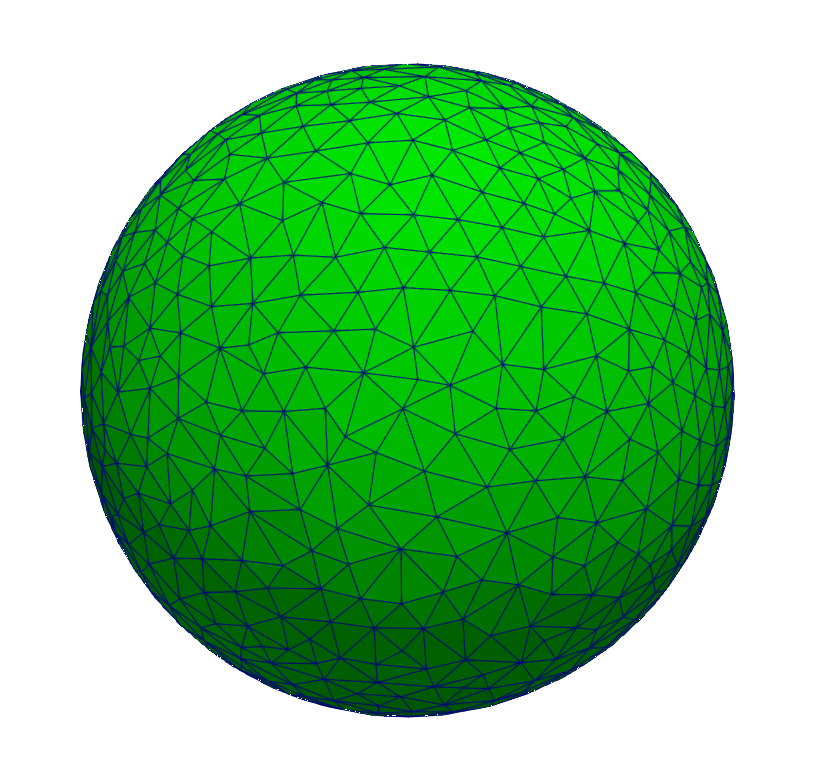}}
\hspace{0.01\linewidth}
\subfloat[Init quality]{
\includegraphics[width=0.2\linewidth]{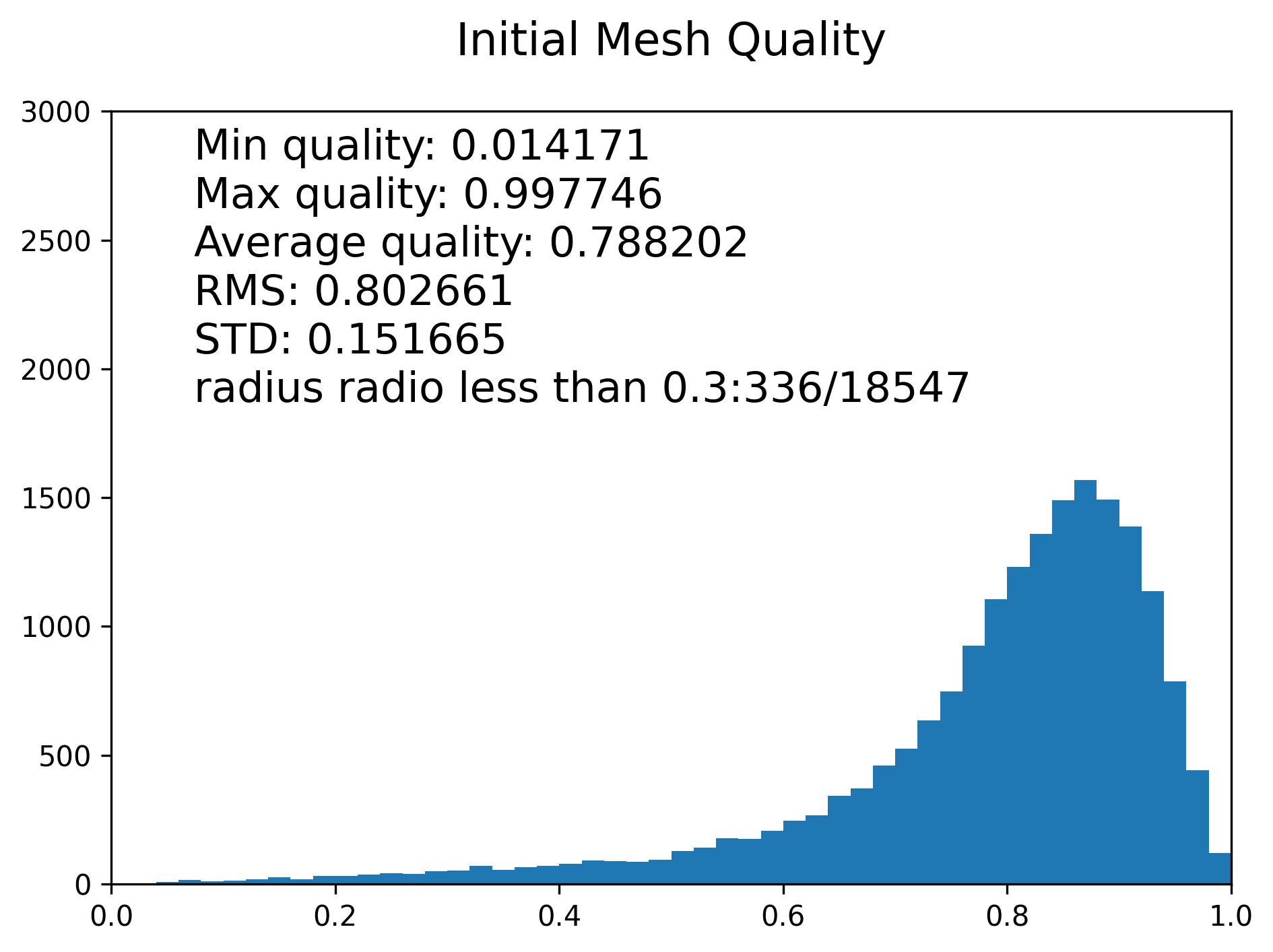}}
\subfloat[RRE]{
\includegraphics[width=0.2\linewidth]{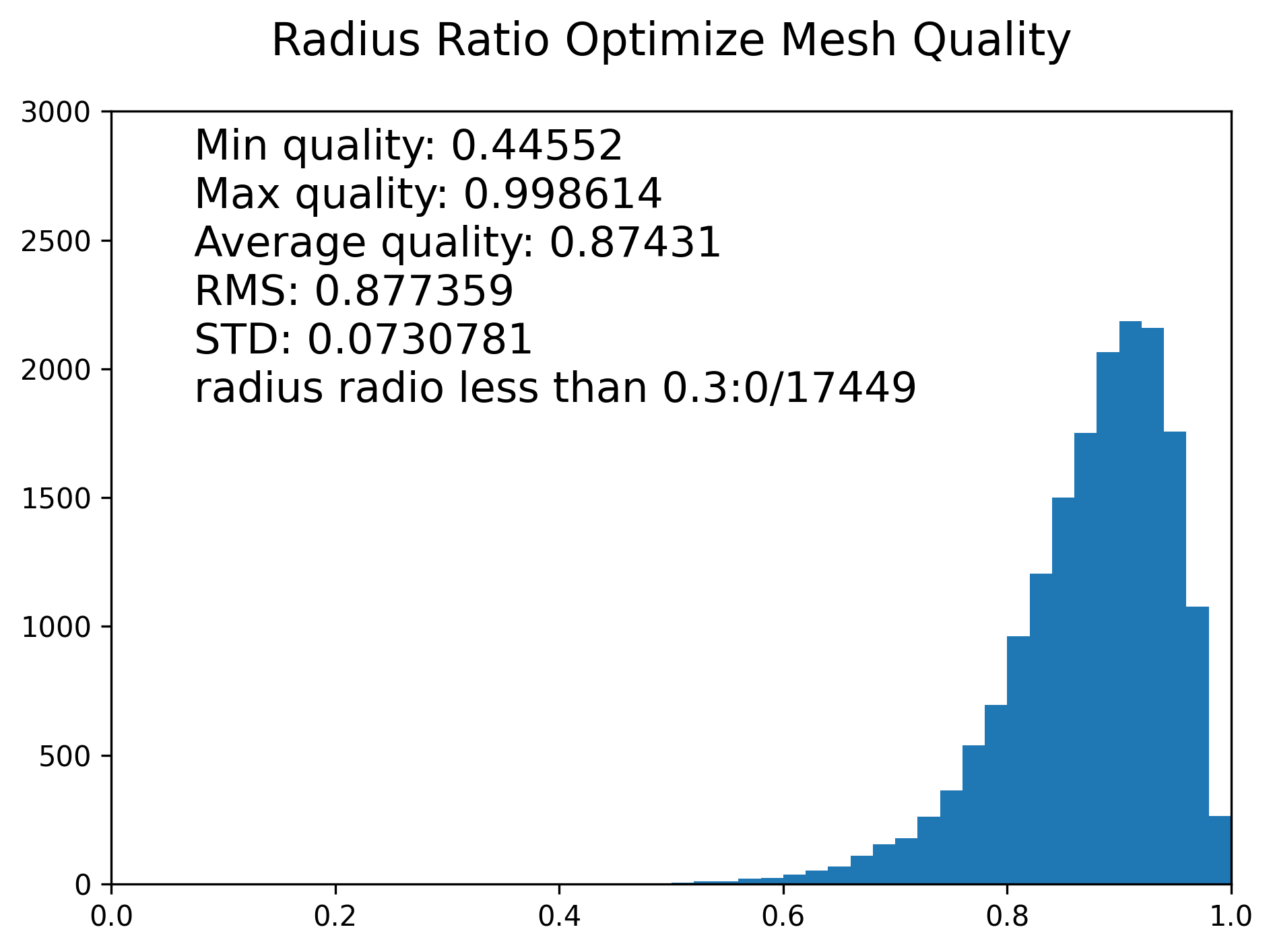}}
\hspace{0.01\linewidth}
\subfloat[Exude]{
\includegraphics[width=0.2\linewidth]{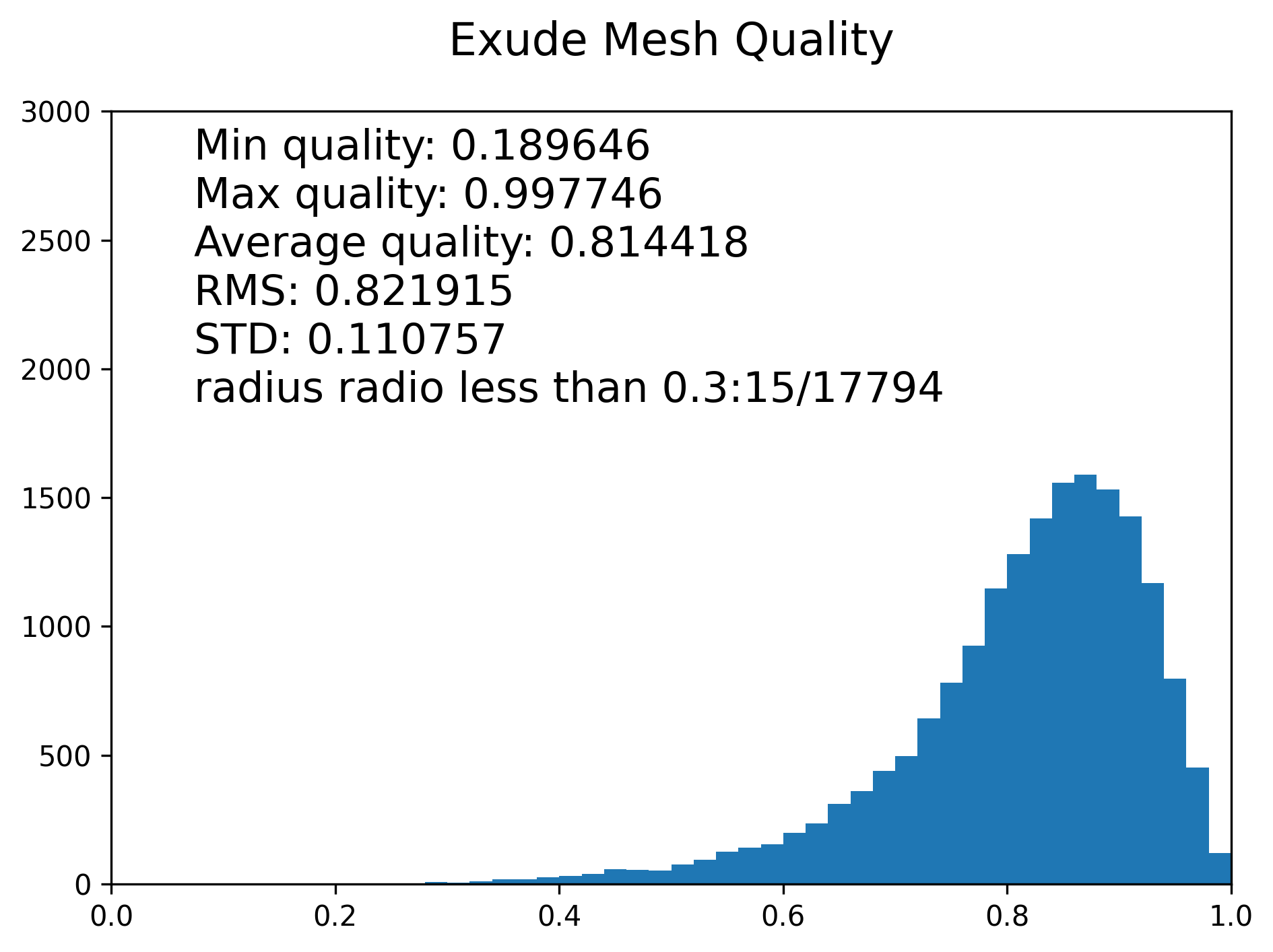}}
\vfill
\subfloat[Precondition RRE]{
\includegraphics[width=0.2\linewidth]{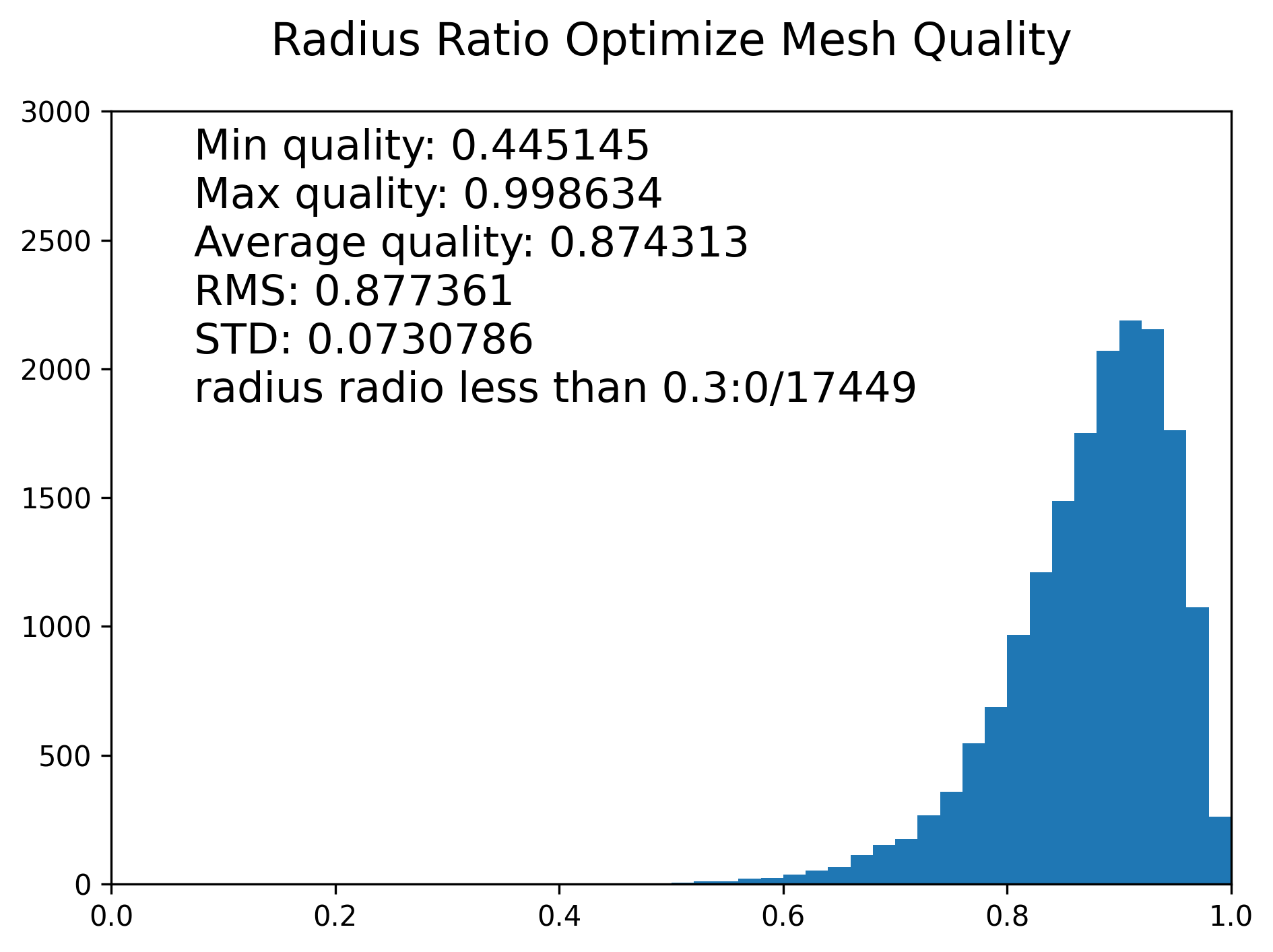}}
\hspace{0.01\linewidth}
\subfloat[ODT]{
\includegraphics[width=0.2\linewidth]{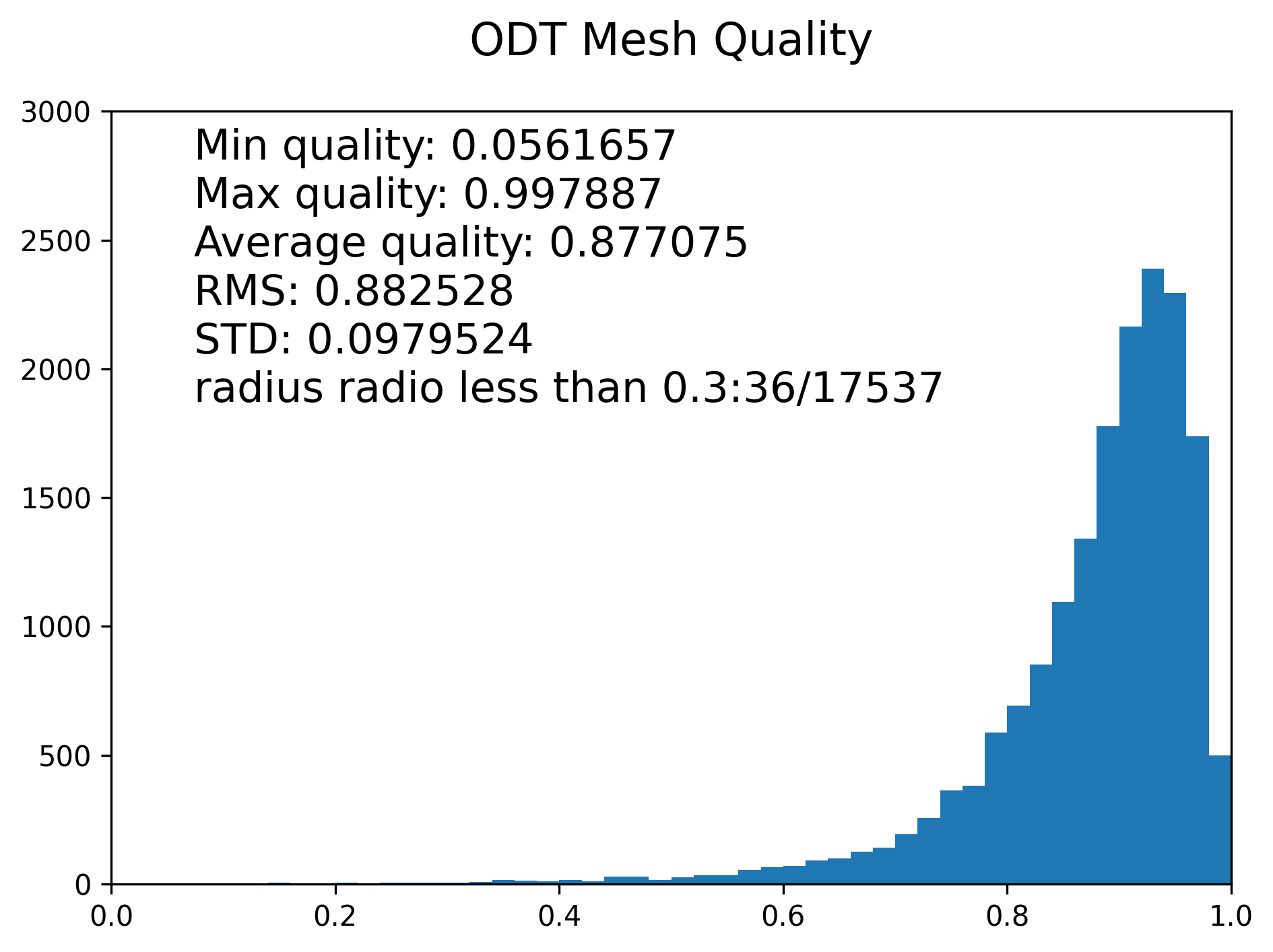}}
\hspace{0.01\linewidth}
\subfloat[ODT+RRE]{
\includegraphics[width=0.2\linewidth]{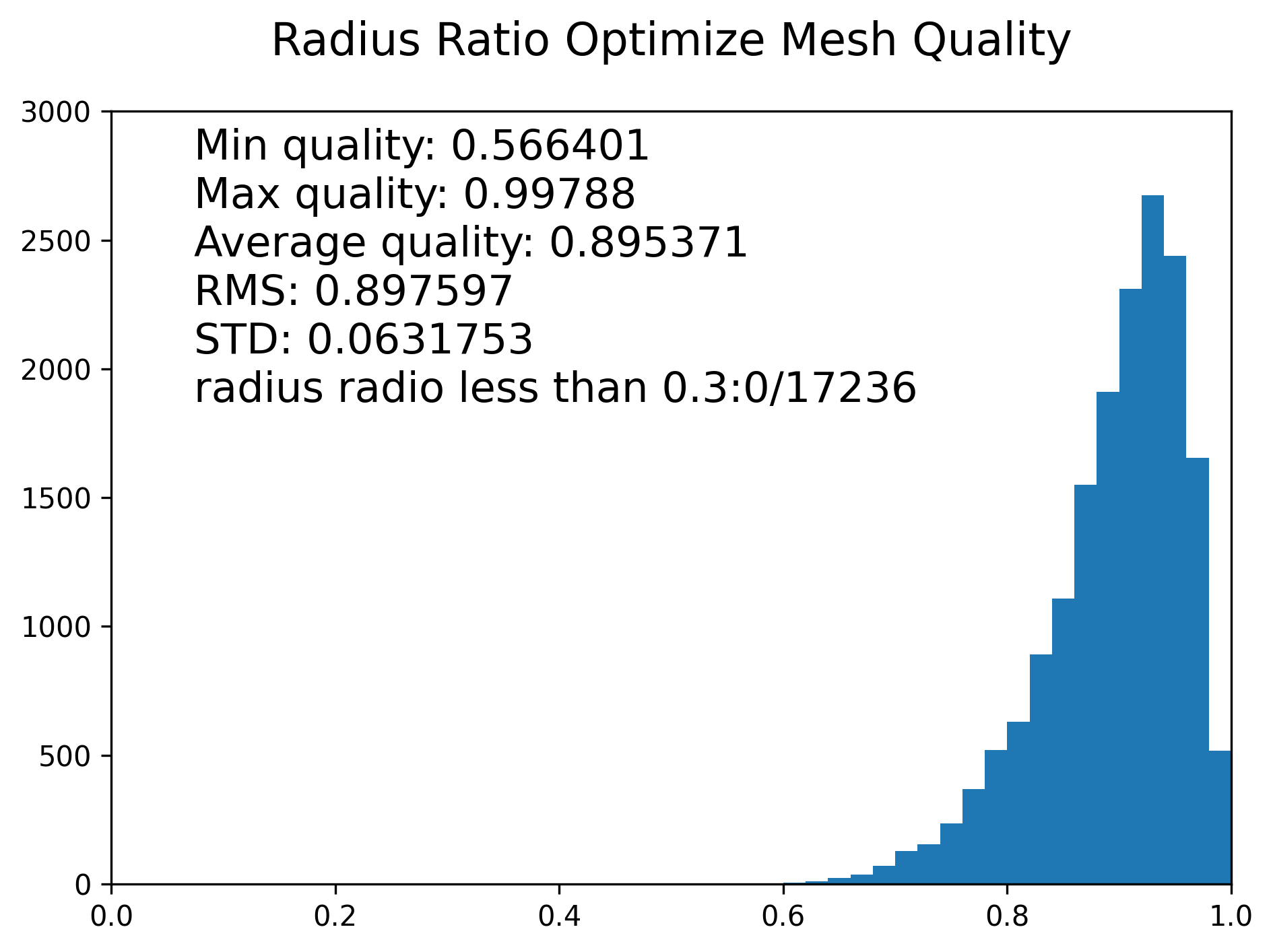}}
\hspace{0.01\linewidth}
\subfloat[ODT+Exude]{
\includegraphics[width=0.2\linewidth]{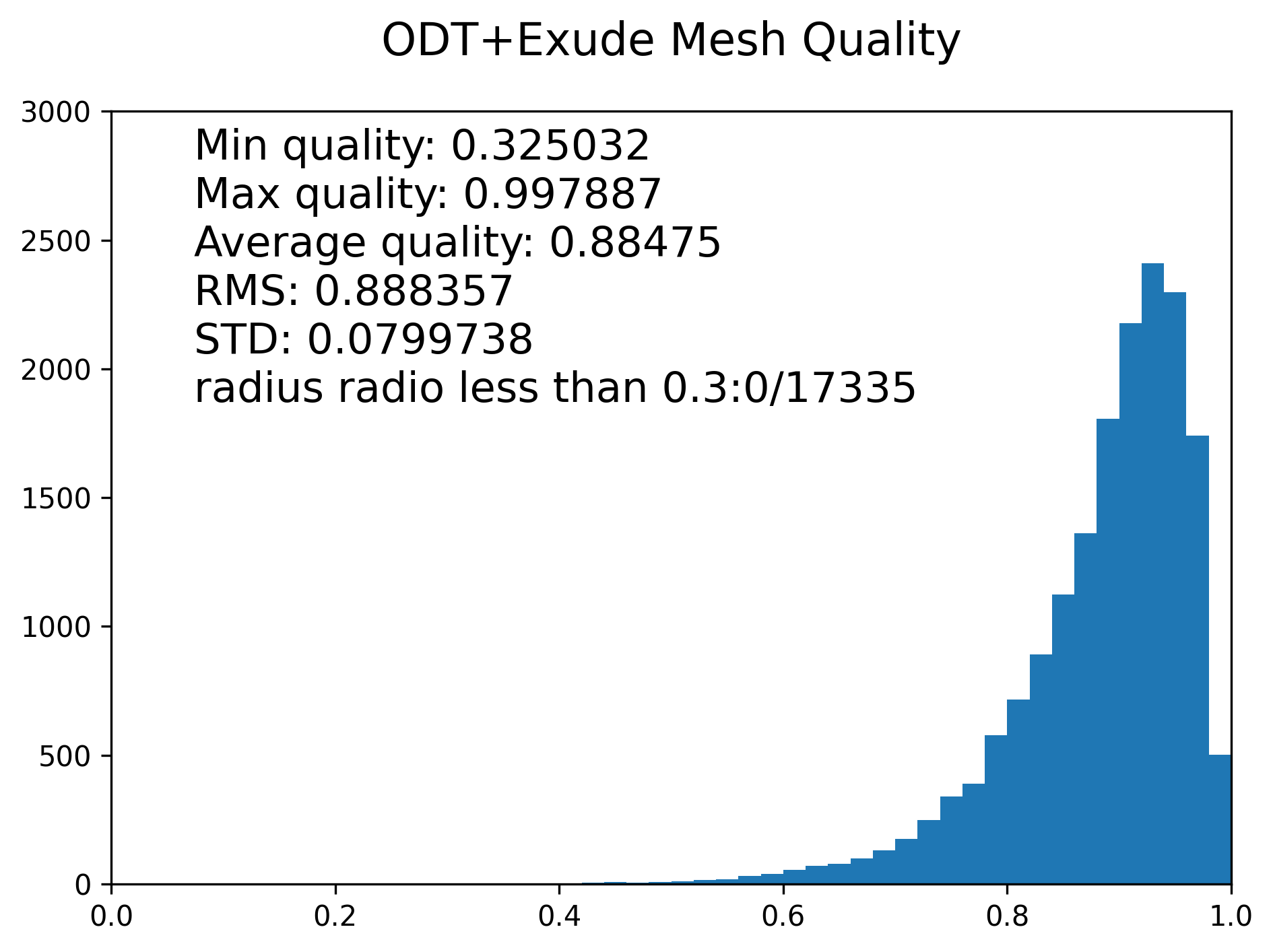}}
\vfill
\subfloat[Init min dihedral angle]{
\includegraphics[width=0.2\linewidth]{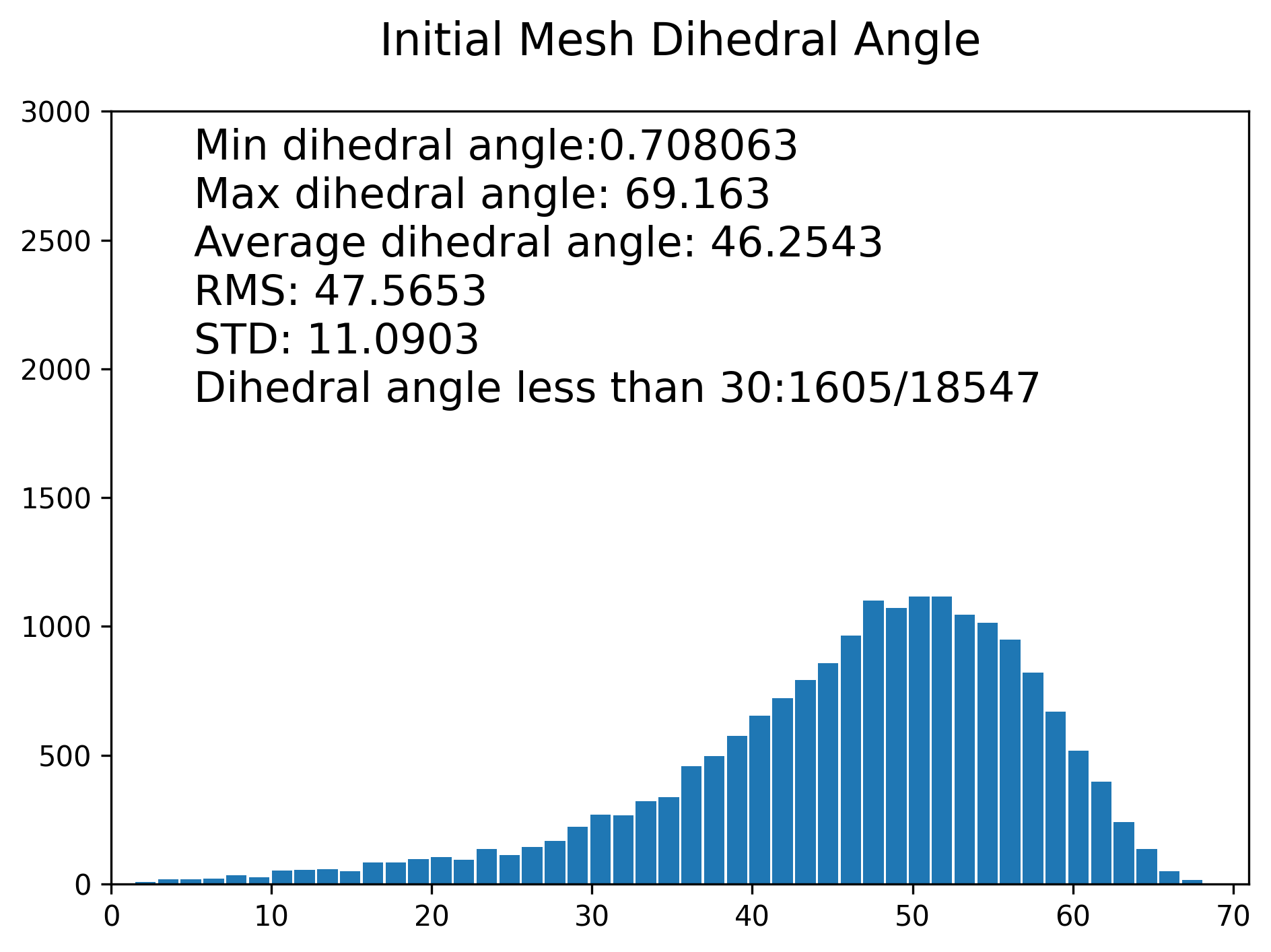}}
\subfloat[ODT min dihedral angle]{
\includegraphics[width=0.2\linewidth]{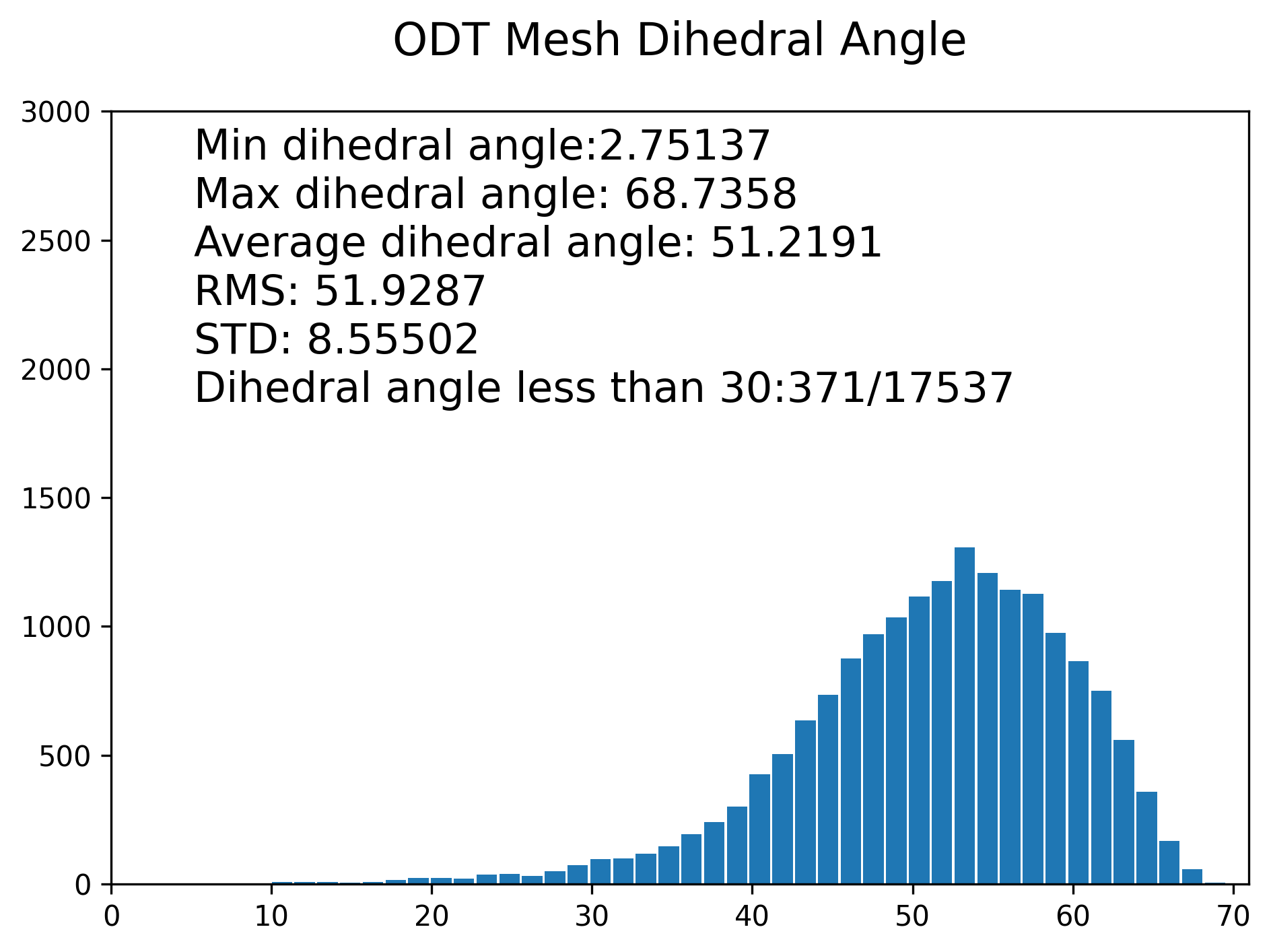}}
\hspace{0.01\linewidth}
\subfloat[RRE min dihedral angle]{
\includegraphics[width=0.2\linewidth]{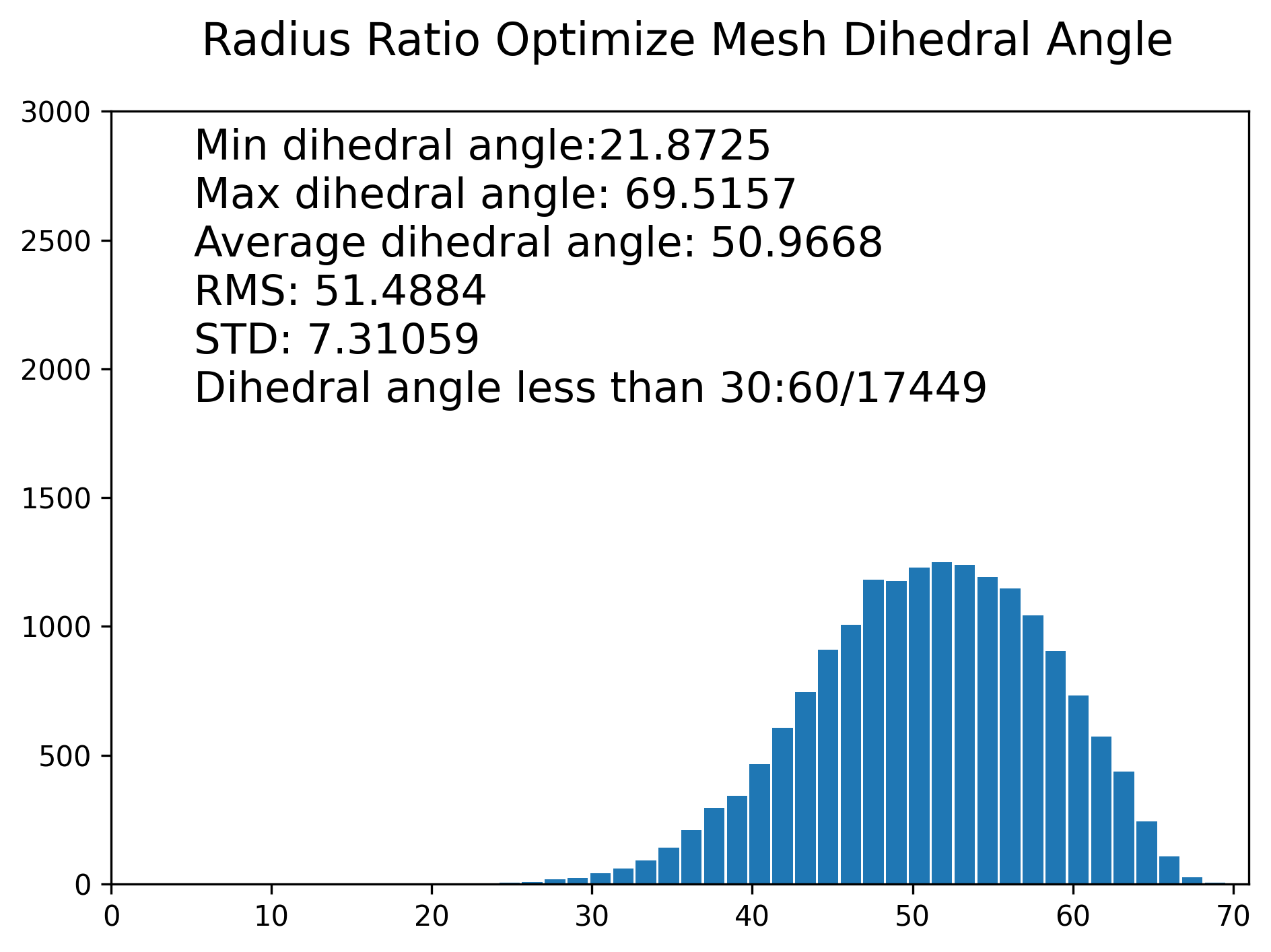}}
\hspace{0.01\linewidth}
\subfloat[Exude min dihedral angle]{
\includegraphics[width=0.2\linewidth]{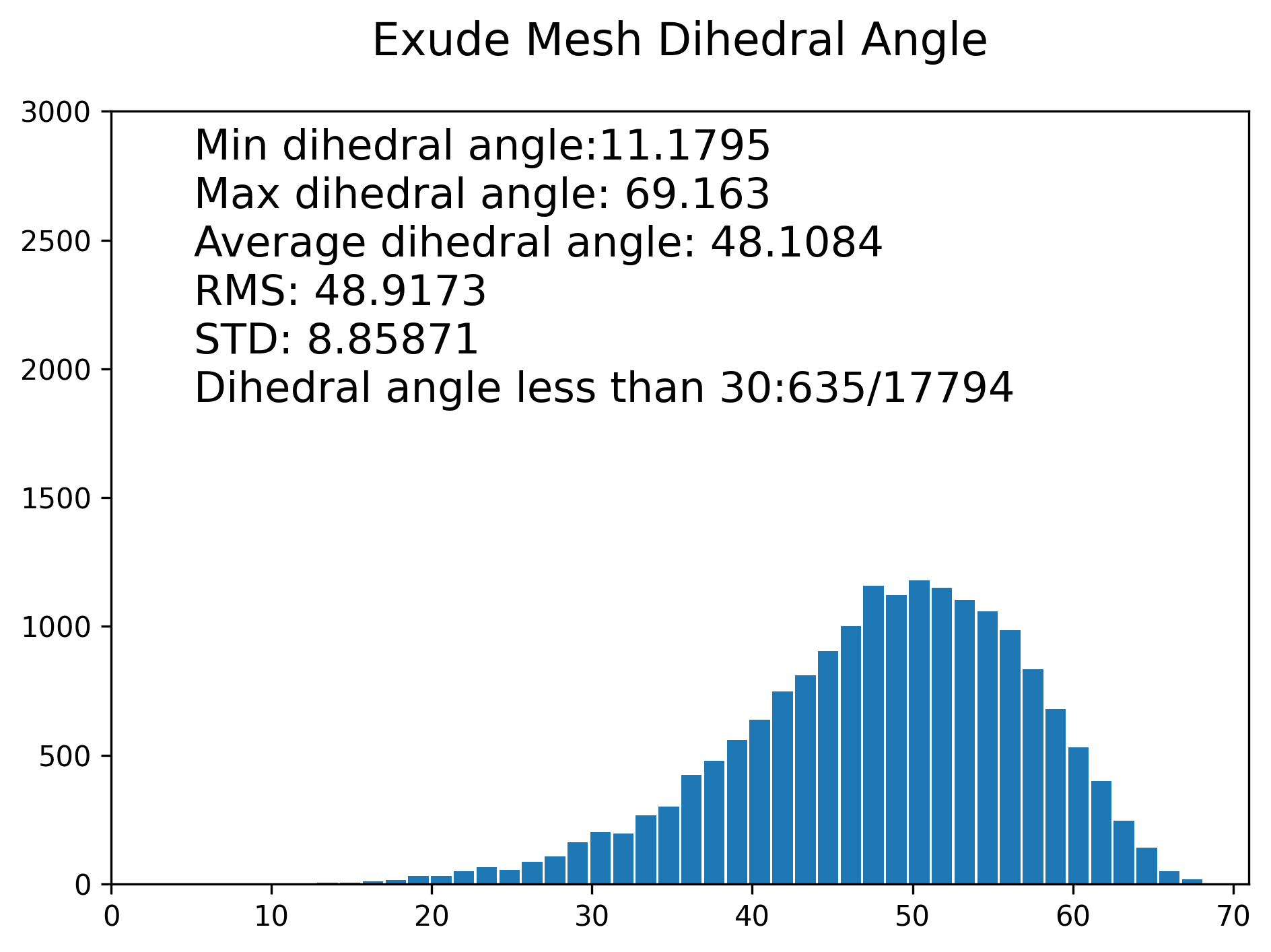}}
\caption{Sphere}
\label{fig:sphere}
\end{figure}

\begin{figure}[htbp]
\centering
\subfloat[Model]{
\includegraphics[width=0.2\linewidth]{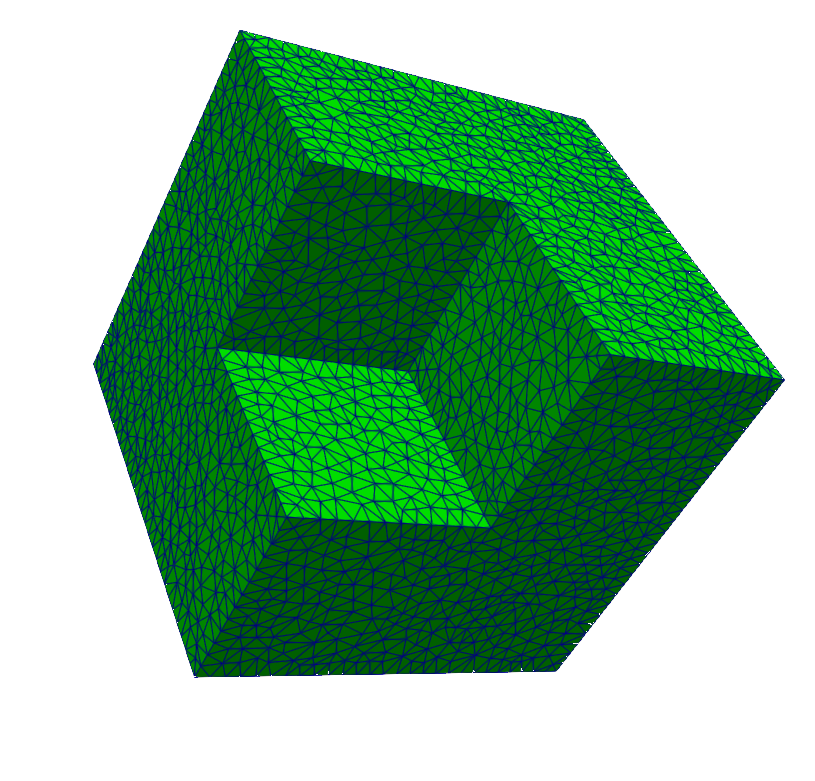}}
\hspace{0.01\linewidth}
\subfloat[Init quality]{
\includegraphics[width=0.2\linewidth]{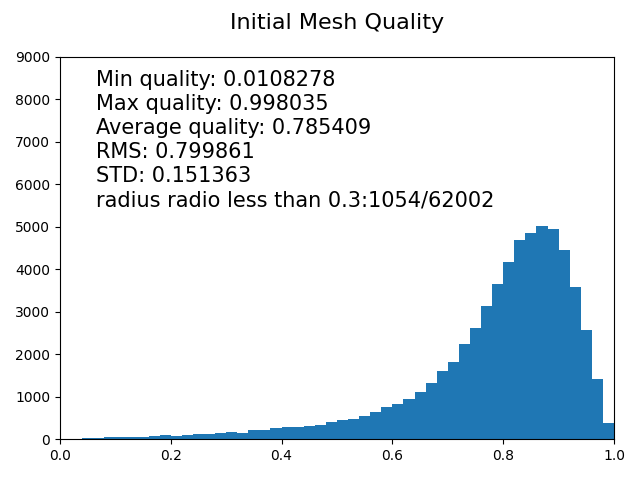}}
\hspace{0.01\linewidth}
\subfloat[RRE]{
\includegraphics[width=0.2\linewidth]{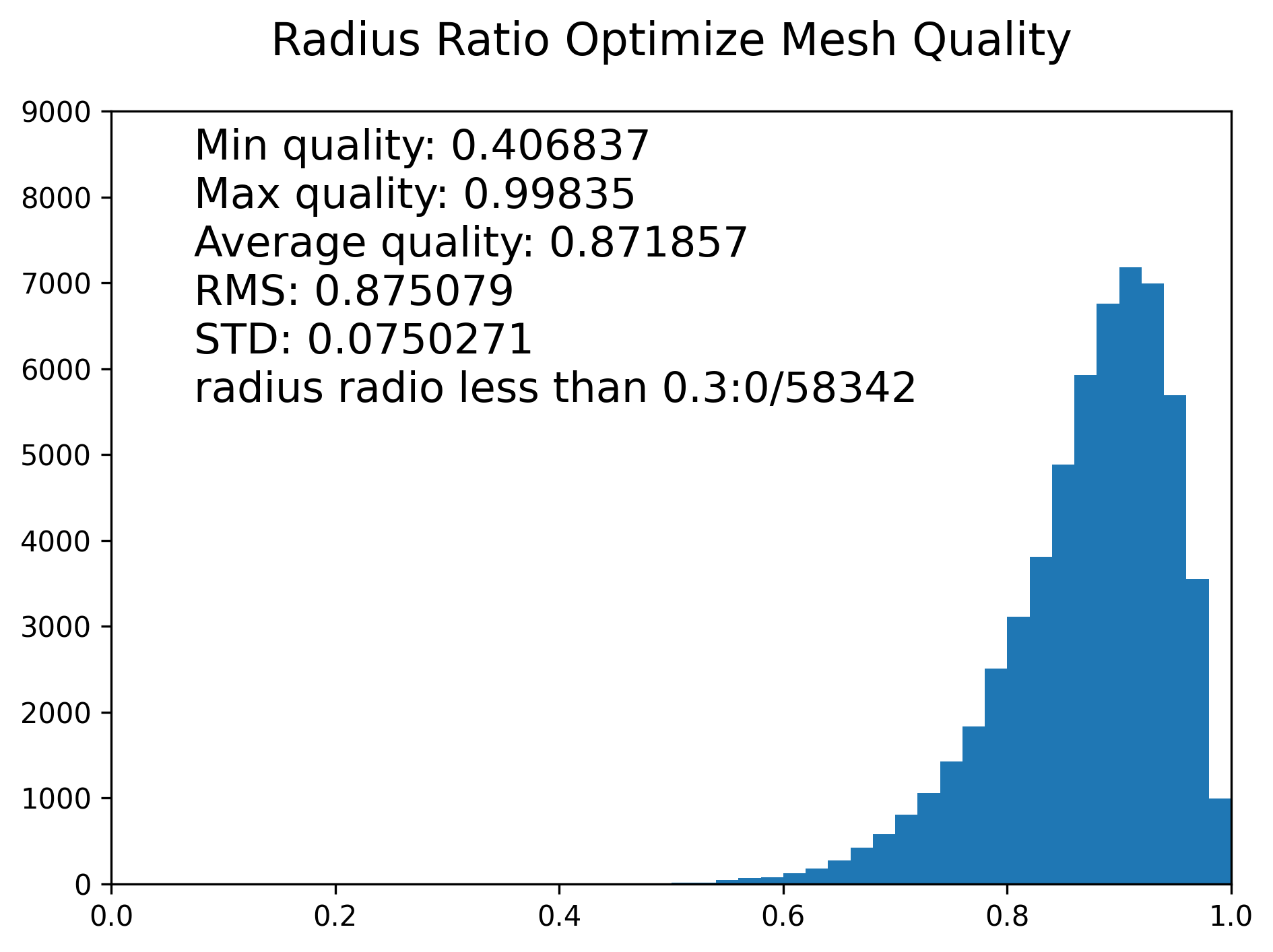}}
\hspace{0.01\linewidth}
\subfloat[Exude]{
\includegraphics[width=0.2\linewidth]{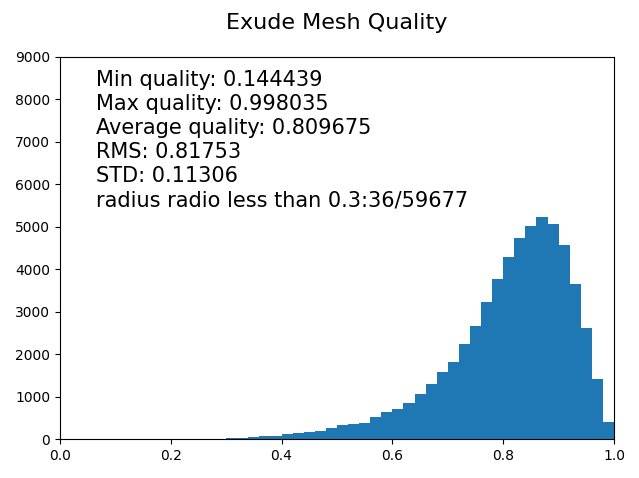}}
\vfill
\subfloat[Precondition RRE]{
\includegraphics[width=0.2\linewidth]{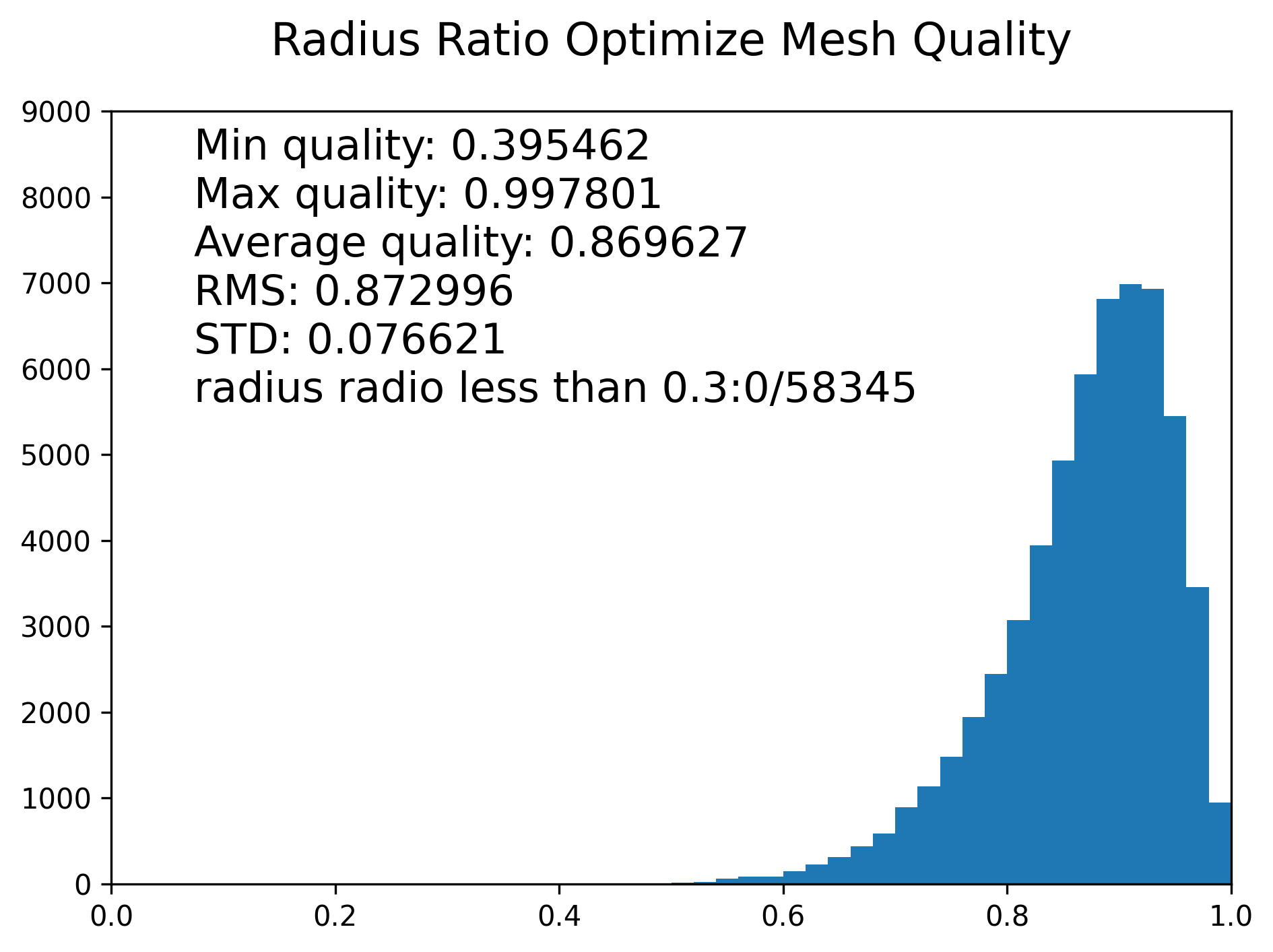}}
\hspace{0.01\linewidth}
\subfloat[ODT]{
\includegraphics[width=0.2\linewidth]{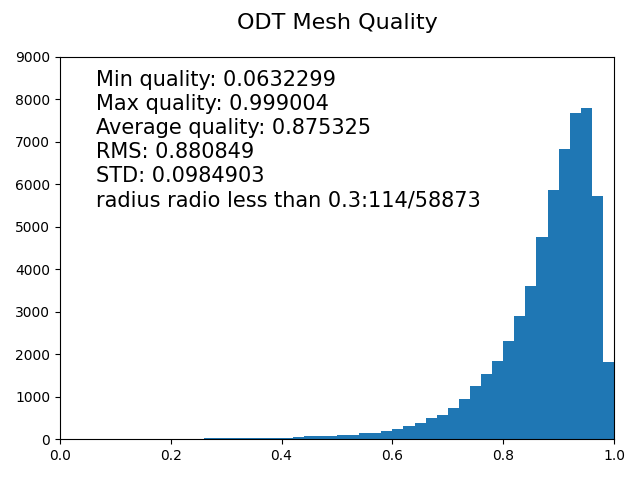}}
\hspace{0.01\linewidth}
\subfloat[ODT+RRE]{
\includegraphics[width=0.2\linewidth]{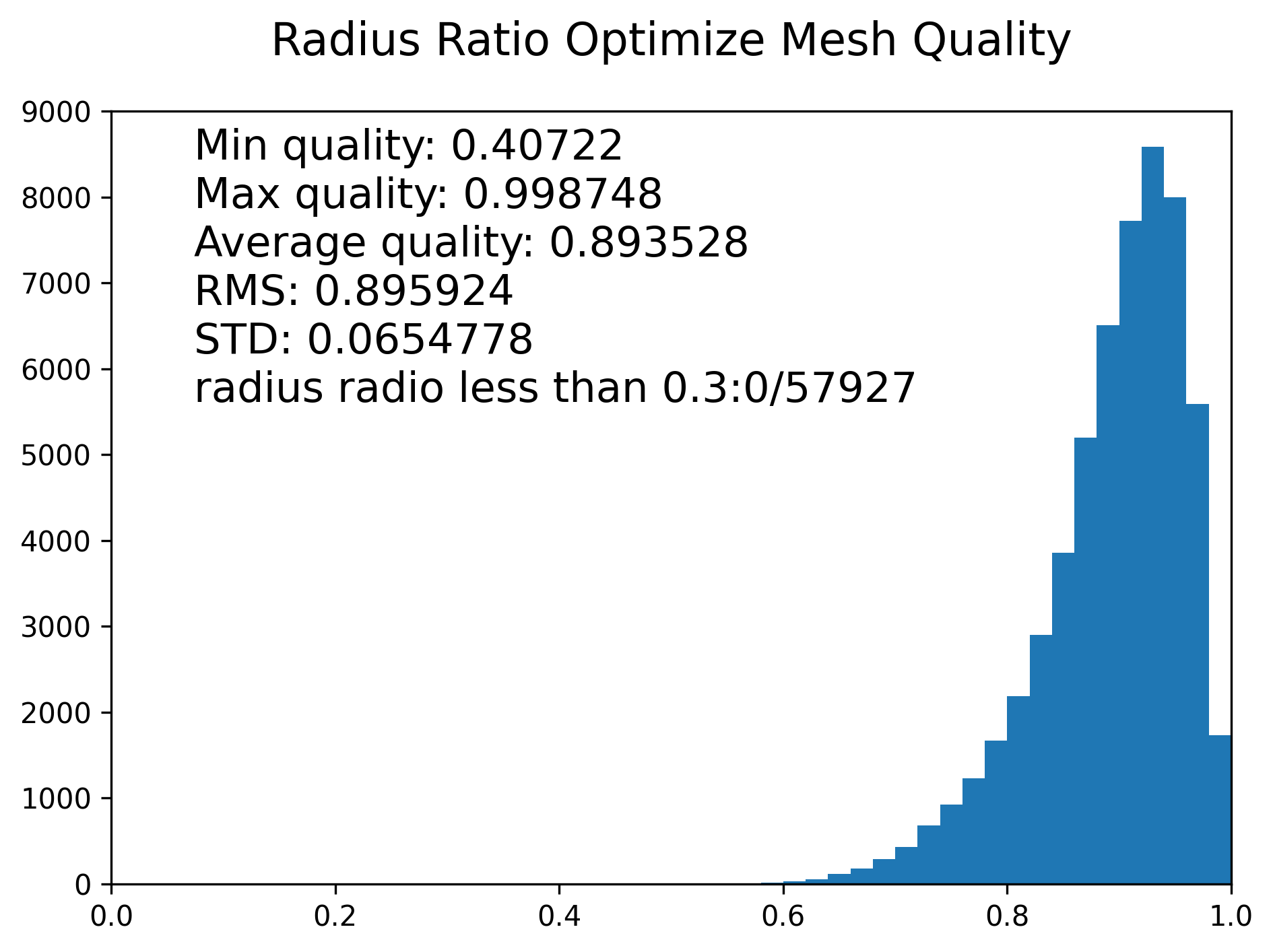}}
\hspace{0.01\linewidth}
\subfloat[ODT+Exude]{
\includegraphics[width=0.2\linewidth]{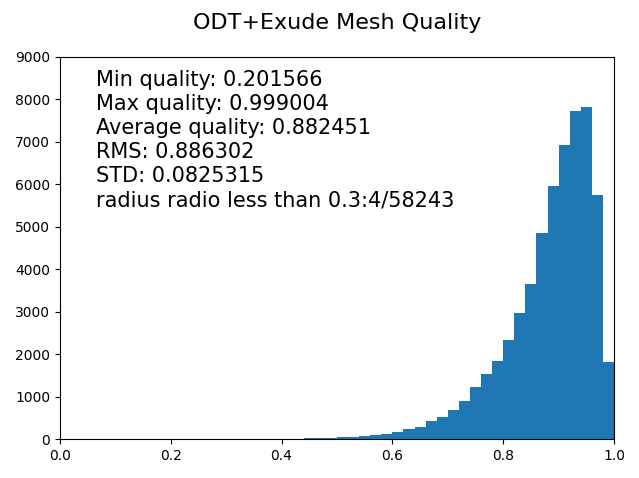}}
\vfill
\subfloat[Init min dihedral angle]{
\includegraphics[width=0.2\linewidth]{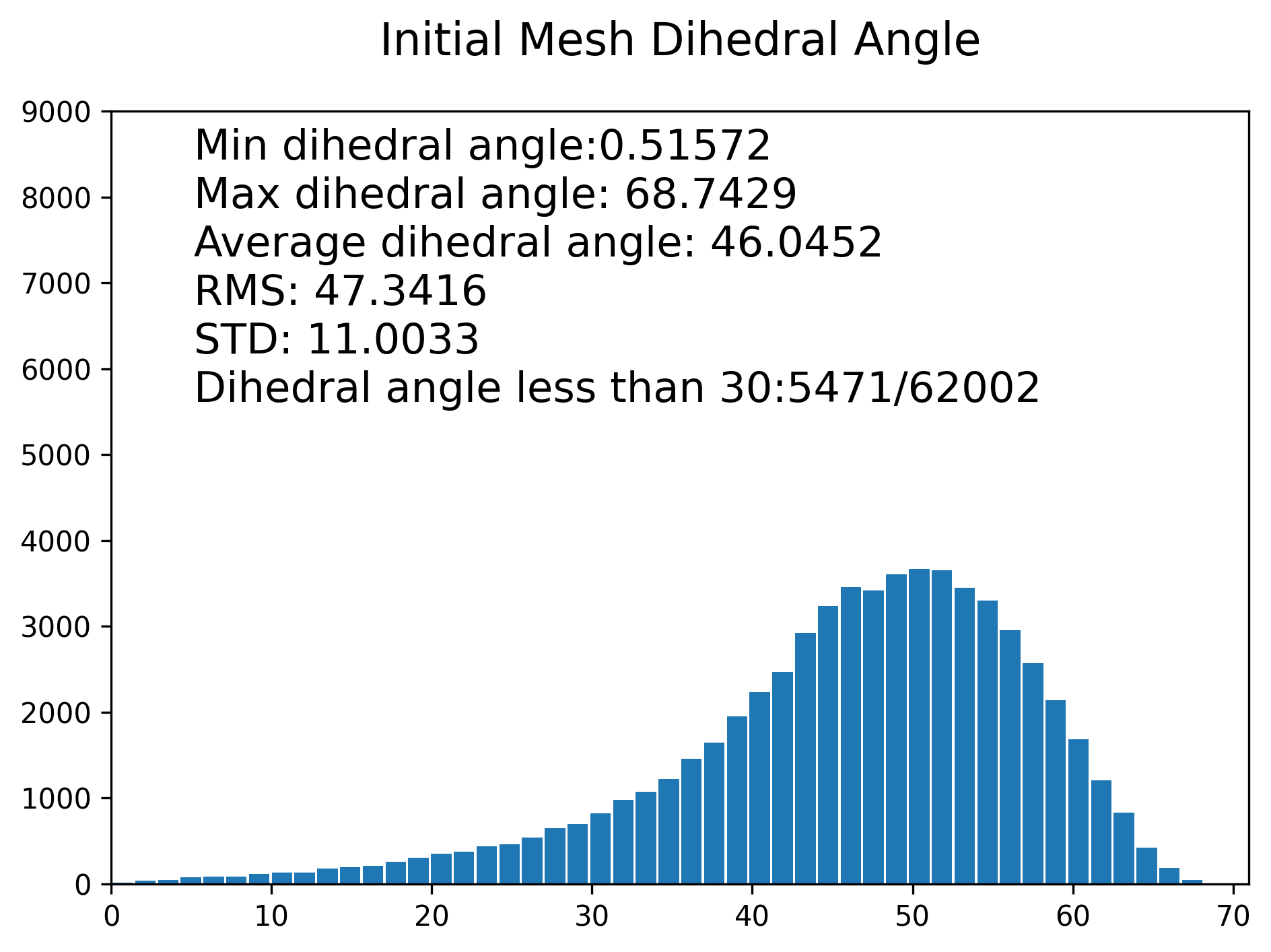}}
\subfloat[ODT min dihedral angle]{
\includegraphics[width=0.2\linewidth]{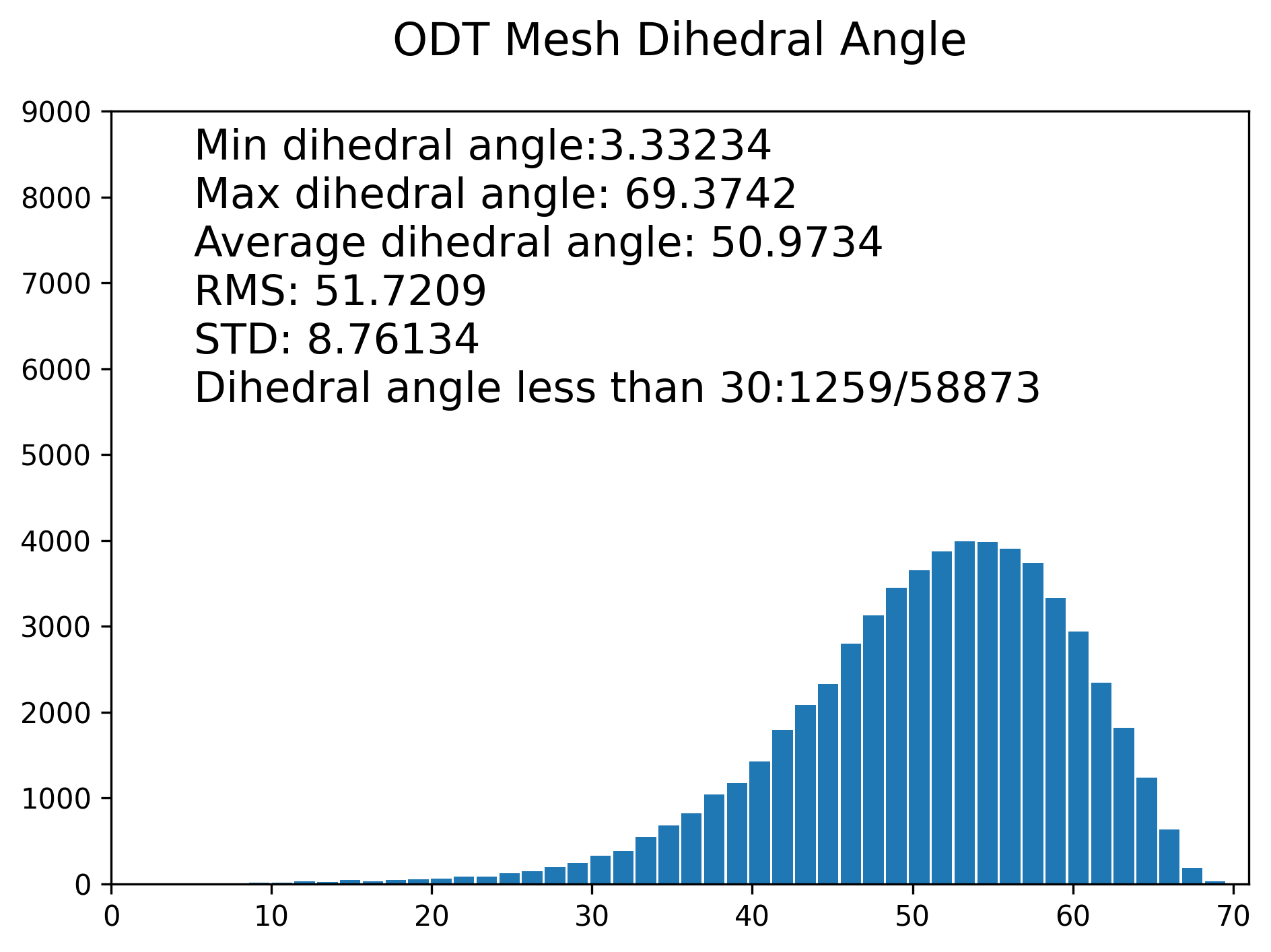}}
\hspace{0.01\linewidth}
\subfloat[RRE min dihedral angle]{
\includegraphics[width=0.2\linewidth]{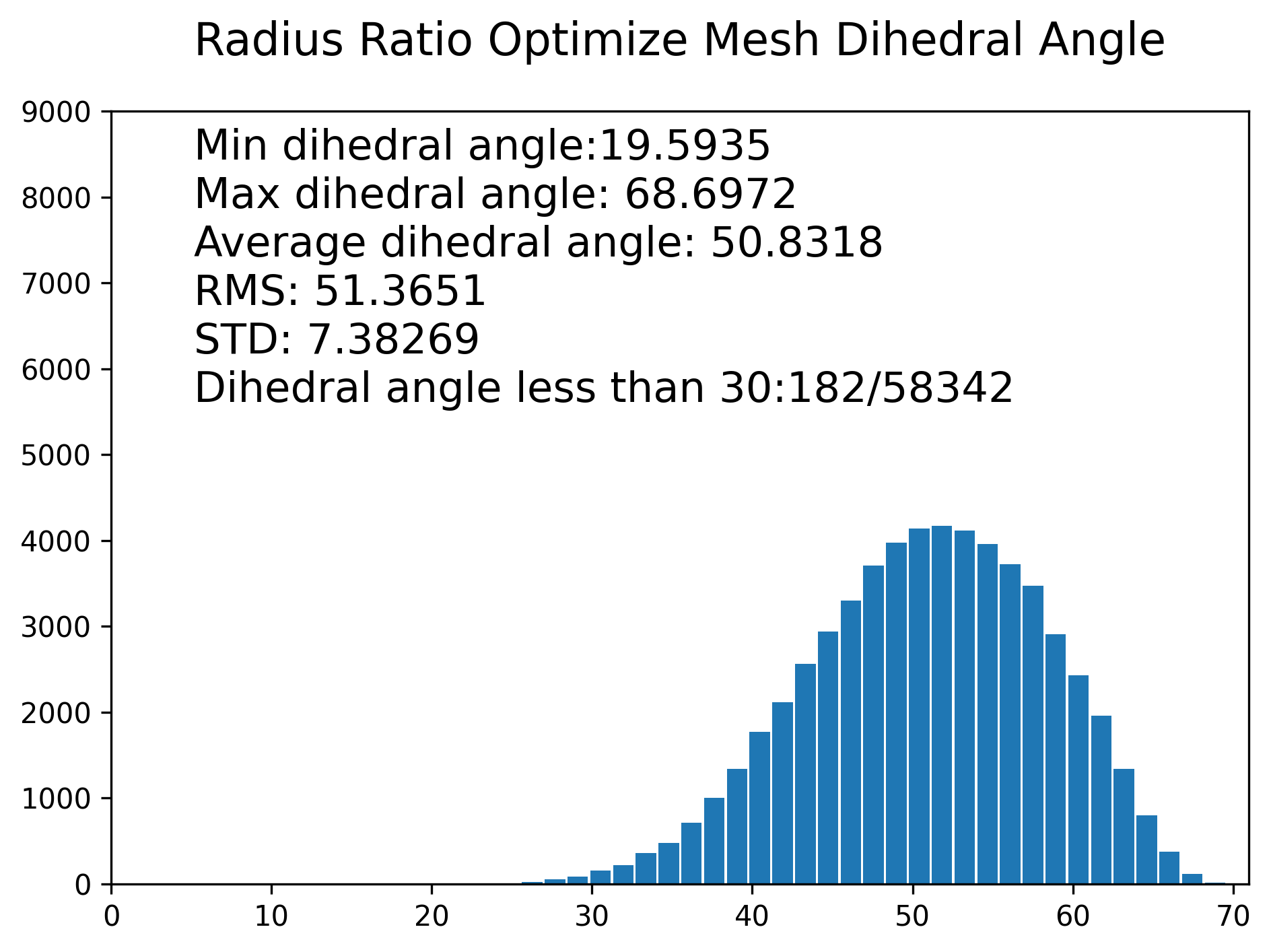}}
\hspace{0.01\linewidth}
\subfloat[Exude min dihedral angle]{
\includegraphics[width=0.2\linewidth]{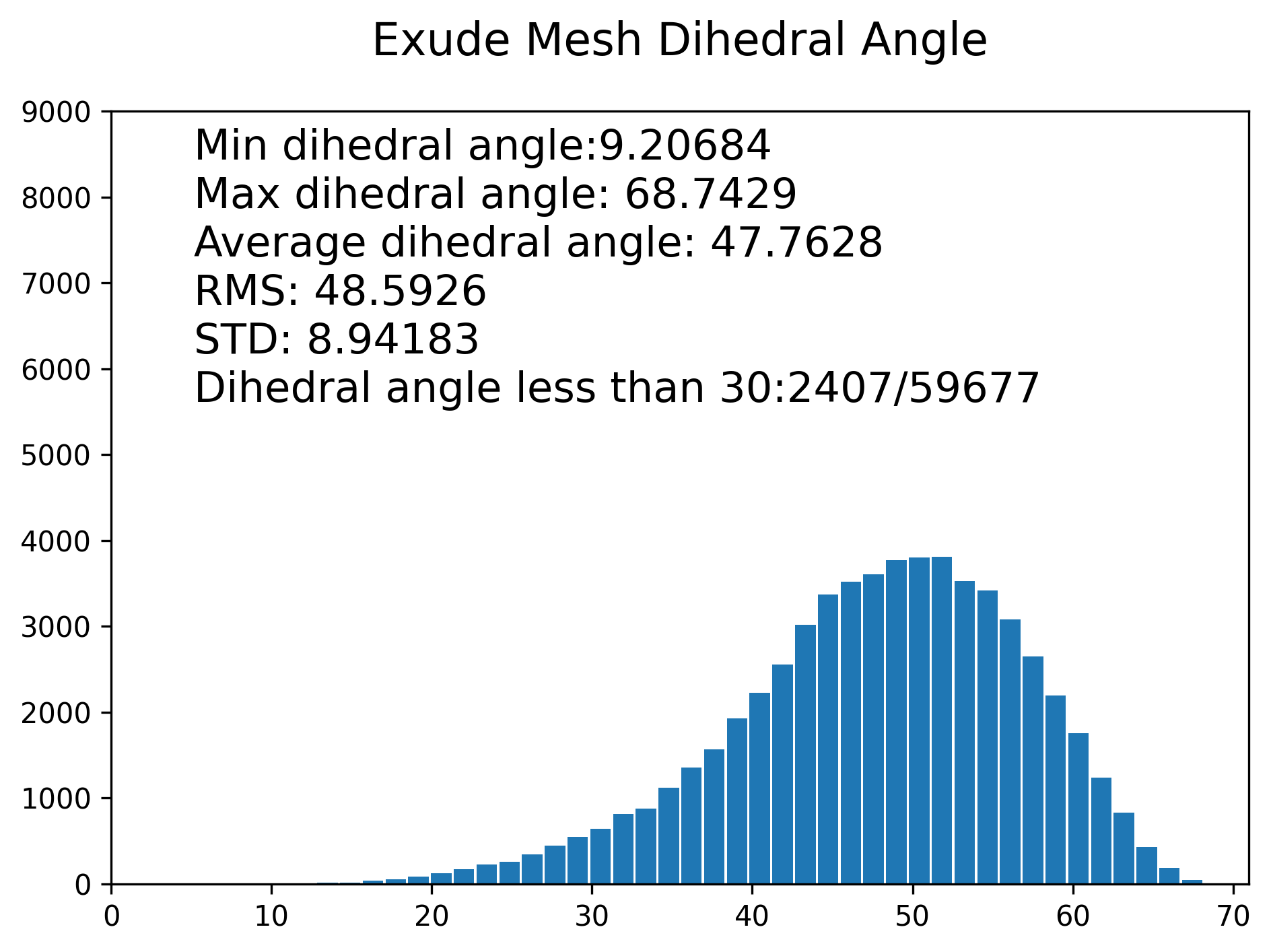}}
\caption{L-shape}
\label{fig:lshape}
\end{figure}

\begin{figure}[htbp]
\centering
\subfloat[Model]{
\includegraphics[width=0.2\linewidth]{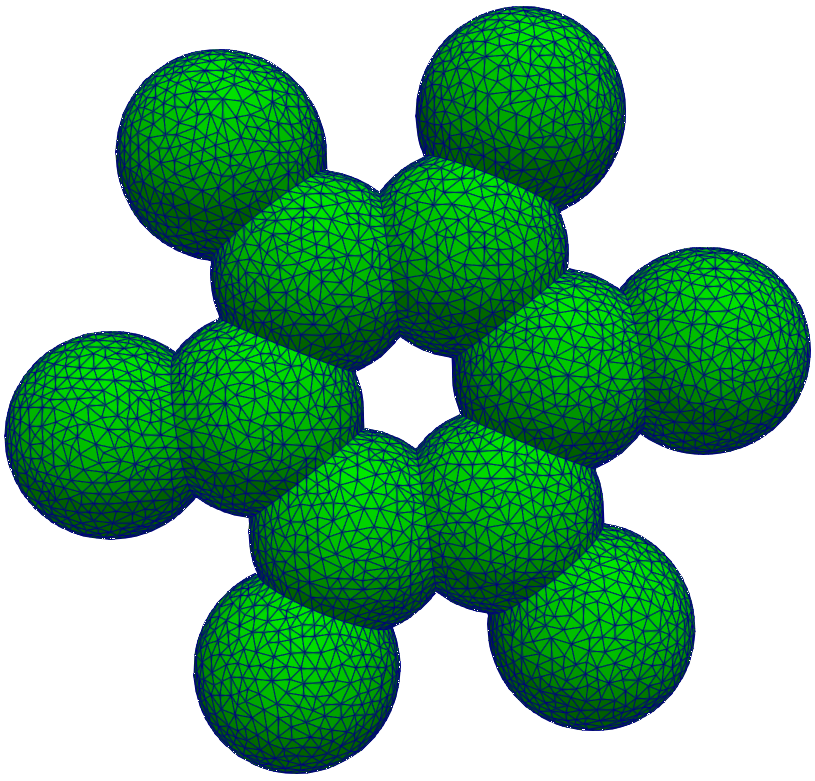}}
\hspace{0.01\linewidth}
\subfloat[Init quality]{
\includegraphics[width=0.2\linewidth]{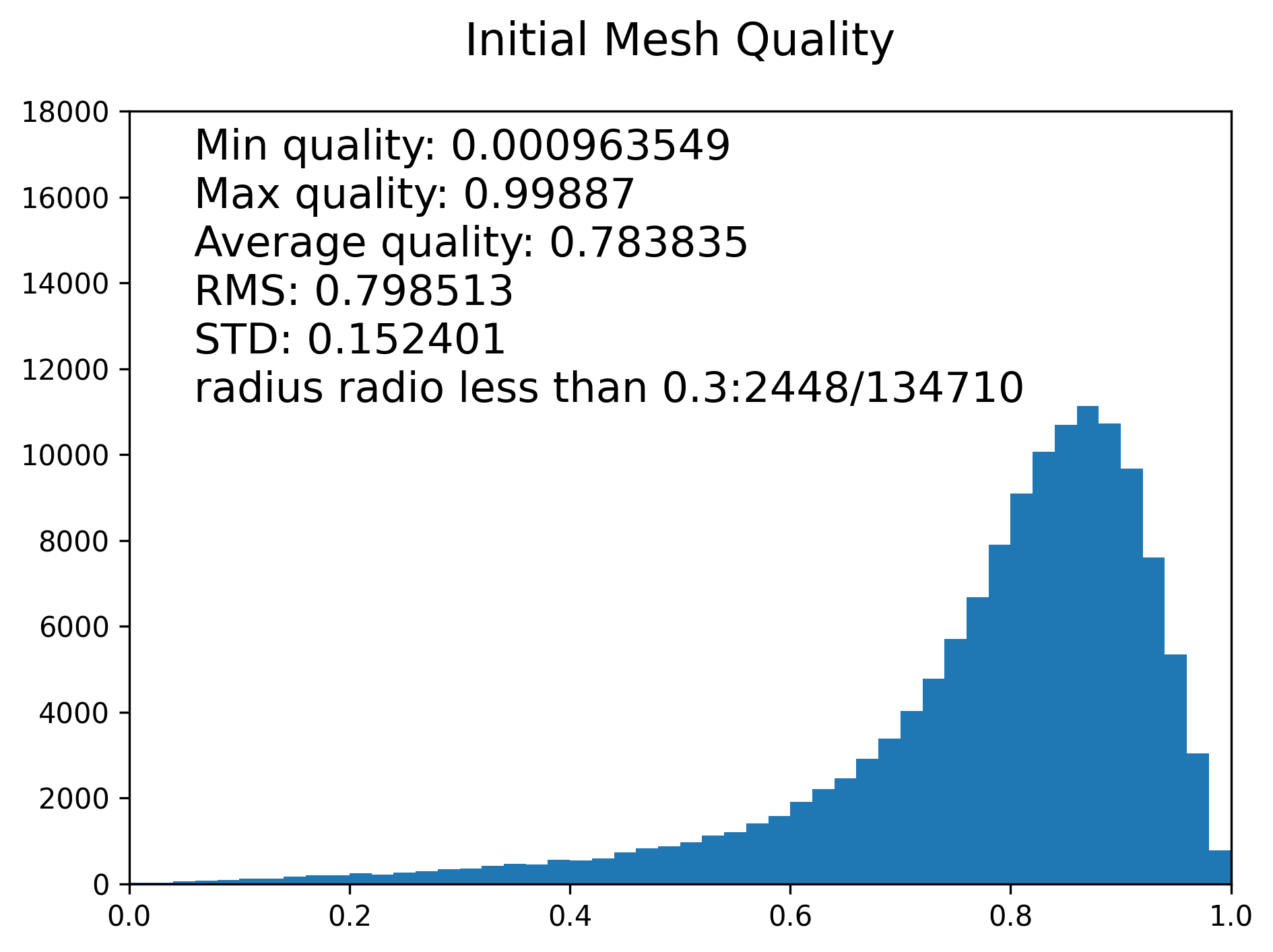}}
\hspace{0.01\linewidth}
\subfloat[RRE]{
\includegraphics[width=0.2\linewidth]{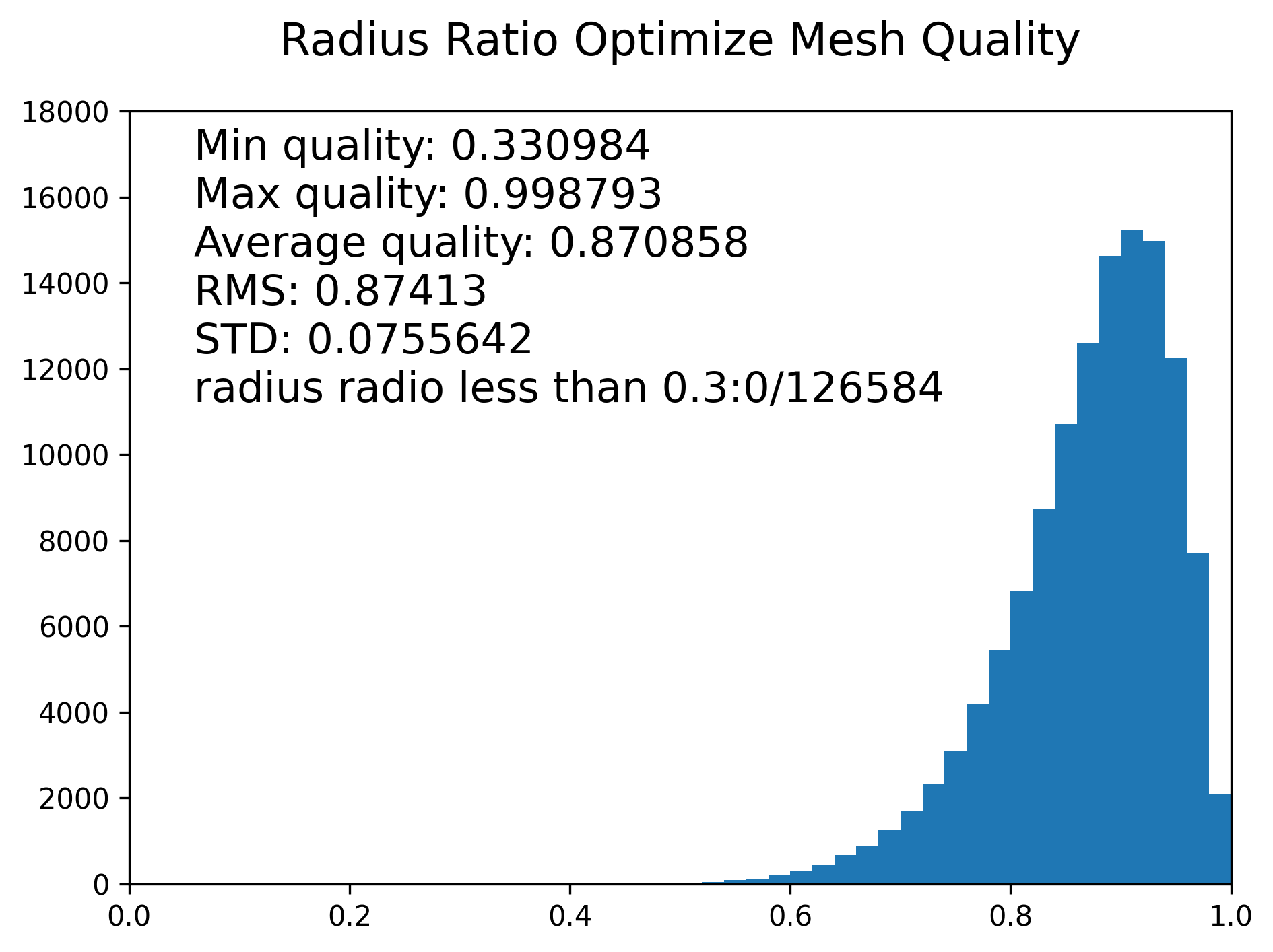}}
\hspace{0.01\linewidth}
\subfloat[Exude]{
\includegraphics[width=0.2\linewidth]{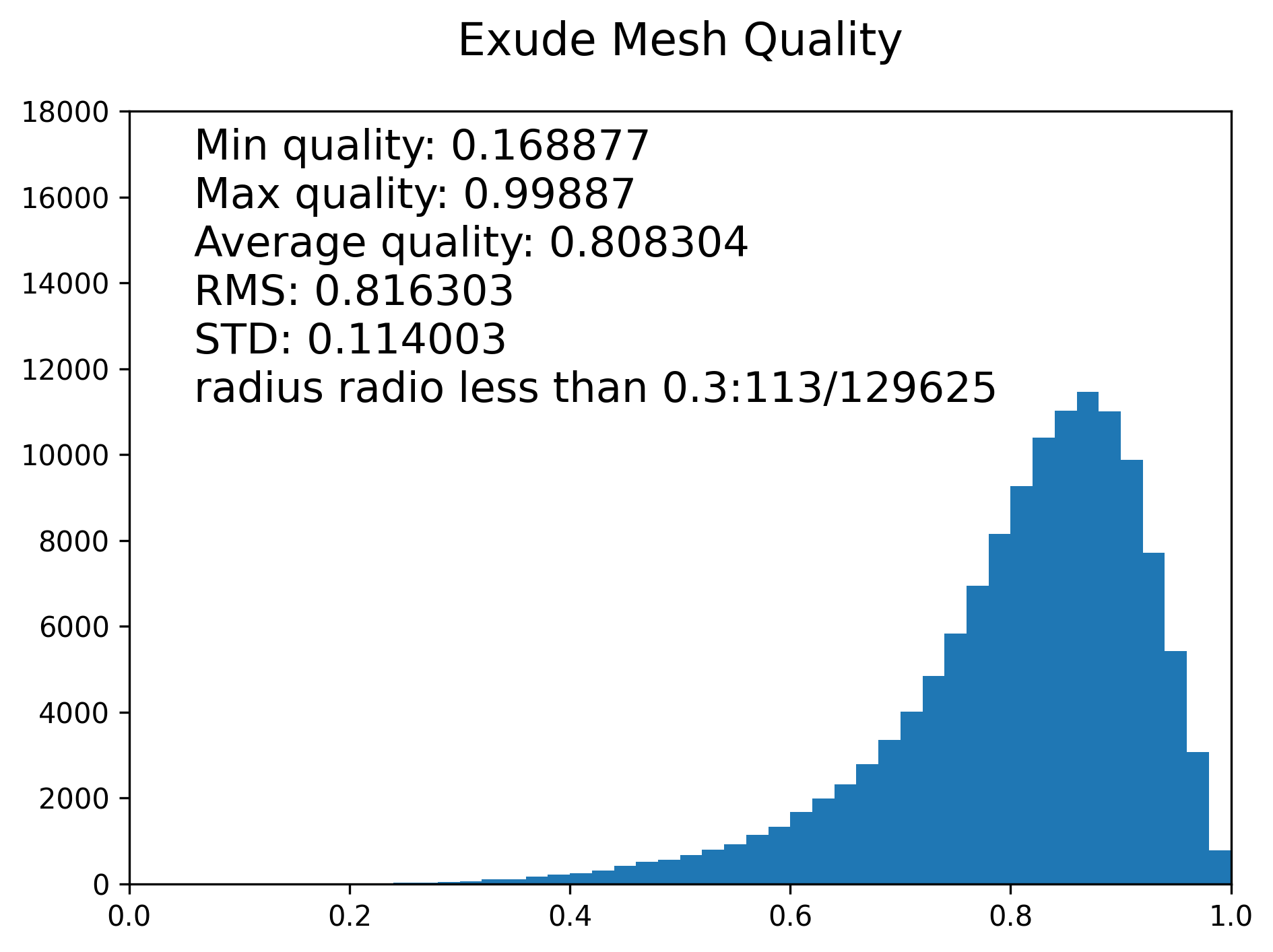}}
\vfill
\subfloat[Precondition RRE]{
\includegraphics[width=0.2\linewidth]{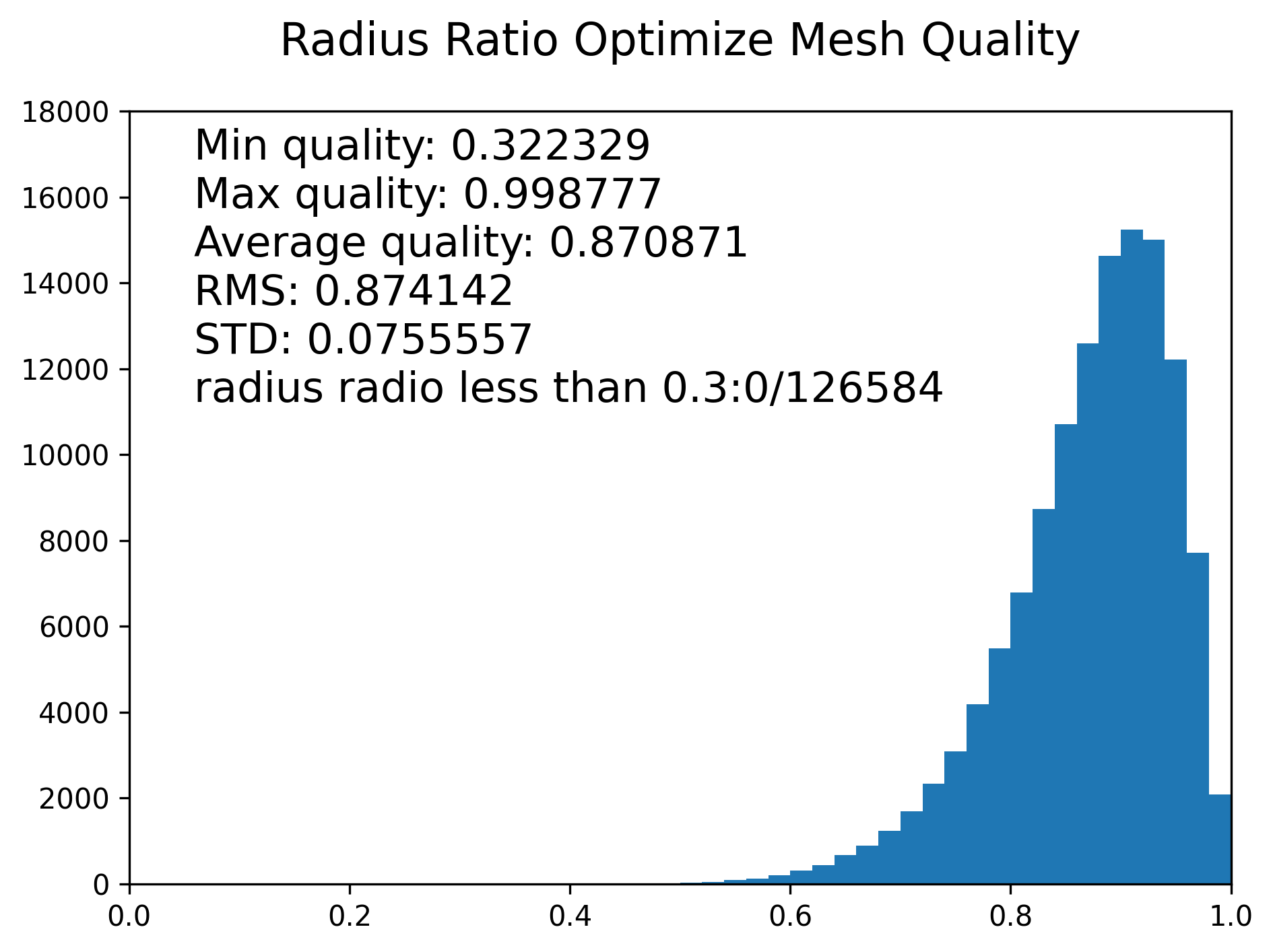}}
\hspace{0.01\linewidth}
\subfloat[ODT]{
\includegraphics[width=0.2\linewidth]{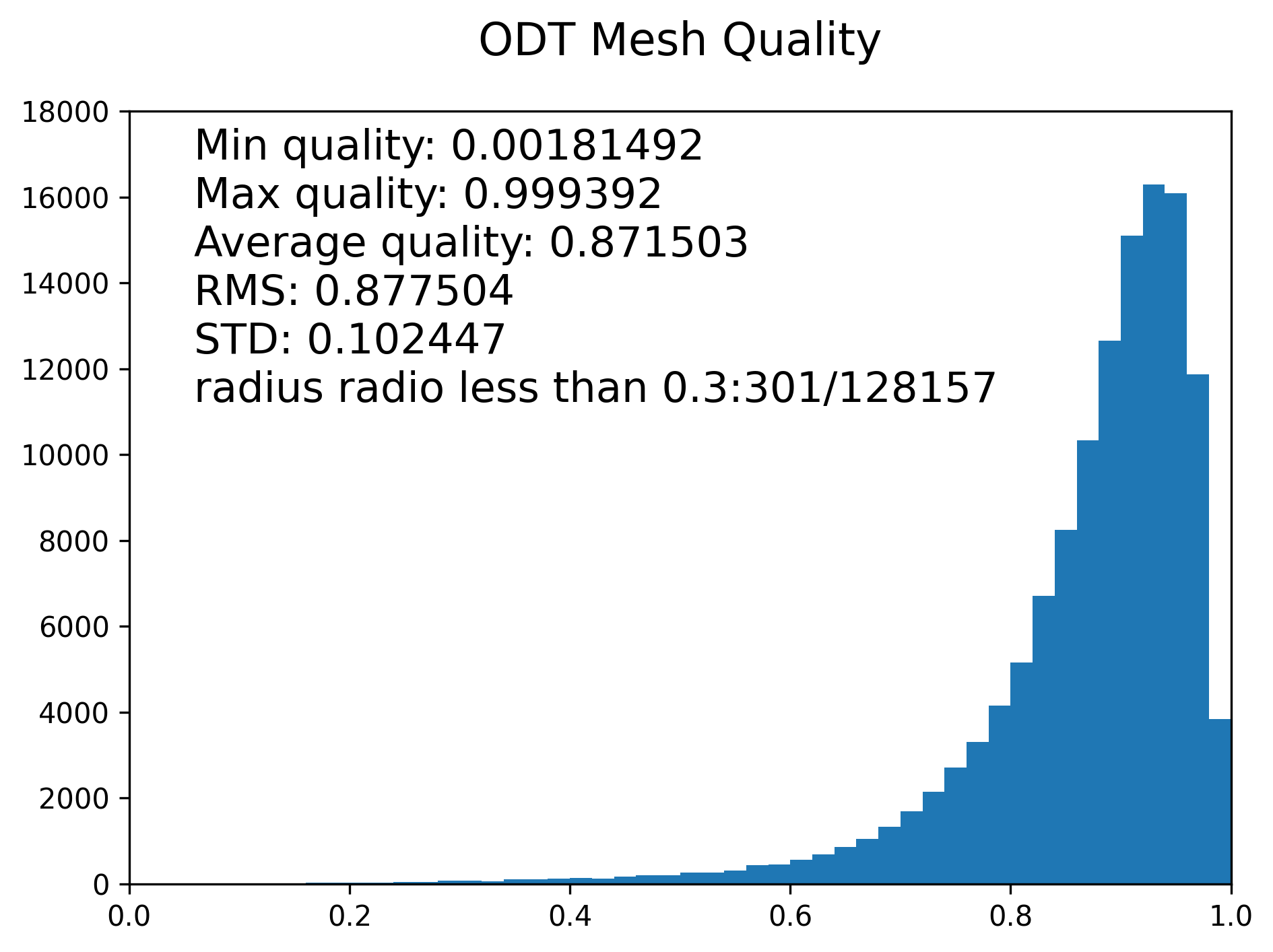}}
\hspace{0.01\linewidth}
\subfloat[ODT + RRE]{
\includegraphics[width=0.2\linewidth]{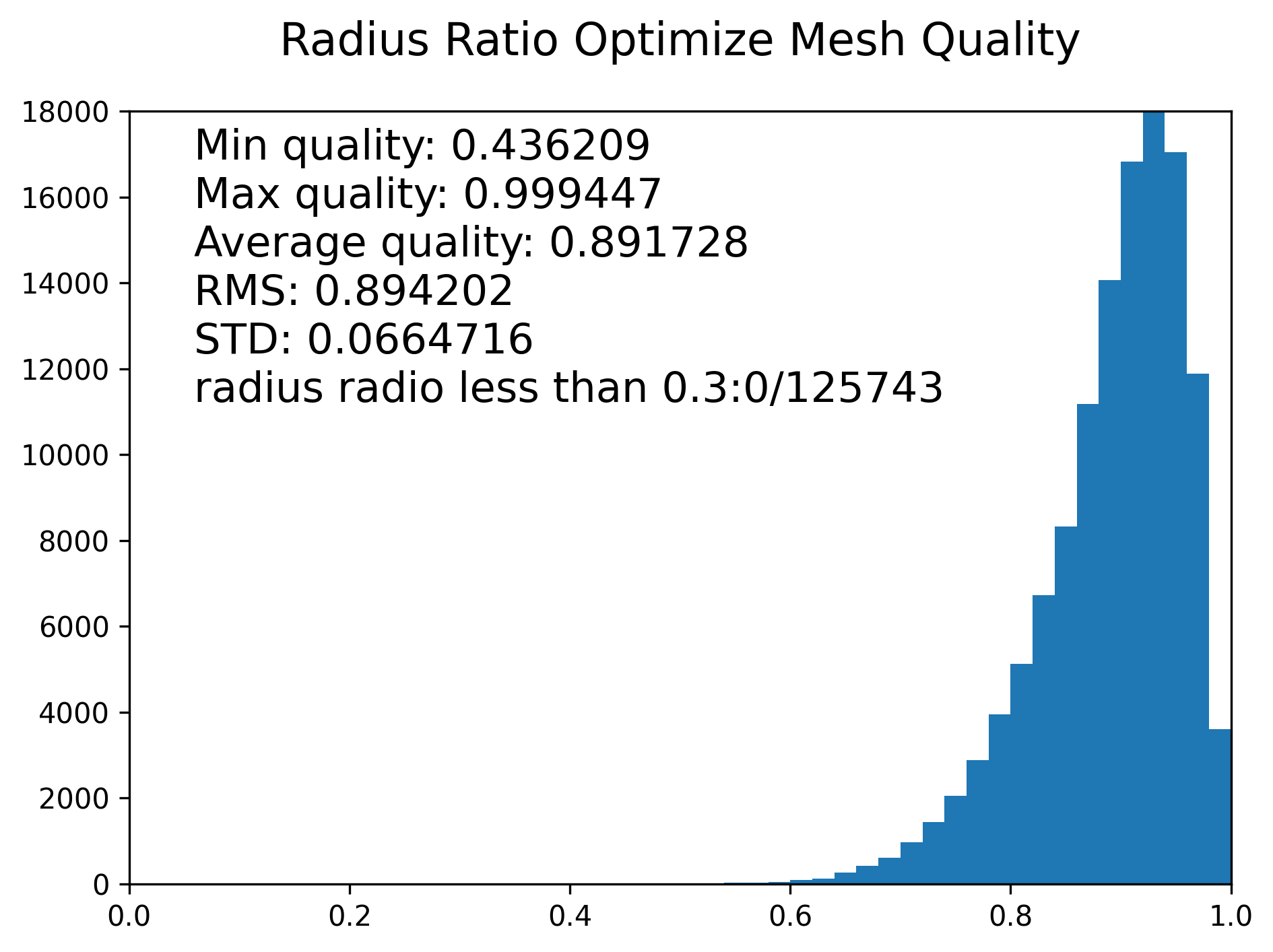}}
\hspace{0.01\linewidth}
\subfloat[ODT + Exude]{
\includegraphics[width=0.2\linewidth]{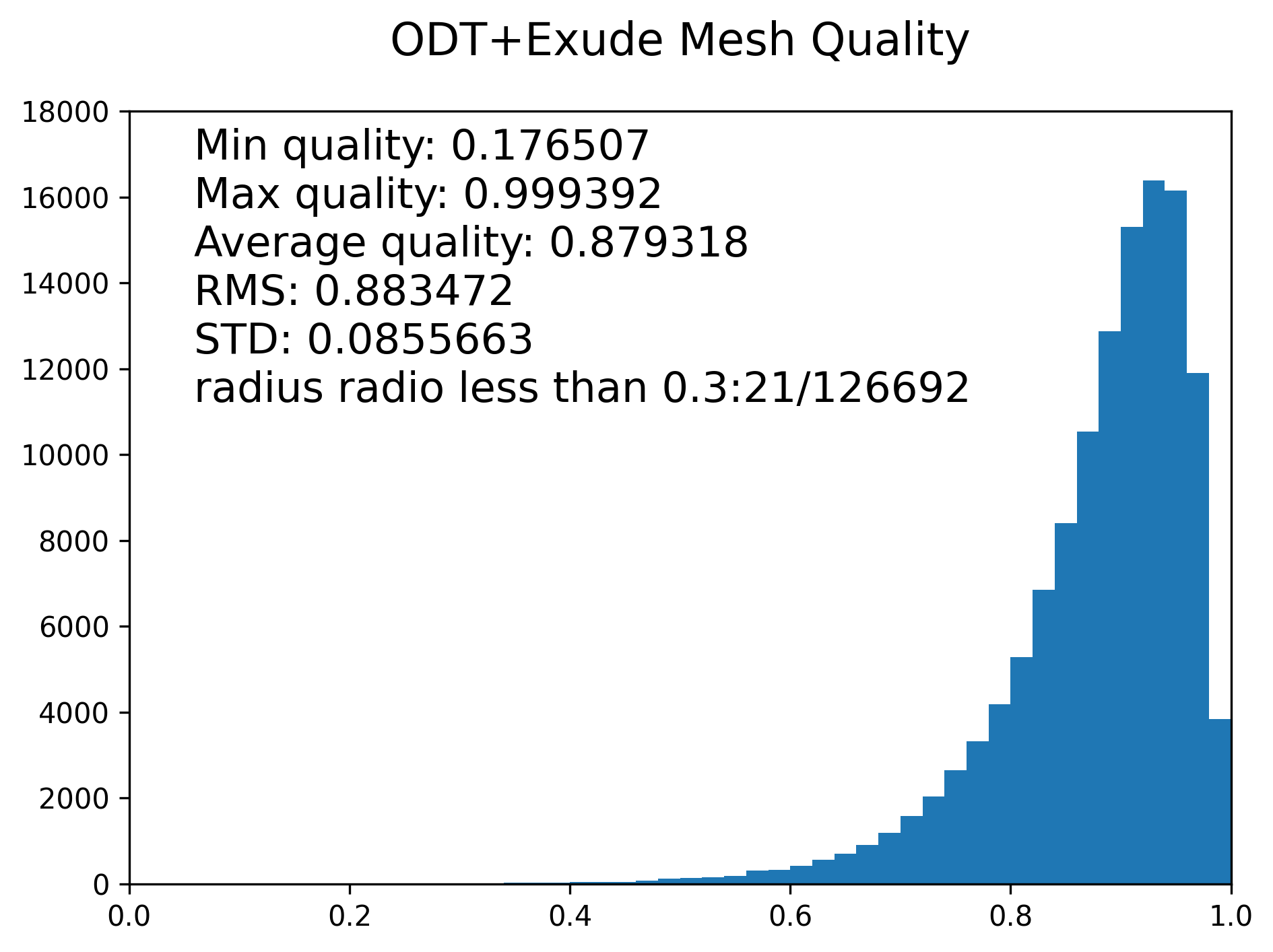}}
\vfill
\subfloat[Init min dihedral angle]{
\includegraphics[width=0.2\linewidth]{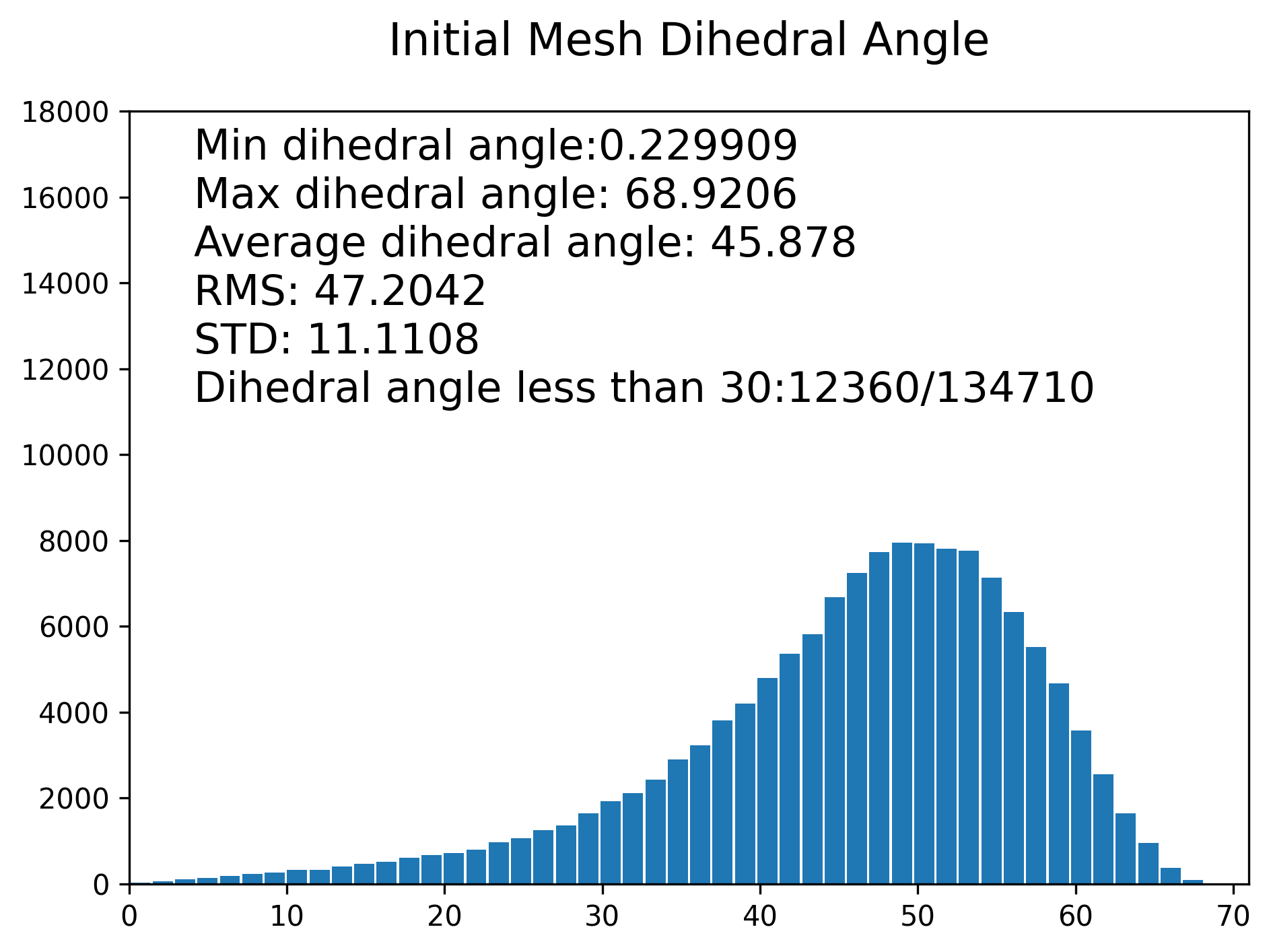}}
\subfloat[ODT min dihedral angle]{
\includegraphics[width=0.2\linewidth]{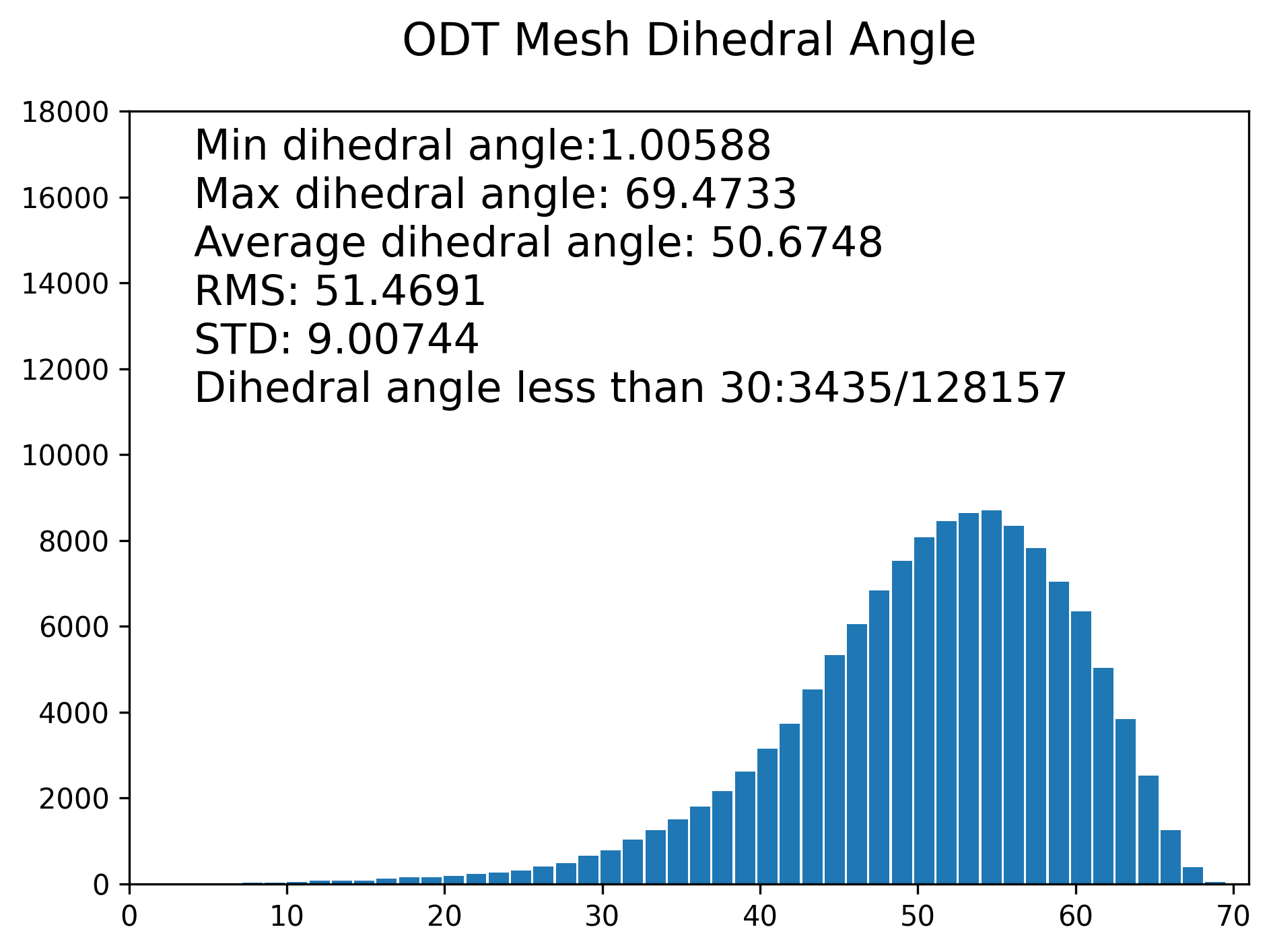}}
\hspace{0.01\linewidth}
\subfloat[RRE min dihedral angle]{
\includegraphics[width=0.2\linewidth]{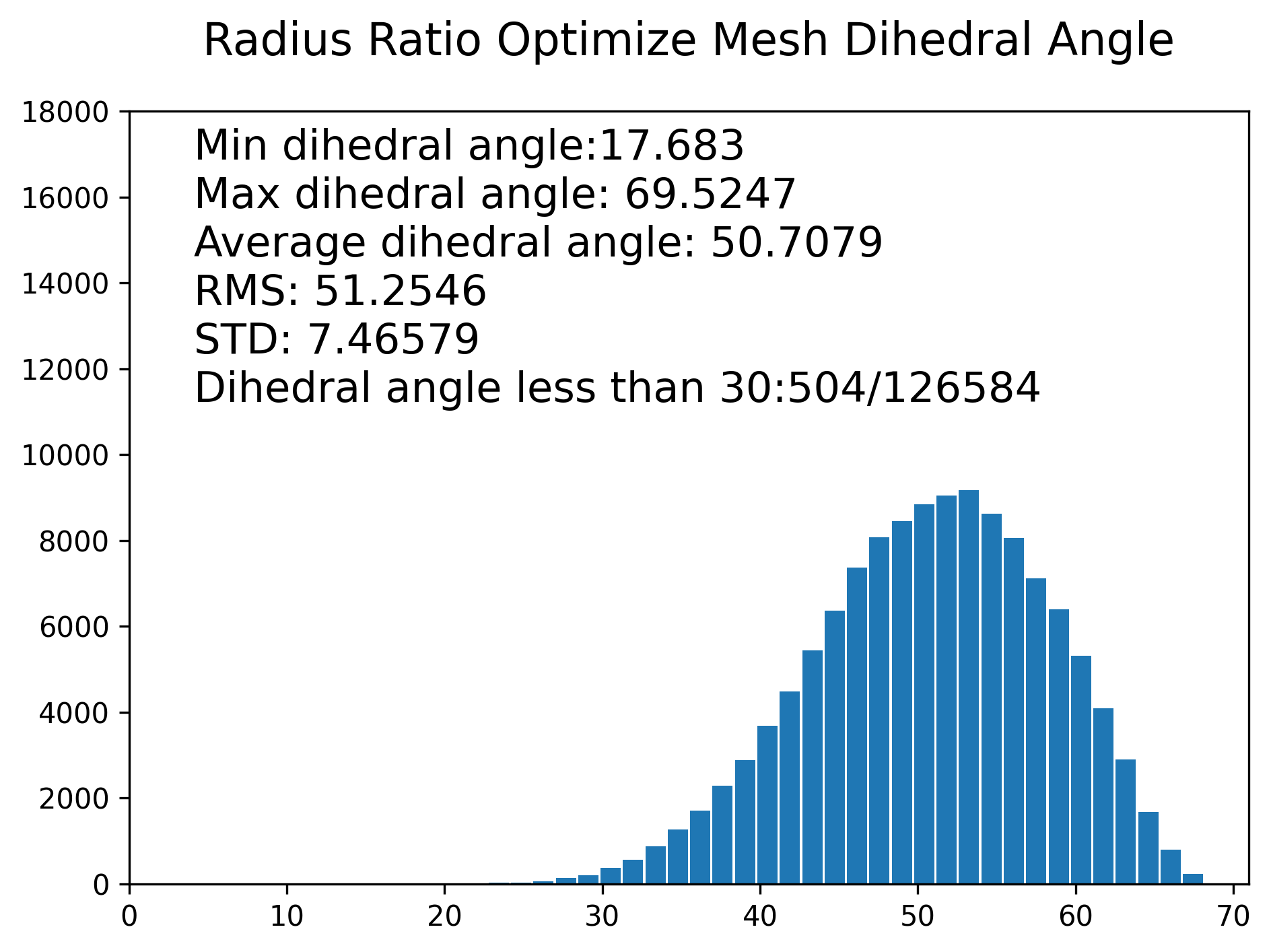}}
\hspace{0.01\linewidth}
\subfloat[Exude min dihedral angle]{
\includegraphics[width=0.2\linewidth]{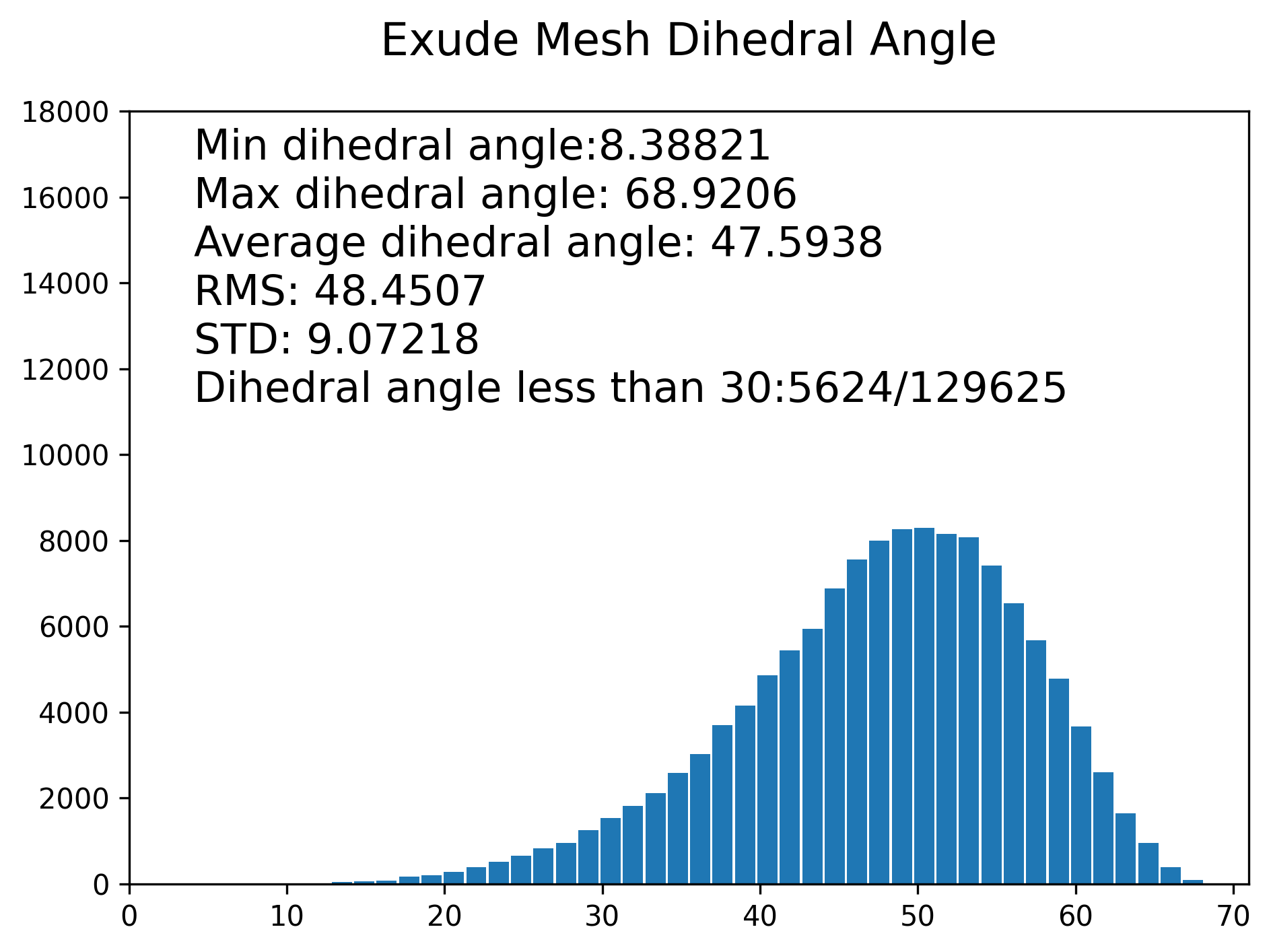}}
\caption{intersection of 12 spheres}
\label{fig:tsp}
\end{figure}
Fig.~\ref{fig:sphere},~\ref{fig:lshape},~\ref{fig:tsp} show three examples: a 
unit sphere, an L-shaped domain, and the intersection of 12 spheres. The initial 
mesh are generated by CGAL~\cite{45ref}. In addition to the RRE method, we compare with the 
ODT method and the Sliver Exude method provided by CGAL. The Exude method is 
designed to remove slivers. These comparisons allow a more complete 
evaluation of the RRE method in improving element quality and removing slivers.

Subfigures (c), (d), and (f) show the distributions of the radius ratio after 
RRE, Exude, and ODT, respectively. The results show that:
\begin{itemize}
\item ODT improves the overall mesh quality and shifts the distribution toward 
    higher values, but its effect on slivers is limited. A long 
    low-quality tail still remains.
\item Exude focuses on removing slivers. It improves the worst elements, 
    but has limited effect on the overall distribution.
\item RRE improves both the overall quality and the removal of slivers. It 
    increases the average quality and reduces slivers at the same time.
\end{itemize}

Subfigures (g) and (h) show the results of combining RRE and Exude with ODT. 
ODT+Exude improves the worst elements, but the overall distribution does not 
change much. In contrast, RRE applied after ODT can further remove remaining 
slivers and continue to improve global distribution.

In addition to the radius ratio, we also use the minimum dihedral angle as a 
secondary quality measure. From subfigures (i), (j), (k), and (l), we see that 
ODT gives limited improvement in the minimum dihedral angle, while Exude and 
RRE both improve it clearly. The improvement by RRE is more significant.

\begin{table}[htbp]
    \caption{Optimization Result Statistics}
    \label{table1}
    %\begin{adjustwidth}{-1.5cm}{0cm}
    \centering
    \begin{tabular}{|c|c|c|c|c|c|c|}
        \hline
        Model & Method & \makecell{RRE\\Iteration count/\\Flip count} & Num of
        cell & \makecell{Min dihedral \\angle(deg.)} & \makecell{Min radius
        \\ratio} & Time(sec.) \\
        \hline
        Sphere & Init & / & 18547  & 0.708 & 0.014 & / \\
        \hline
        & ODT & /  & 17537   & 2.751 & 0.056 & 16.96\\
        \hline
        & Exude &  / &  17794   & 11.18 & 0.190 & 8.63 \\
        \hline
        & RRE & 112/2 &  17449   & 21.87 & 0.446 & 7.88\\
        \hline
        & Precondition RRE & 50/2 & 17449 & 21.87 & 0.445 & 4.64 \\
        \hline
        & ODT+RRE & 68/1 & 17236 & 27.76 & 0.566 & 21.43 \\
        \hline
        & ODT+Exude & / & 17335 & 15.13 & 0.325 & 23.35 \\
        \hline
        Lshape & Init & / & 62002 & 0.516 & 0.011 & / \\
        \hline
        & ODT & / & 58873 & 3.332 & 0.063 & 59.56 \\
        \hline
        & Exude & / & 59677 & 9.207 & 0.202 & 31.20 \\
        \hline
        & RRE & 222/3 & 58342 & 19.59 & 0.407 & 50.98 \\
        \hline
        & Precondition RRE & 131/3 & 58345 & 18.87 & 0.395 & 33.66 \\
        \hline
        & ODT+RRE & 93/2 & 57927 & 18.76 & 0.407 & 81.34\\
        \hline
        & ODT+Exude & / & 58243 & 9.25 & 0.202 & 81.80 \\
        \hline
        \makecell{Intersection of \\12 spheres} &  Init & / & 134710 & 0.230 & 9.63e-4& /\\
        \hline
            & ODT & / & 128157 & 1.006 & 0.002 & 123.71 \\
        \hline
          & Exude & / & 129625 & 8.388 & 0.169 & 69.24 \\
        \hline
        & RRE & 225/3 & 126584 & 17.68 & 0.331 & 116.43 \\
        \hline
        & Precondition RRE & 131/3 & 126584 & 16.95 & 0.322 & 100.10 \\
        \hline
        & ODT+RRE & 118/2 & 125743 & 20.27 & 0.436 & 183.17\\
        \hline
        & ODT+Exude & / & 126692 & 7.57 & 0.177 & 174.94\\
        \hline
    \end{tabular}
%\end{adjustwidth}
\end{table}
Table~\ref{table1} gives a further comparison in terms of quality improvement 
and efficiency. In all three examples, RRE gives higher minimum radius ratio 
and minimum dihedral angle than the other methods. In terms of efficiency, the 
running time of RRE is close to that of ODT. The Exude method is faster, but it 
focuses only on sliver removal. RRE improves both global quality and sliver 
removal, and its effect is close to the combined ODT+Exude method, with good 
overall efficiency.

Table~\ref{table1} also compares the results with and without the preconditioner. 
As in the two-dimensional case, the “iteration count” refers to the number of 
calls to the energy function during the Wolfe line search. This measure better 
reflects the computational cost. With preconditioner, the number of iterations 
is clearly reduced, and the running time is also shorter. This shows that the 
preconditioner improves convergence. The reduction in time is smaller than the 
reduction in iteration count, because the problem size is small and the cost of 
building the preconditioner is relatively high. It can also be seen that the 
final results are similar with and without preconditioner. This means that 
the preconditioner does not change the solution itself. Its role is to improve 
the search direction and accelerate convergence. Overall, the preconditioner 
improves efficiency while keeping the same optimization result.

\begin{figure}[htbp]
\centering
\subfloat[Model]{
\includegraphics[width=0.2\linewidth]{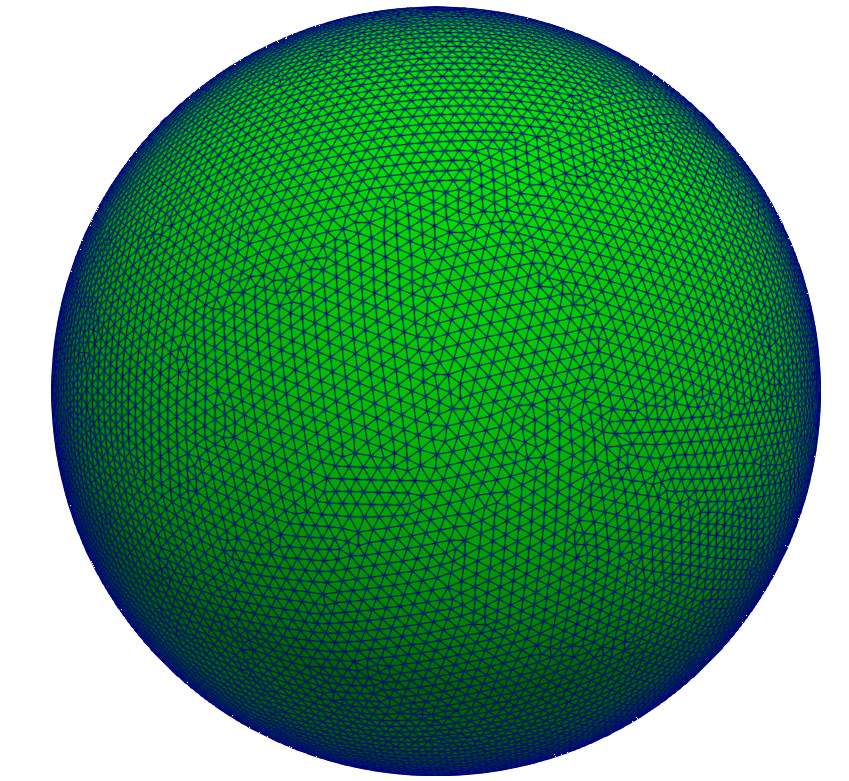}}
\hspace{0.01\linewidth}
\subfloat[Init min dihedral angle]{
\includegraphics[width=0.2\linewidth]{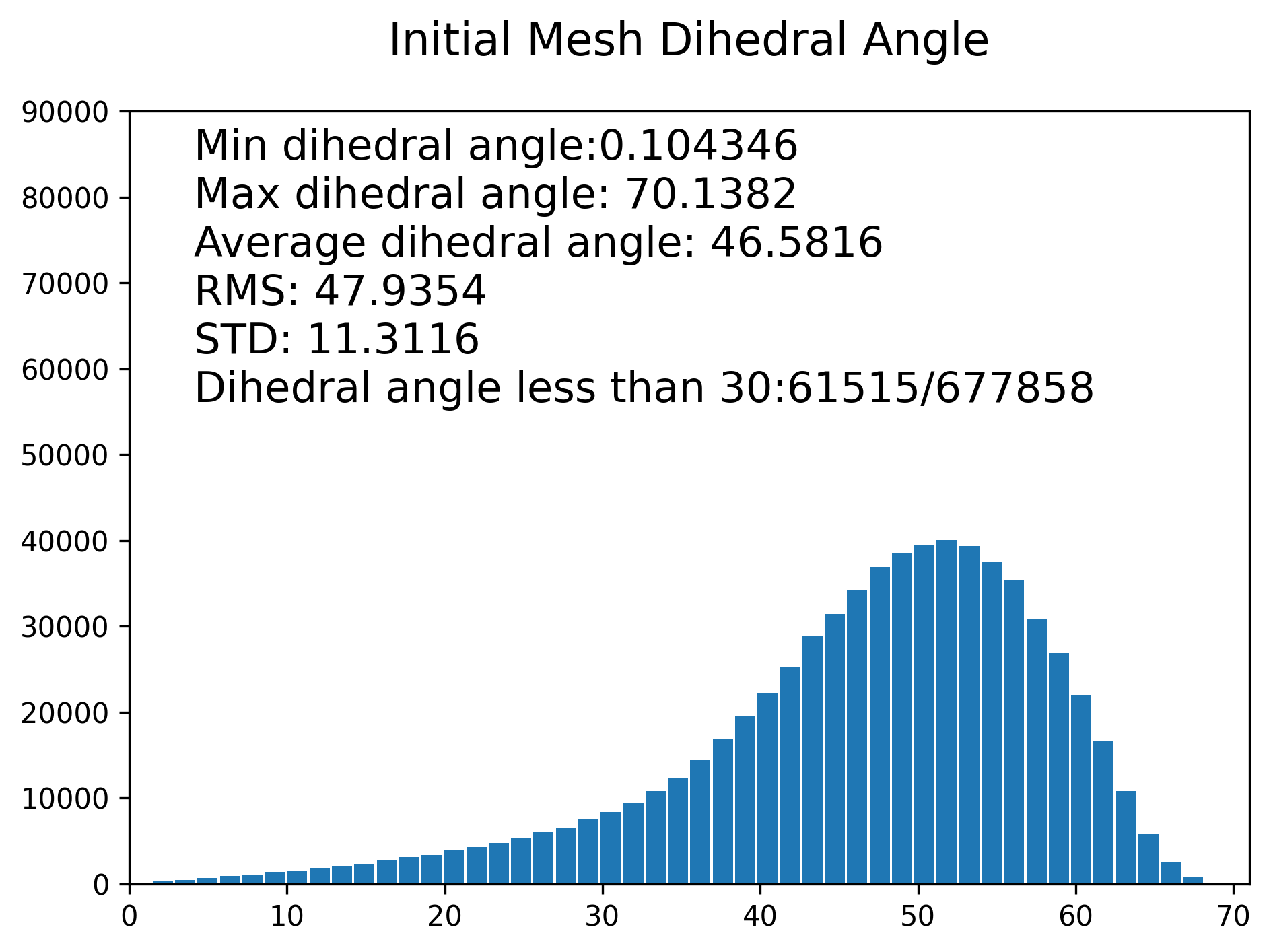}}
\hspace{0.01\linewidth}
\subfloat[RRE min dihedral angle]{
\includegraphics[width=0.2\linewidth]{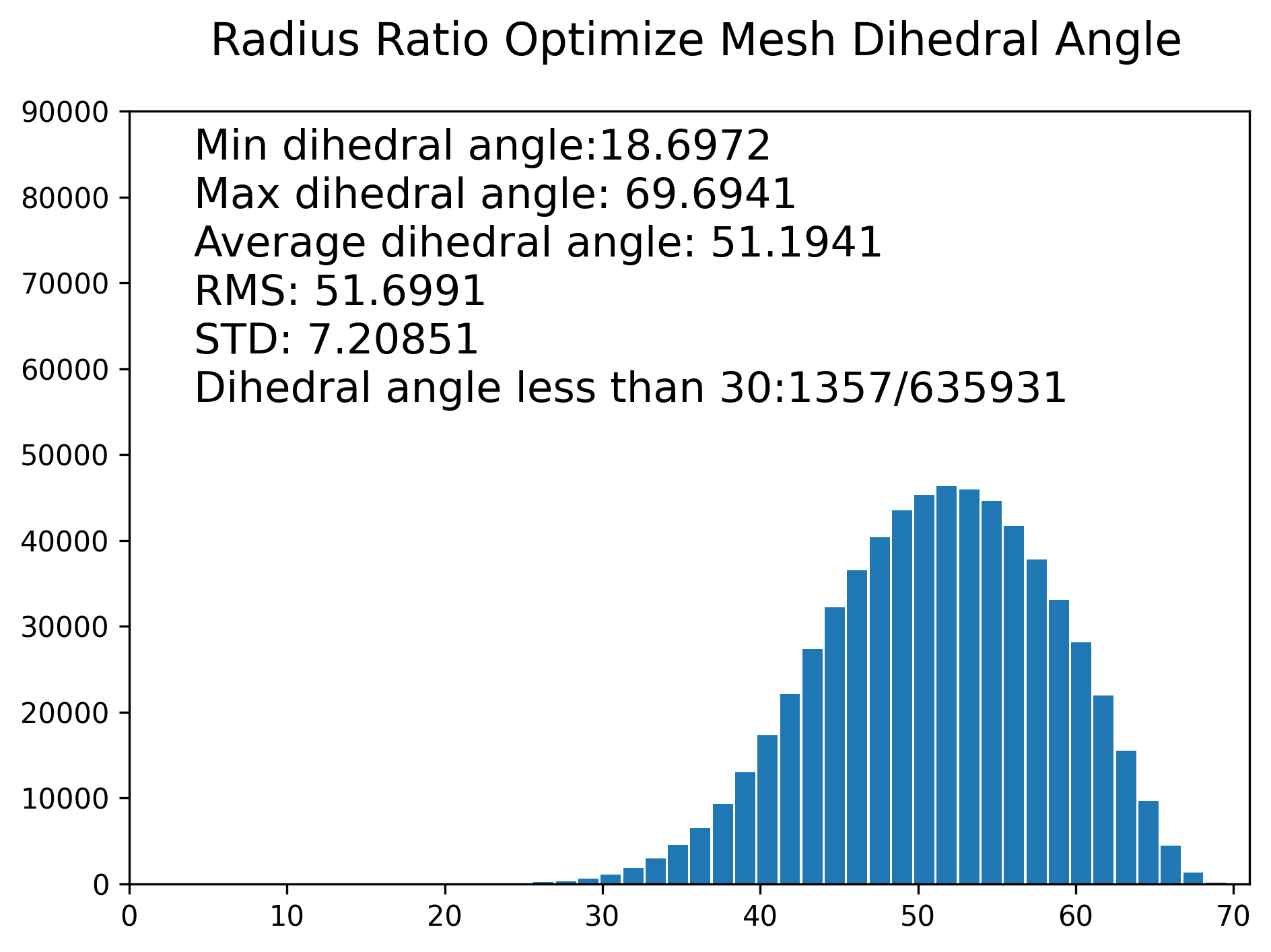}}
\hspace{0.01\linewidth}
\subfloat[Precondition RRE min dihedral angle]{
\includegraphics[width=0.2\linewidth]{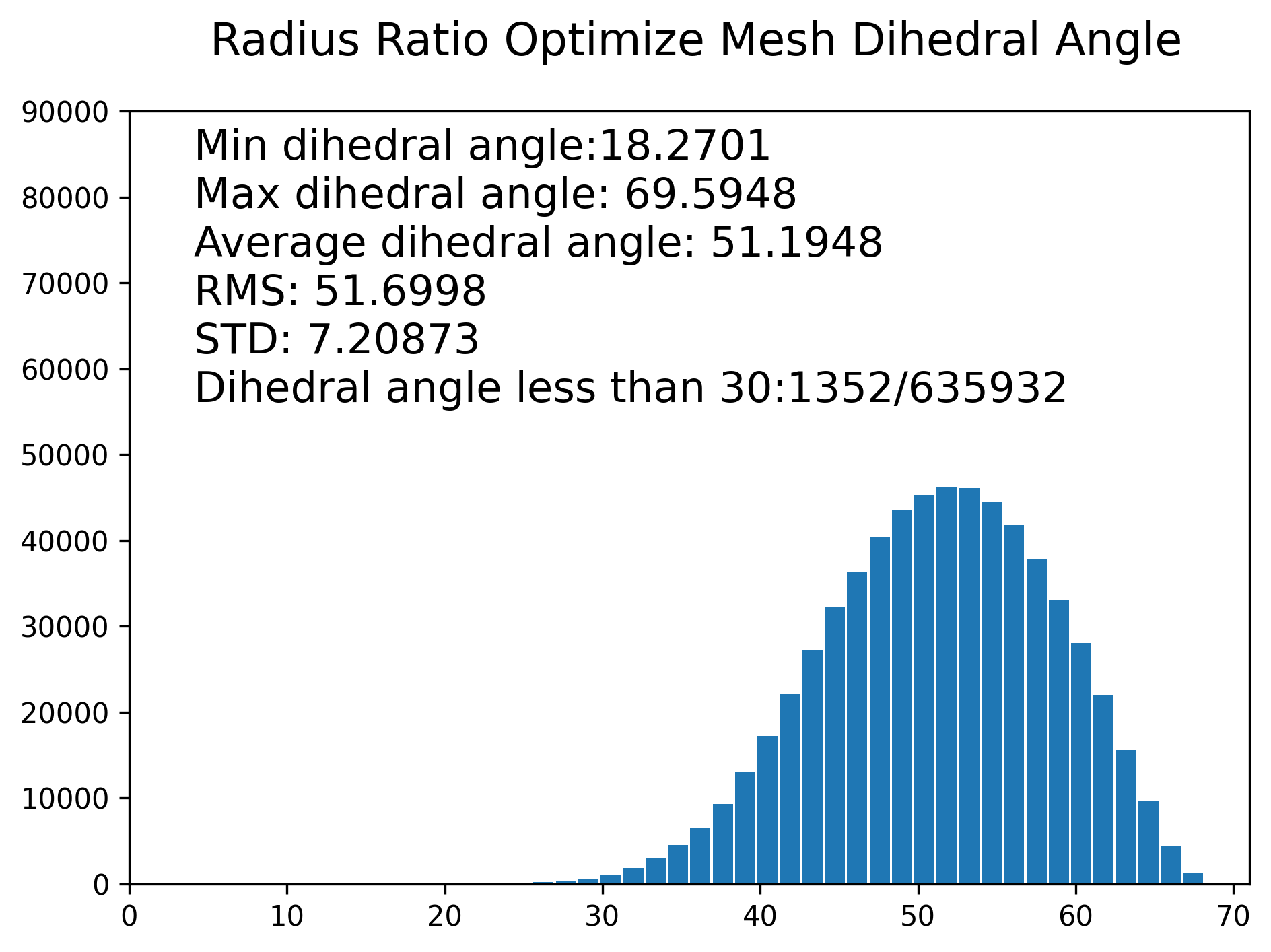}}
\vfill
\subfloat[Init quality]{
\includegraphics[width=0.3\linewidth]{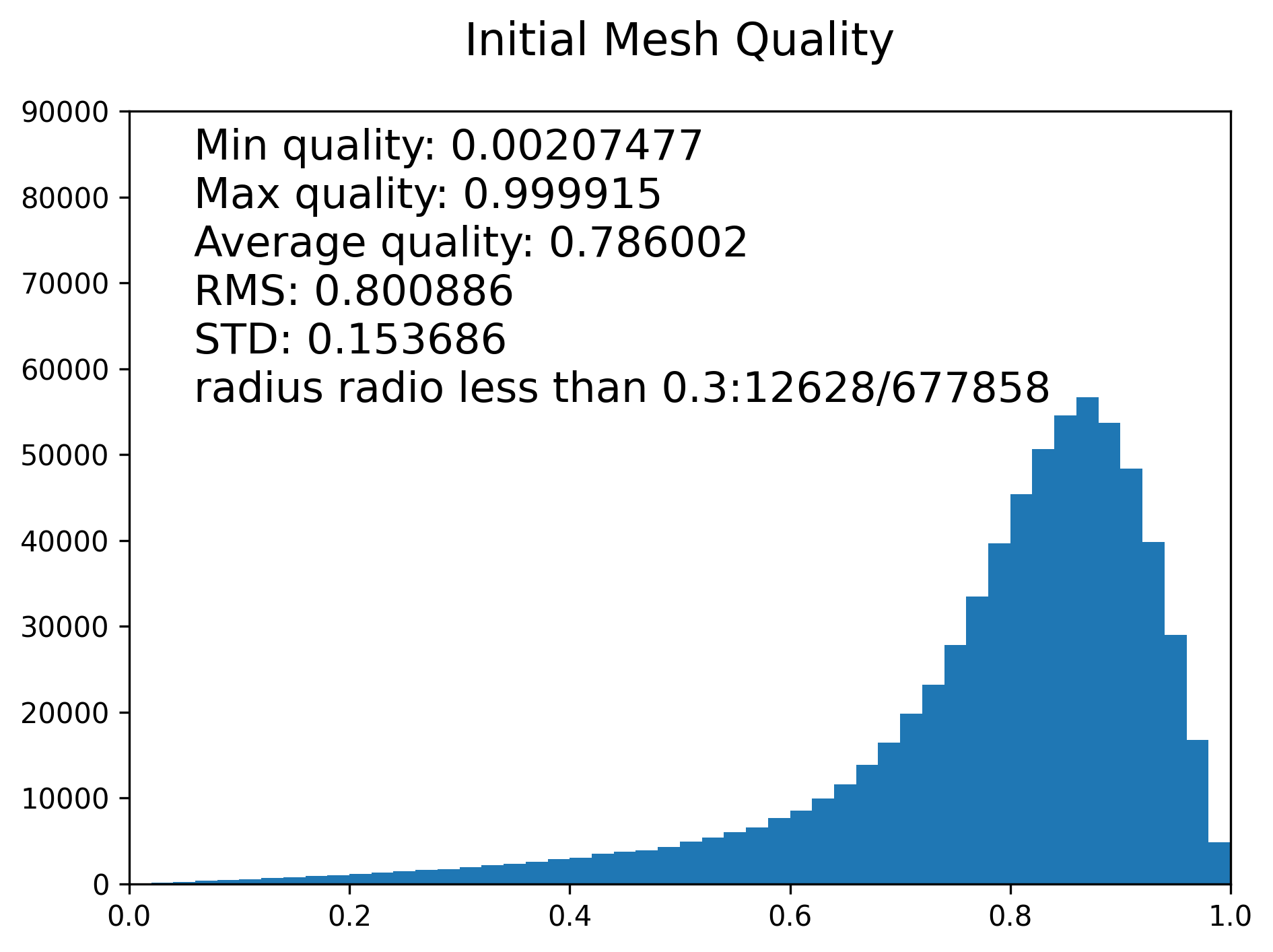}}
\hspace{0.01\linewidth}
\subfloat[RRE]{
\includegraphics[width=0.3\linewidth]{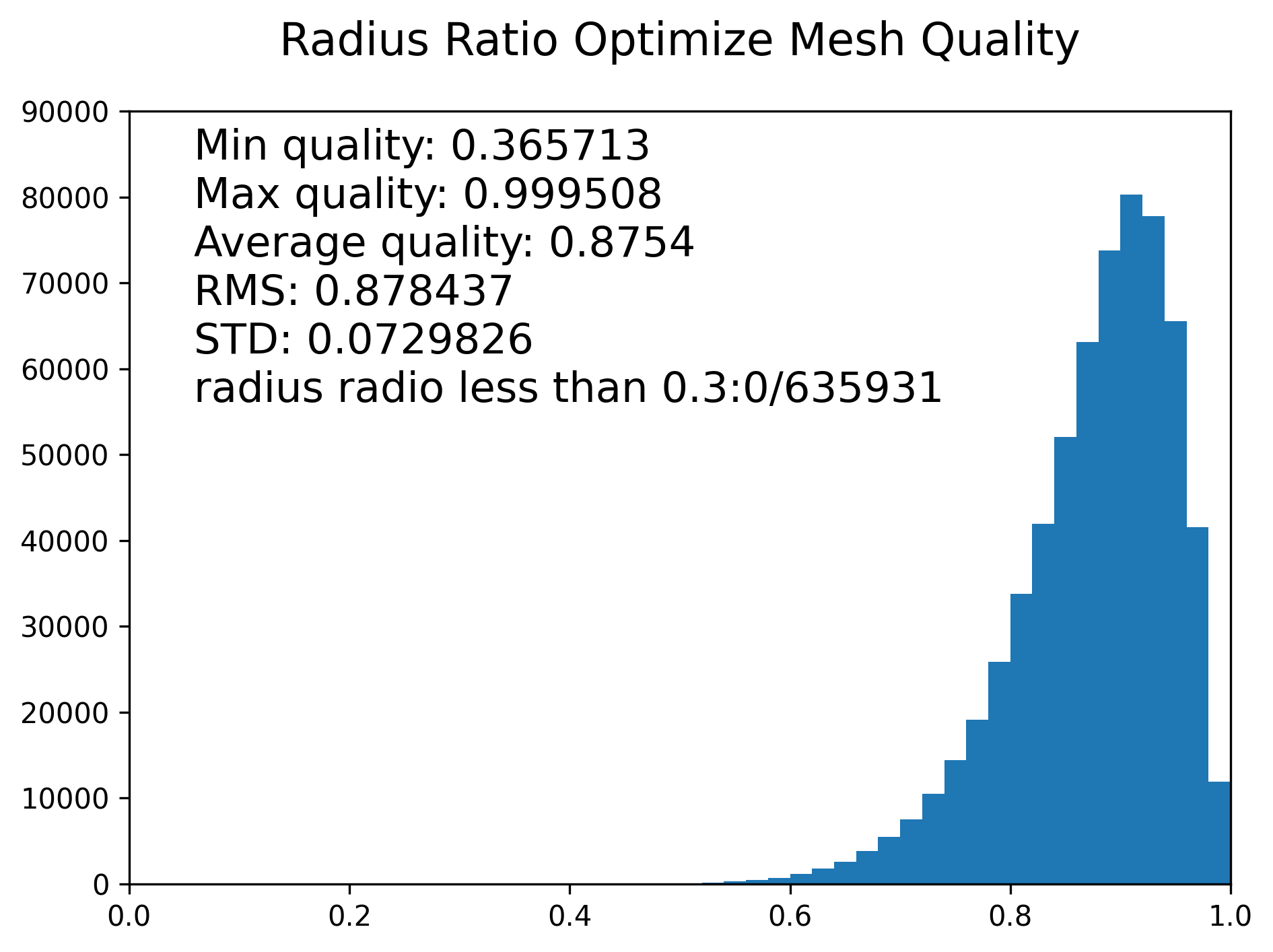}}
\hspace{0.01\linewidth}
\subfloat[Precondition RRE]{
\includegraphics[width=0.3\linewidth]{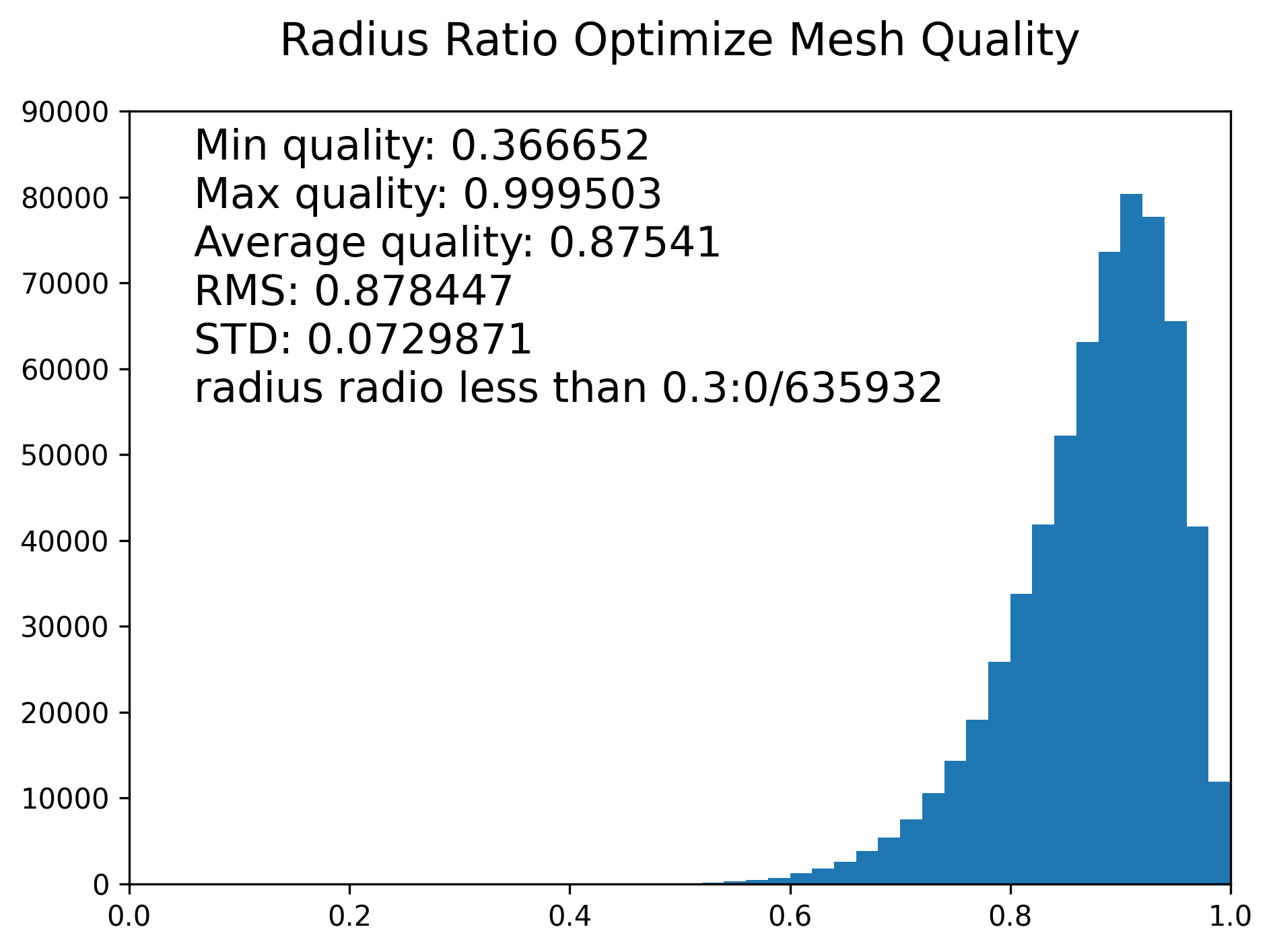}}
\caption{Sphere large scale mesh}
\label{fig:spherebig}
\end{figure}

\begin{figure}[htbp]
\centering
\subfloat[Model]{
\includegraphics[width=0.2\linewidth]{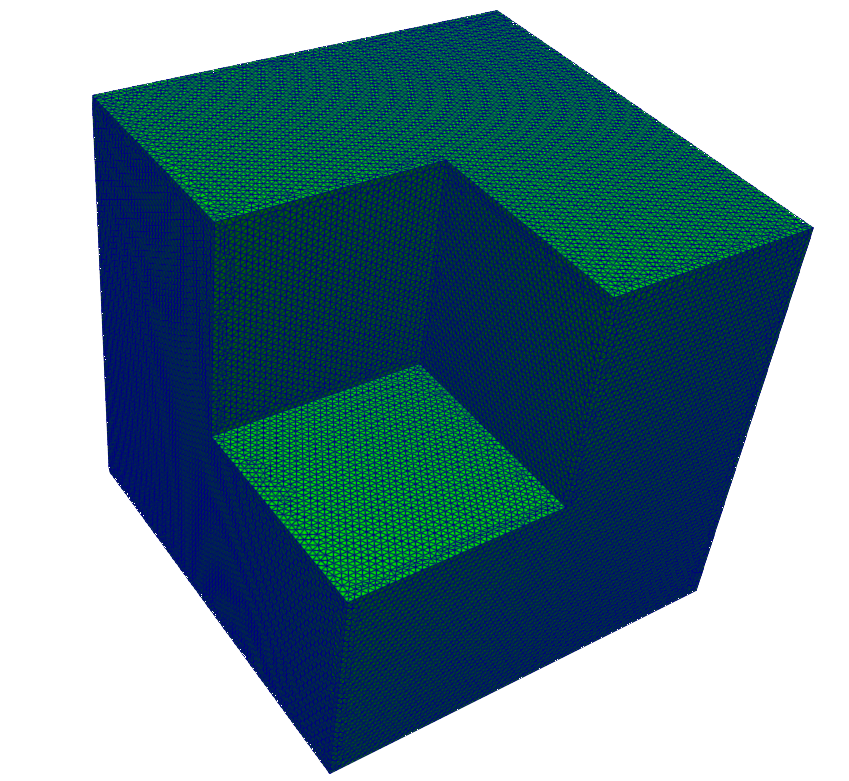}}
\hspace{0.01\linewidth}
\subfloat[Init min dihedral angle]{
\includegraphics[width=0.2\linewidth]{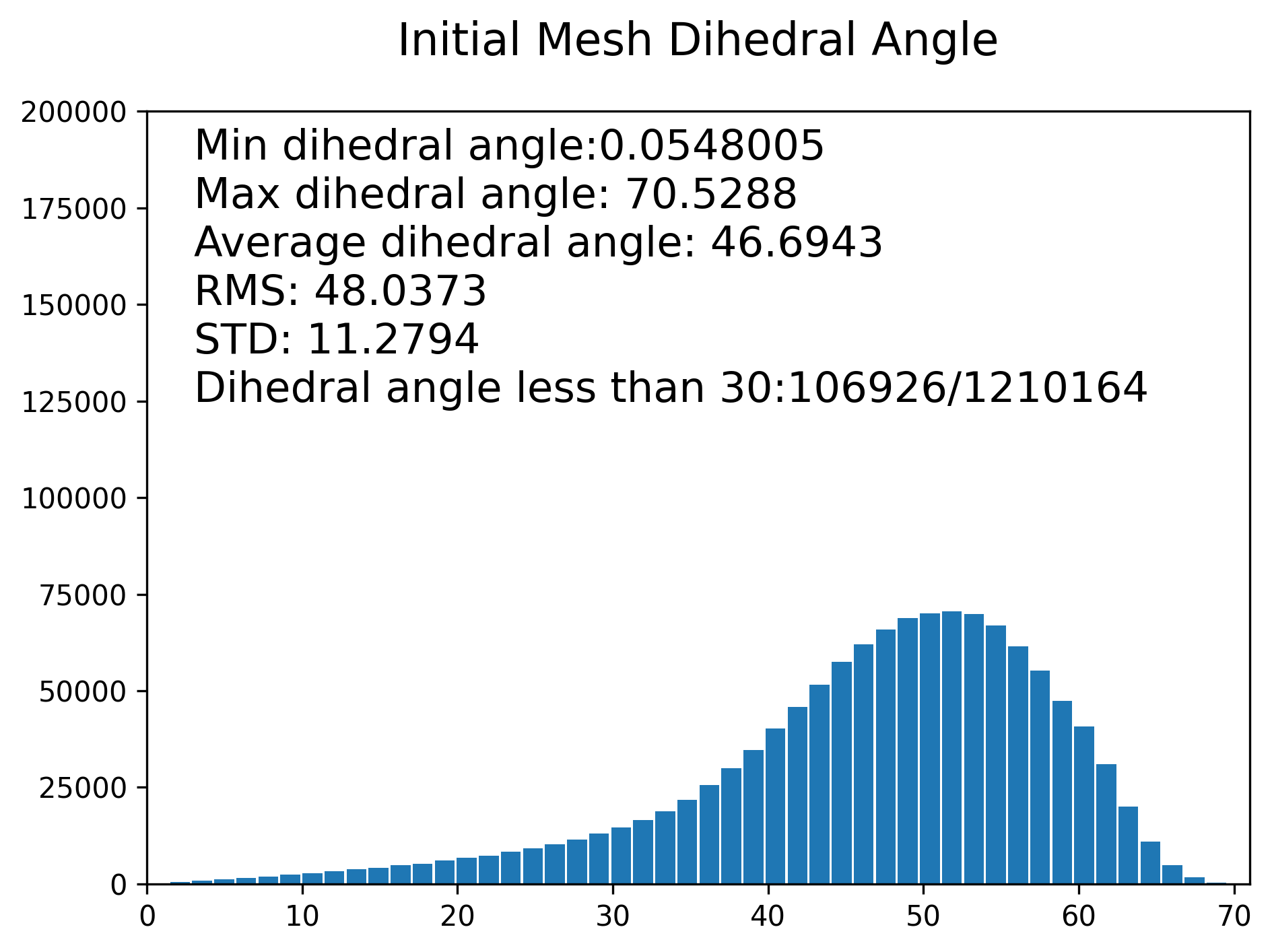}}
\hspace{0.01\linewidth}
\subfloat[RRE min dihedral angle]{
\includegraphics[width=0.2\linewidth]{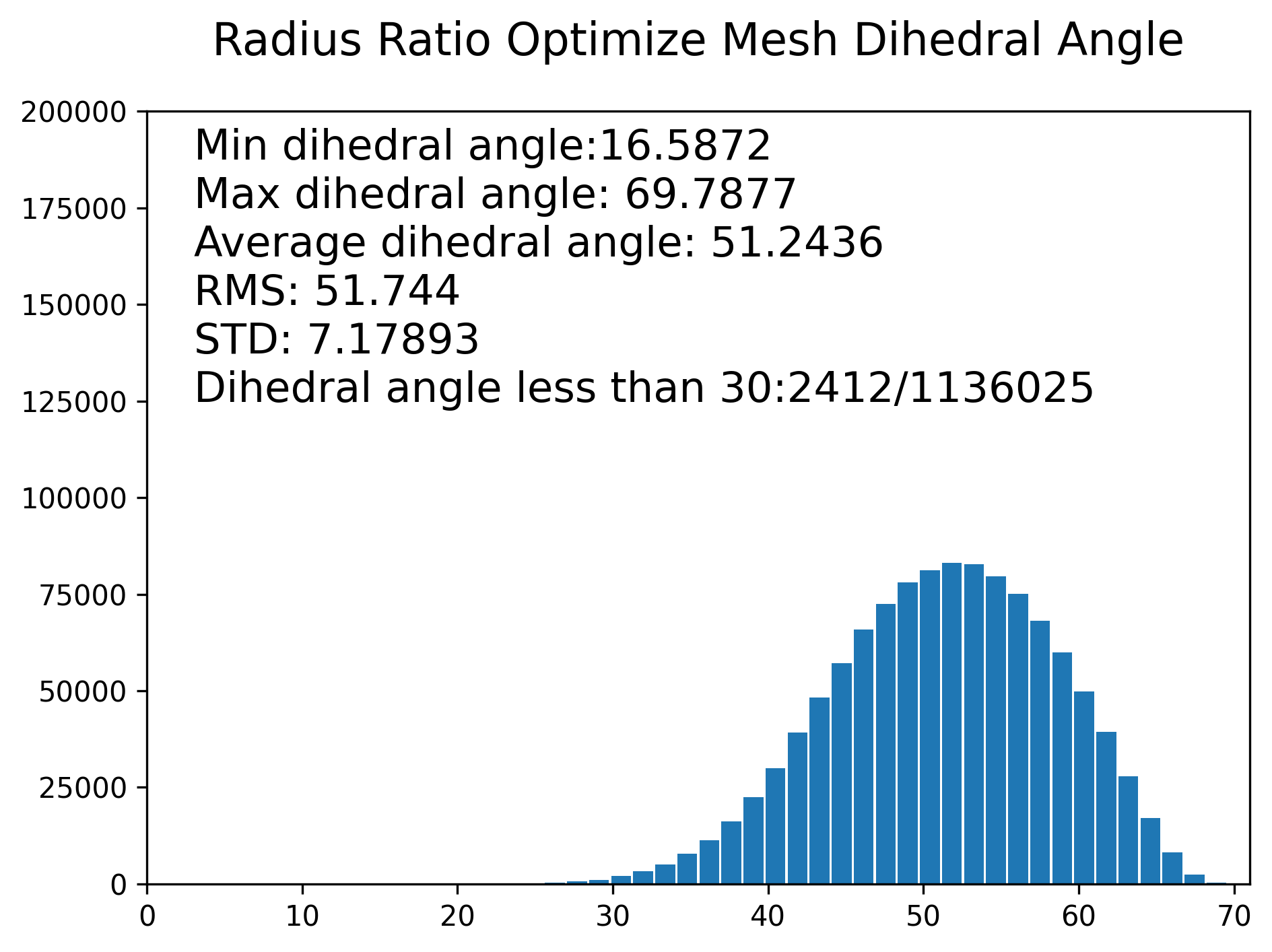}}
\hspace{0.01\linewidth}
\subfloat[Precondition RRE min dihedral angle]{
\includegraphics[width=0.2\linewidth]{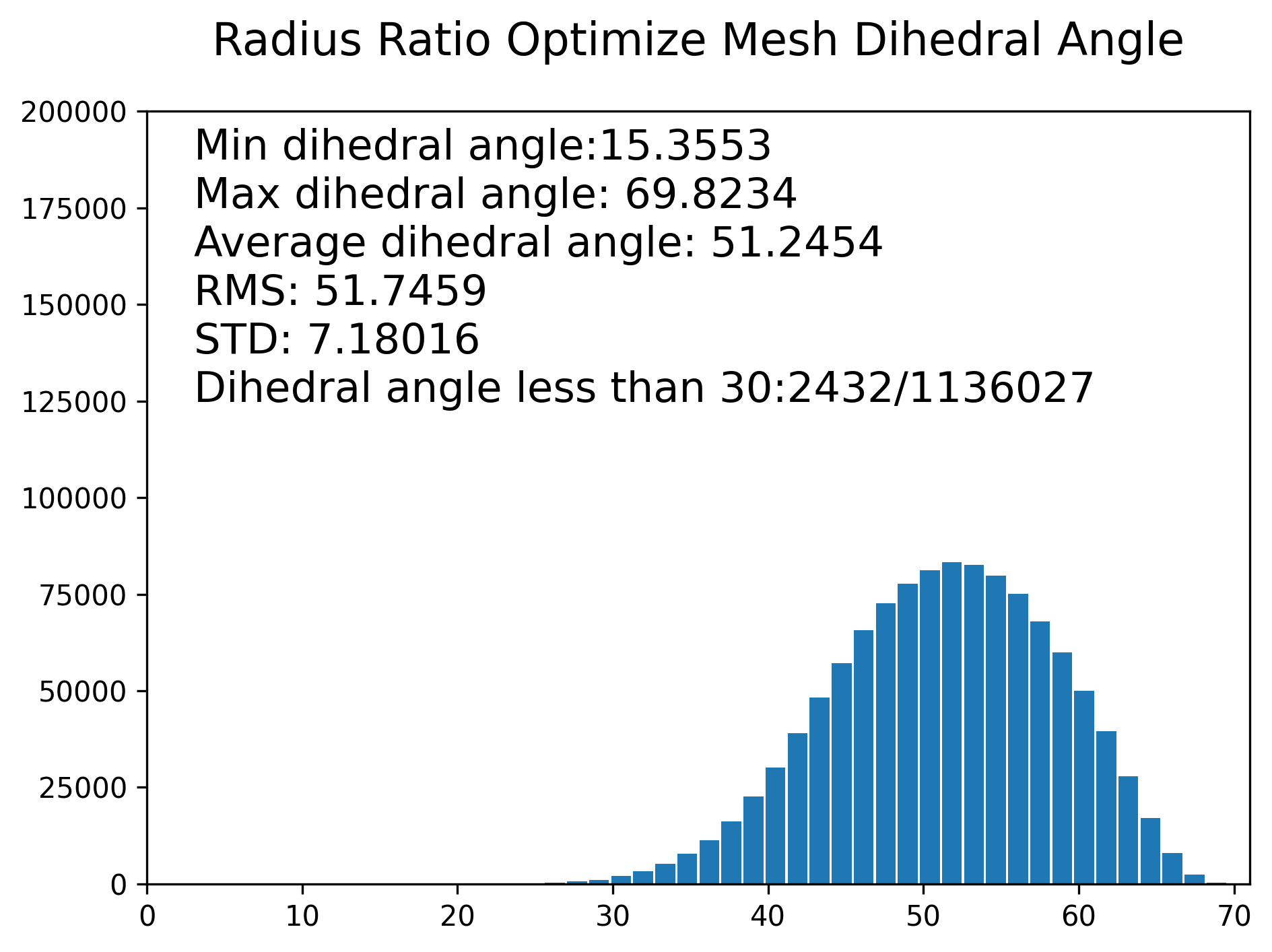}}
\vfill
\subfloat[Init quality]{
\includegraphics[width=0.3\linewidth]{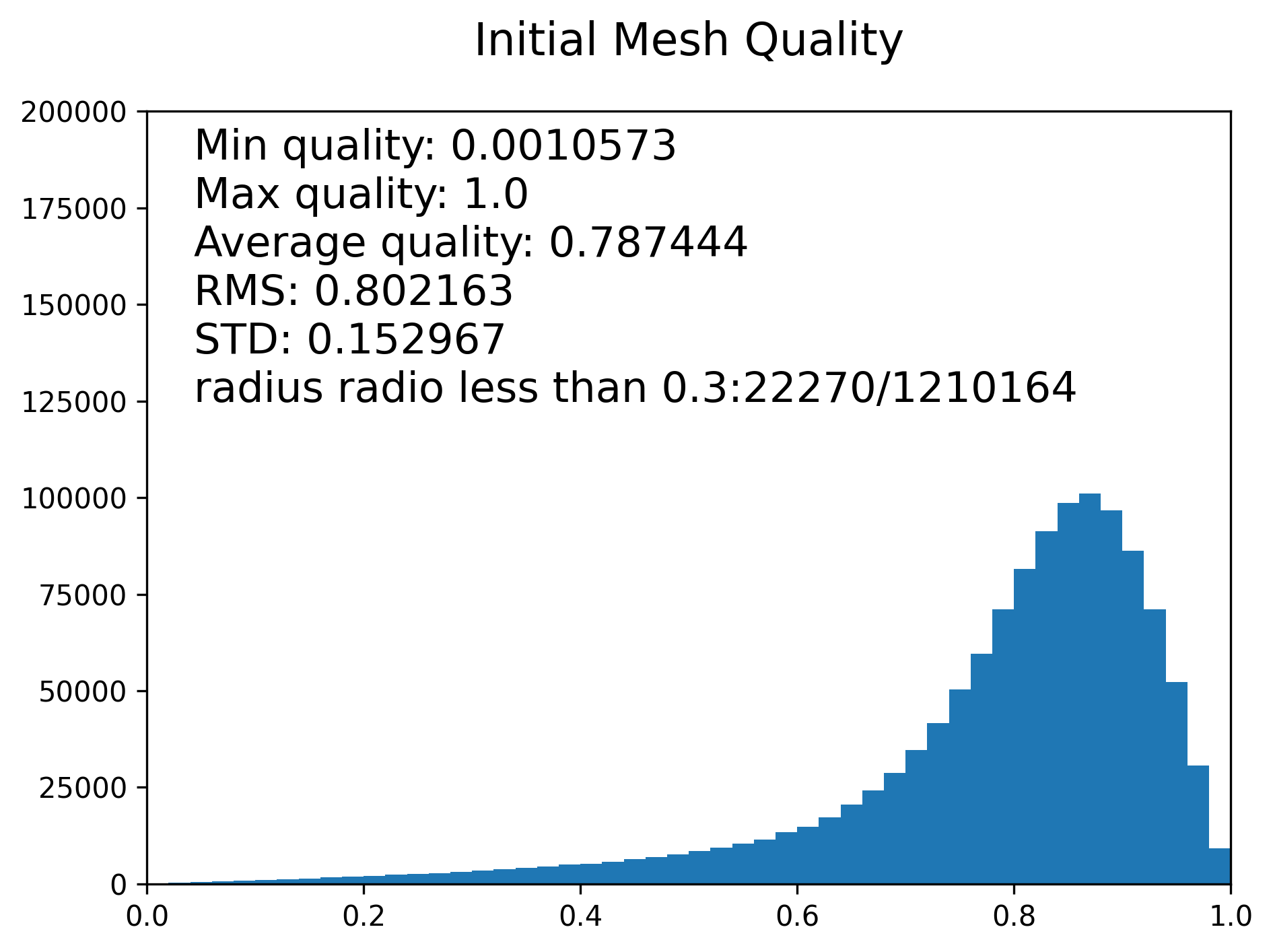}}
\hspace{0.01\linewidth}
\subfloat[RRE]{
\includegraphics[width=0.3\linewidth]{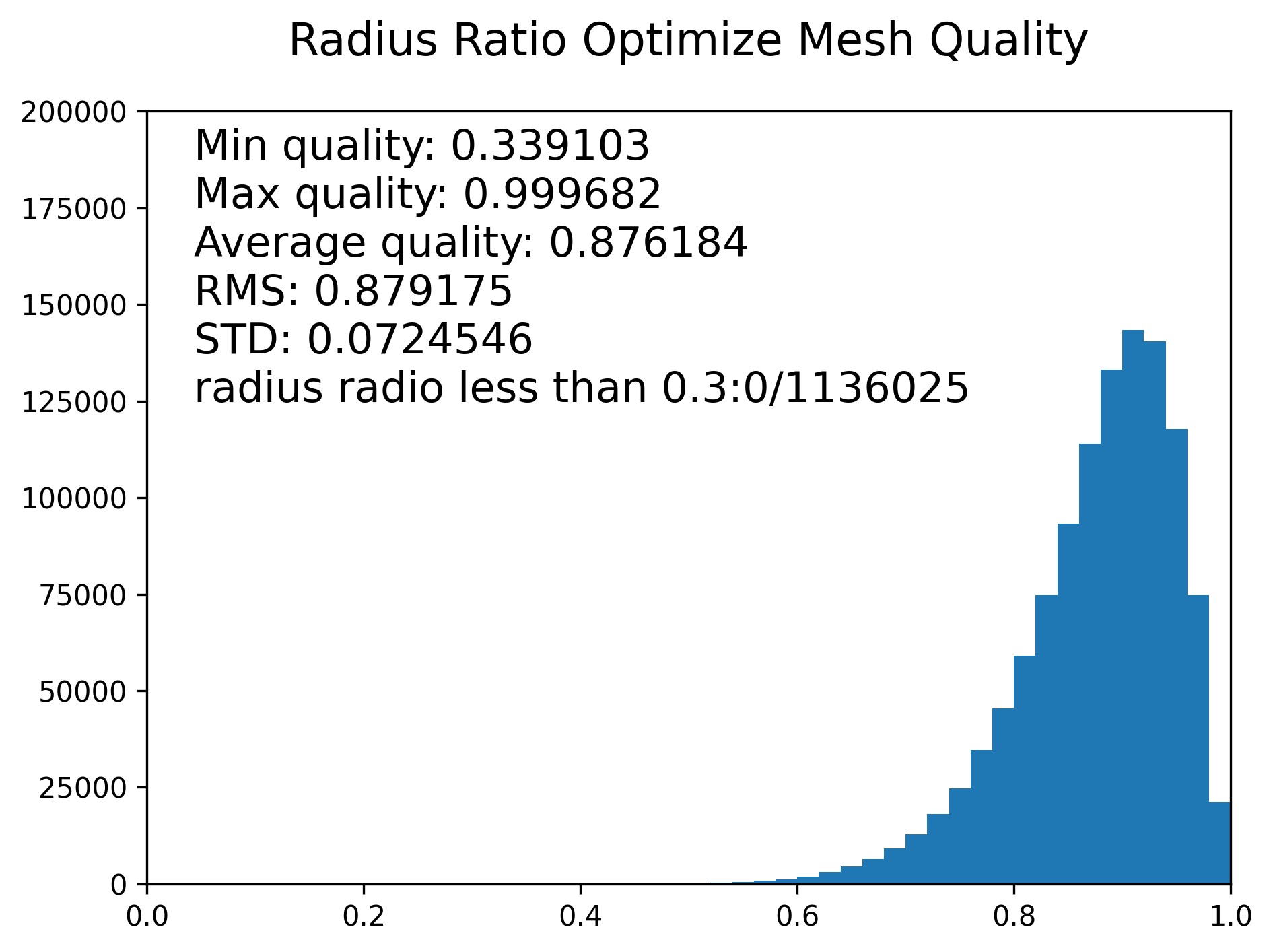}}
\hspace{0.01\linewidth}
\subfloat[Precondition RRE]{
\includegraphics[width=0.3\linewidth]{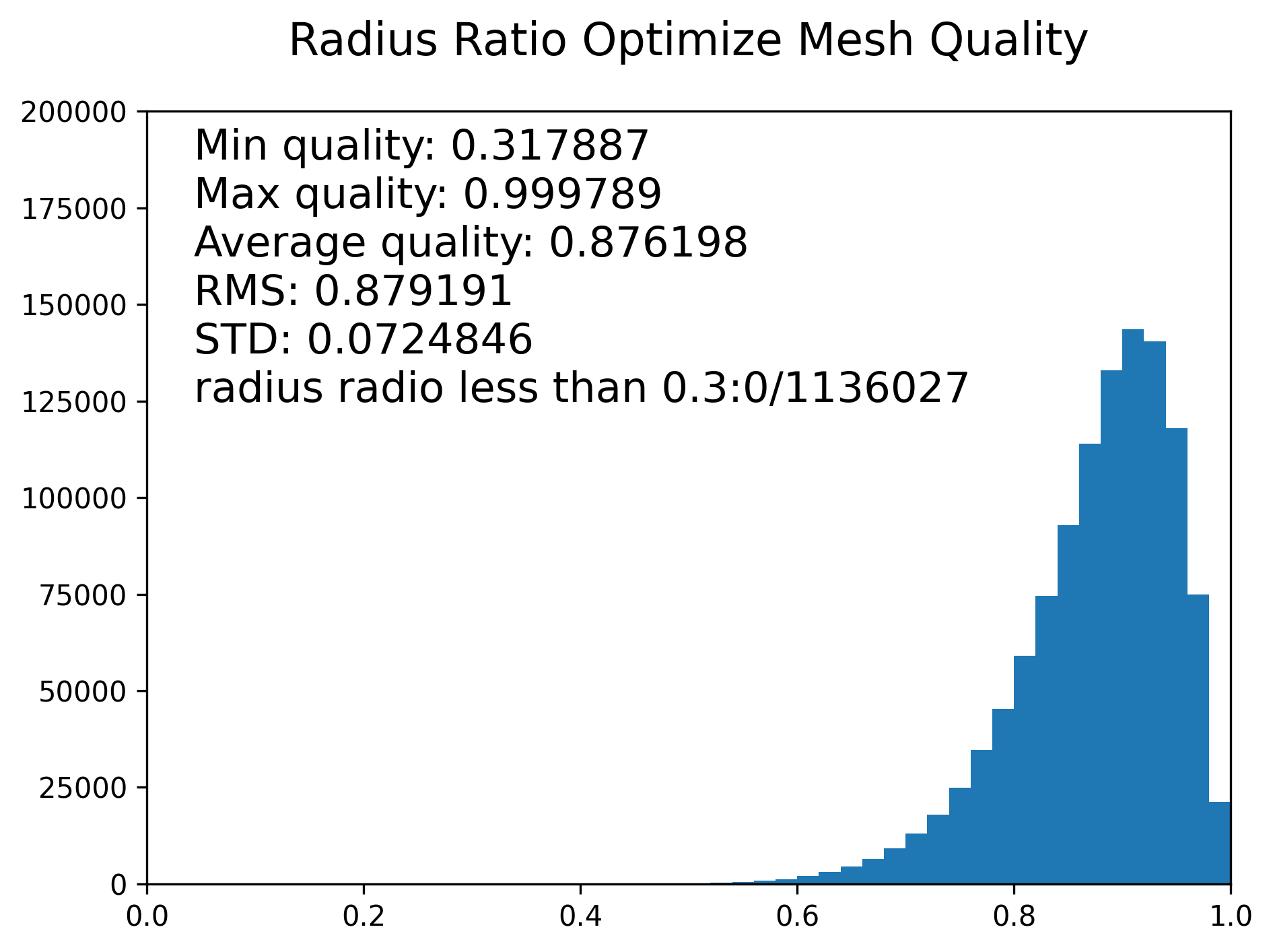}}
\caption{L-shape large scale mesh}
\label{fig:lshapebig}
\end{figure}

\begin{figure}[htbp]
\centering
\subfloat[Model]{
\includegraphics[width=0.2\linewidth]{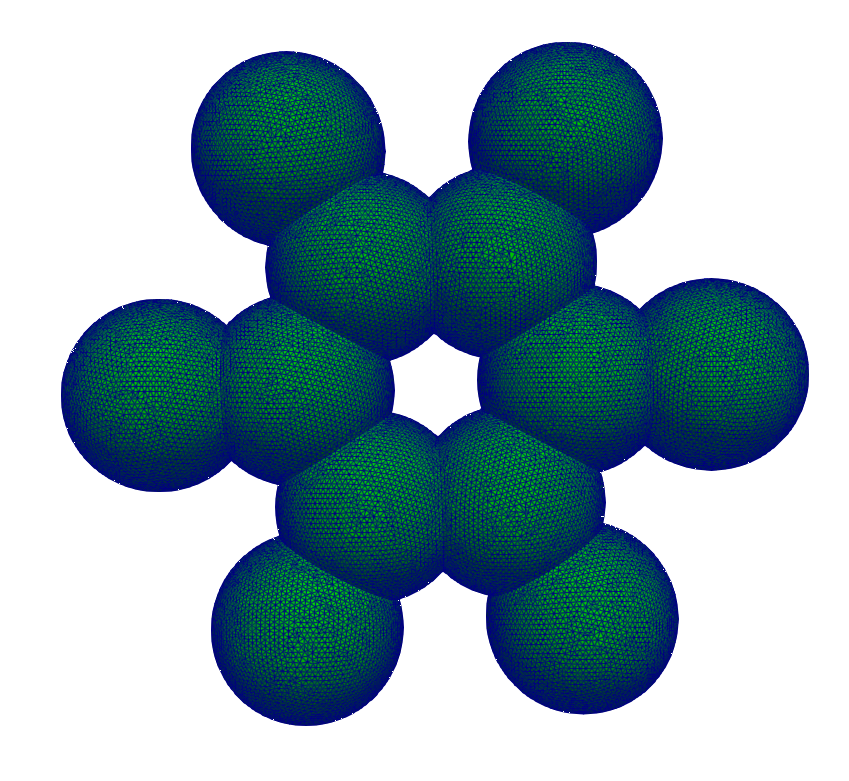}}
\hspace{0.01\linewidth}
\subfloat[Init min dihedral angle]{
\includegraphics[width=0.2\linewidth]{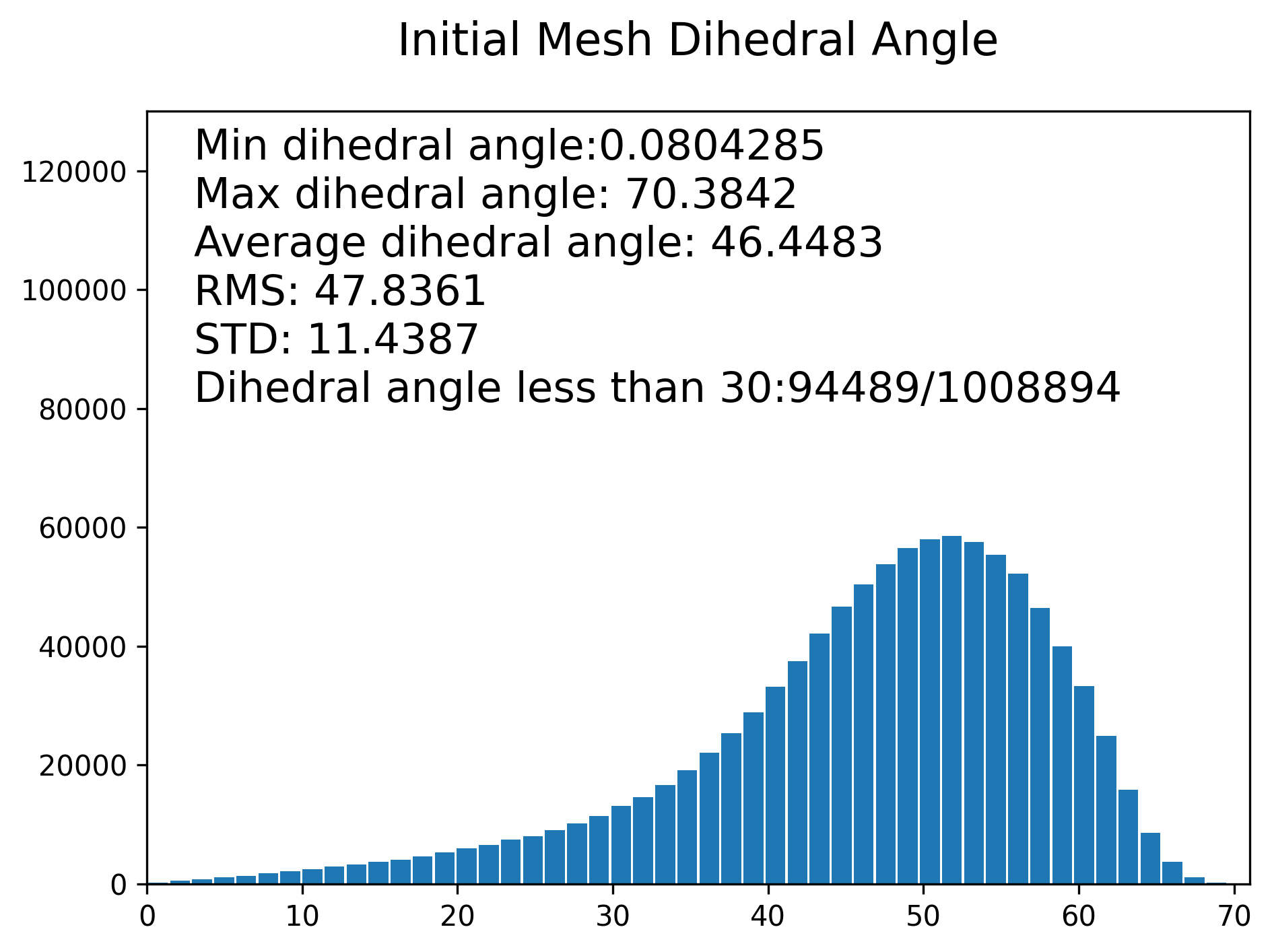}}
\subfloat[RRE min dihedral angle]{
\includegraphics[width=0.2\linewidth]{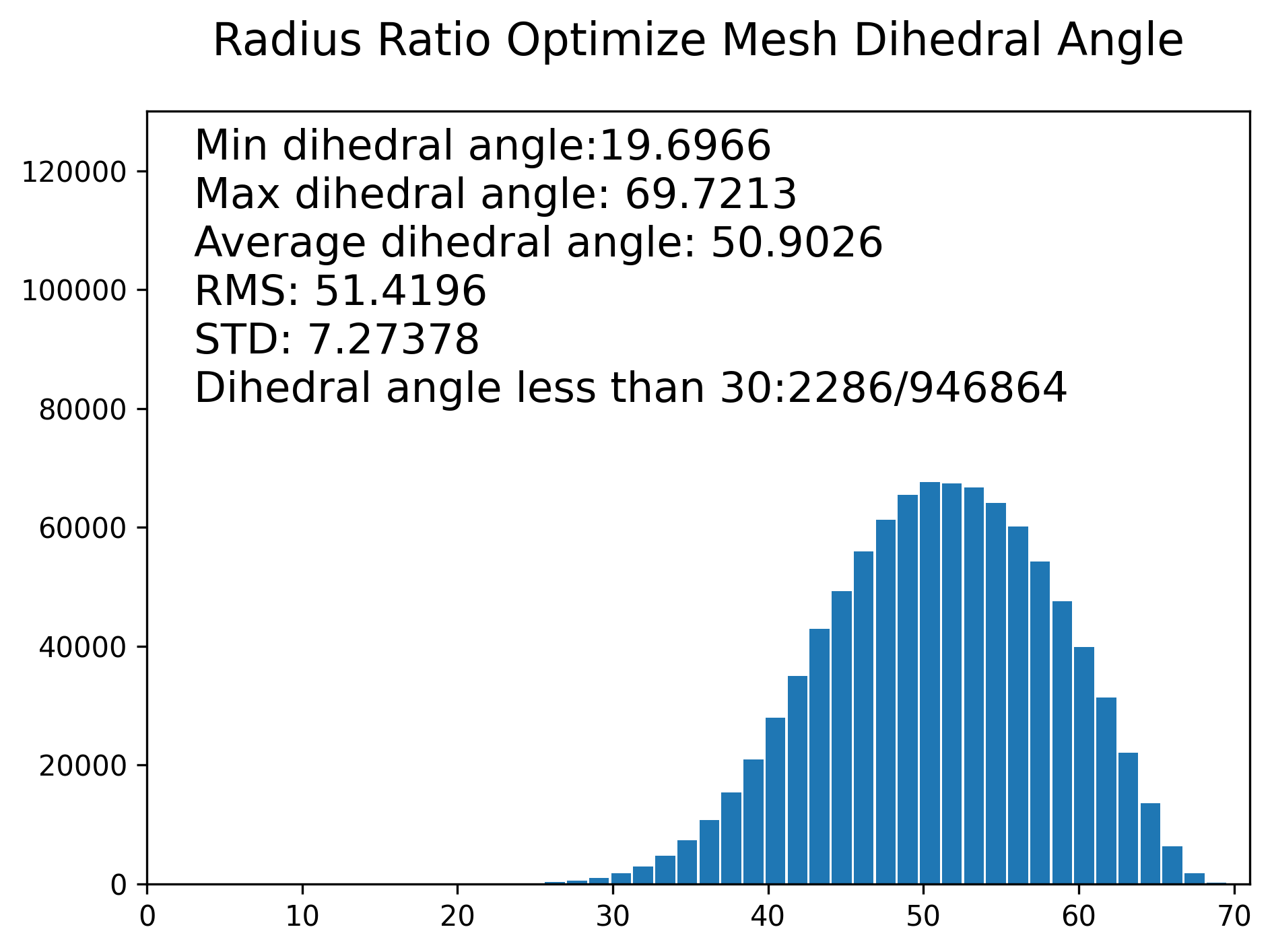}}
\hspace{0.01\linewidth}
\subfloat[Precondition RRE min dihedral angle]{
\includegraphics[width=0.2\linewidth]{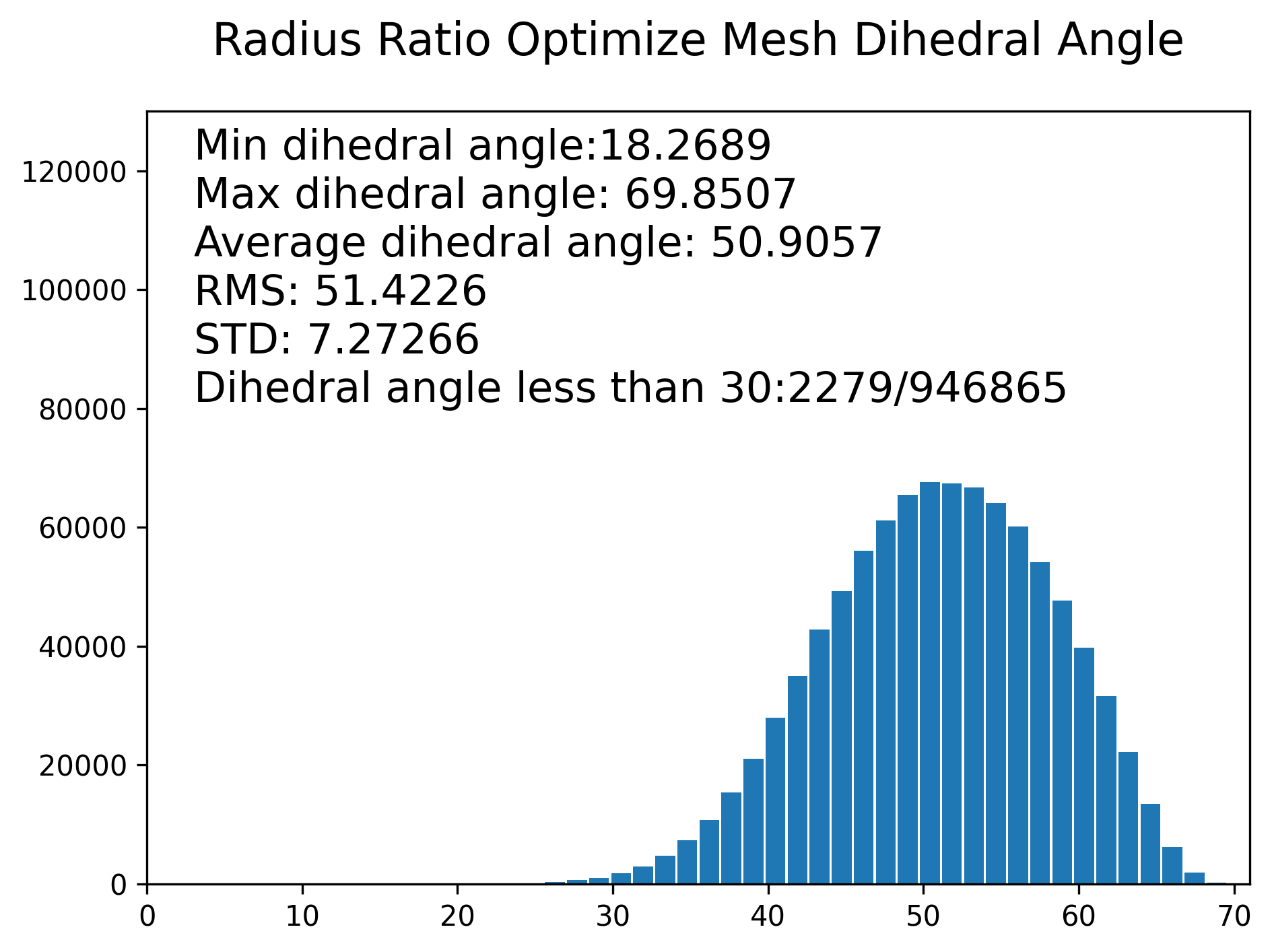}}
\vfill
\subfloat[Init quality]{
\includegraphics[width=0.3\linewidth]{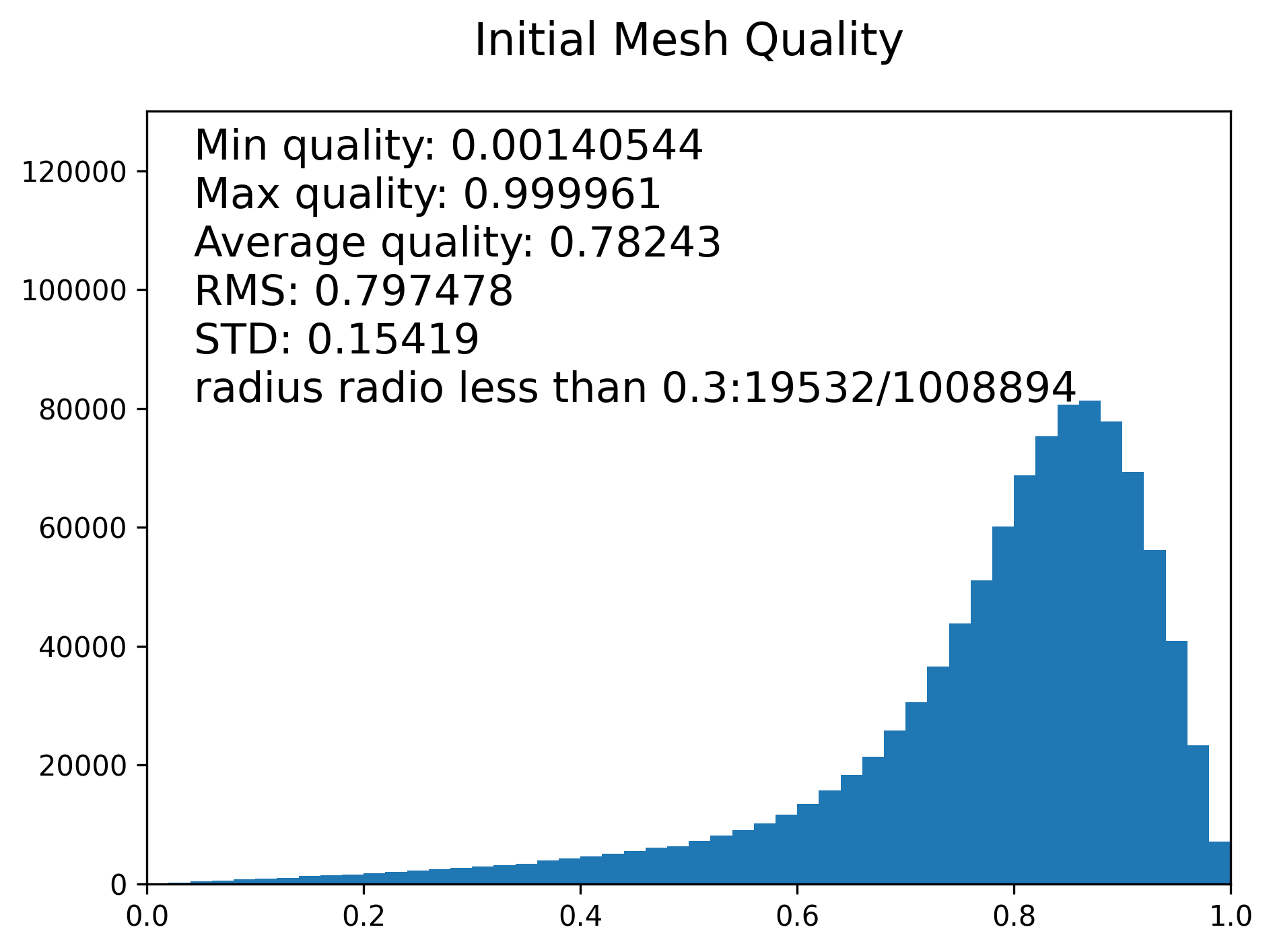}}
\subfloat[RRE]{
\includegraphics[width=0.3\linewidth]{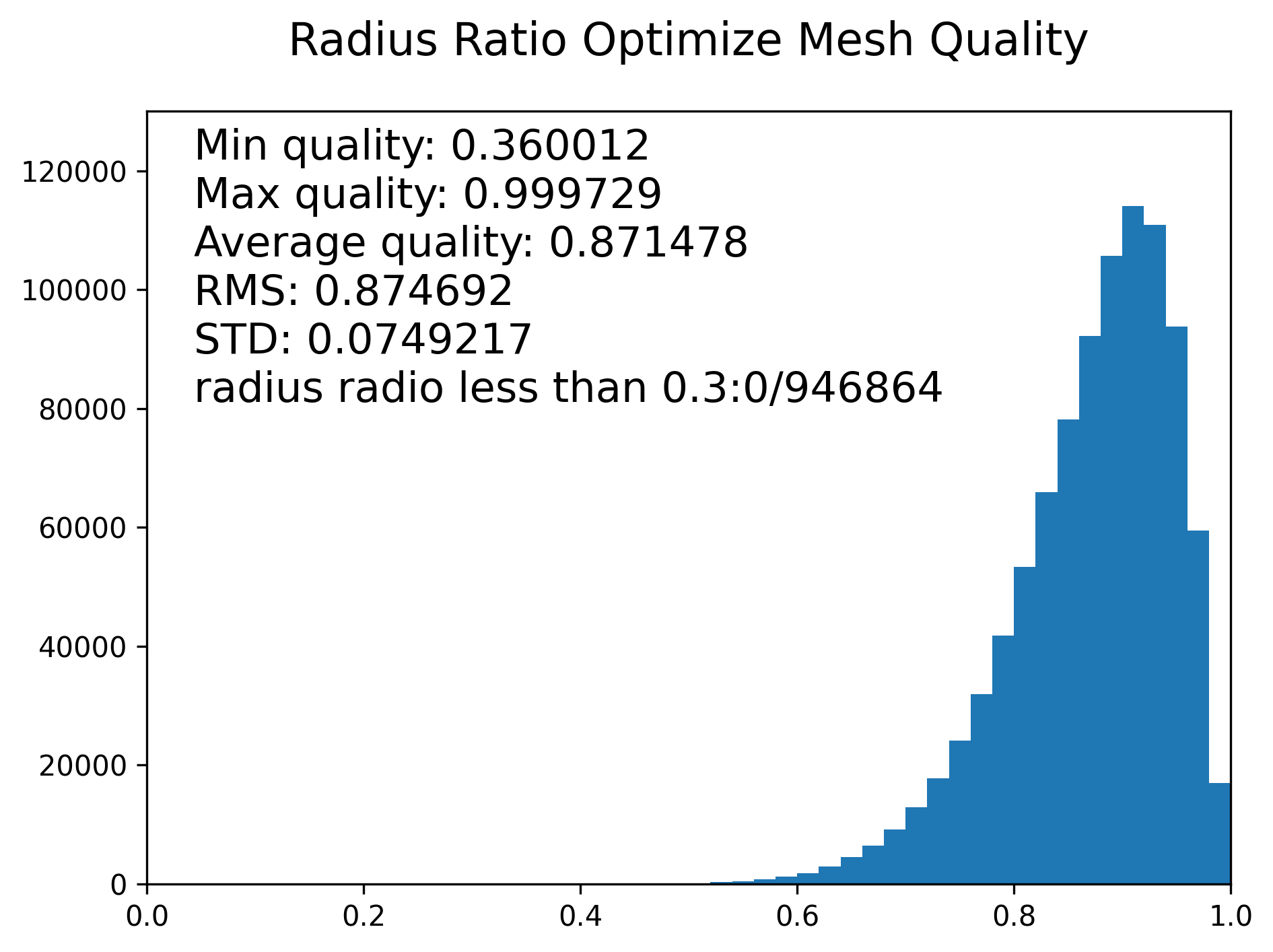}}
\hspace{0.01\linewidth}
\subfloat[Precondition RRE]{
\includegraphics[width=0.3\linewidth]{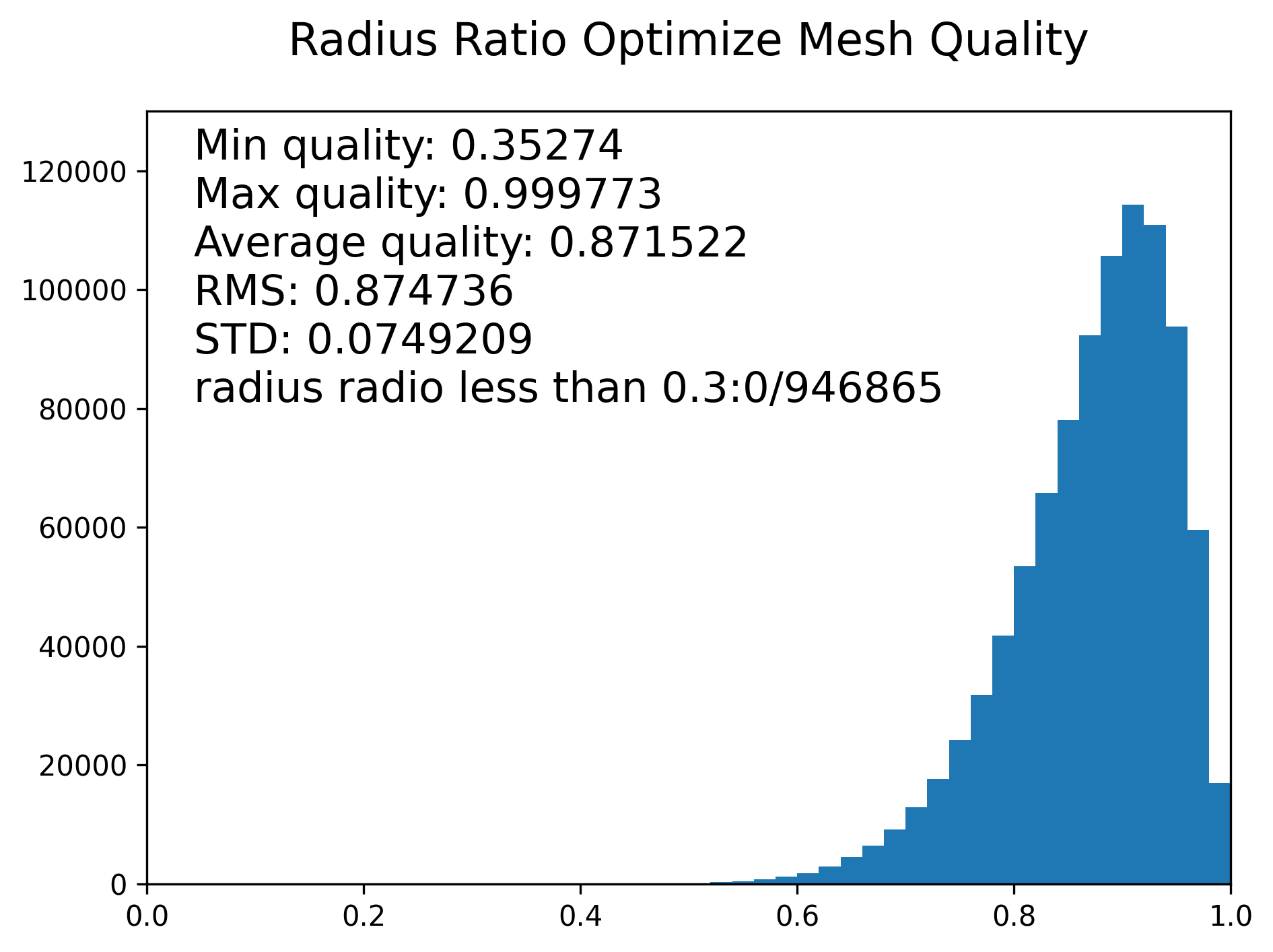}}
\caption{Intersection of 12 spheres}
\label{fig:12spherebig}
\end{figure}

\begin{table}[htbp]
    \caption{Comparison of RRE and precondition RRE}
    \label{table2}
    \centering
    %\begin{adjustwidth}{-1.5cm}{0cm}
    \begin{tabular}{|c|c|c|c|c|c|c|}
        \hline
        Model & Method & \makecell{Iteration count/\\Flip count}  &Num of cell &
        \makecell{min dihedral \\angle (deg.)} & \makecell{Min radius\\ ratio} & Time(sec.) \\
        \hline
        Sphere & Init & / & 677858 & 0.104 & 0.002 & /\\
        \hline
        & RRE &  385/3  & 635931 & 18.70 & 0.366 & 1006.18 \\
        \hline
        & Precondition RRE & 137/3 & 635932 & 18.27 & 0.367 & 453.80 \\
        \hline
        Lshape & Init & /  & 1210164 & 0.055 & 0.001 & / \\
        \hline
        & RRE & 407/4 & 1136025 & 16.59 & 0.339 & 1987.53 \\
        \hline
        & Precondition RRE& 204/4 & 1136027 & 15.36 & 0.318 & 1291.91 \\
        \hline
        \makecell{Intersection of\\ 12 spheres} &  Init & / & 1008894 & 0.804 & 0.001 & / \\
        \hline
        & RRE & 424/5 & 946864 & 19.70 & 0.360 & 1721.65 \\
        \hline
        & Precondition RRE & 206/5 & 946865 & 18.27 & 0.353 & 1165.99 \\
        \hline
\end{tabular}
%\end{adjustwidth}
\end{table}

Fig.~\ref{fig:spherebig},~\ref{fig:lshapebig},~\ref{fig:12spherebig}~show the 
results of RRE and preconditioned RRE on larger mesh, with up to millions of 
elements, for the three test cases. Table~\ref{table2} reports both the quality and the 
efficiency. From the figures, the two methods produce very similar distributions 
of the radius ratio and the dihedral angles. Although the minimum radius ratio 
and the minimum dihedral angle exhibit slight differences, the overall quality 
remains comparable. These discrepancies are mainly attributed to the different 
optimization trajectories induced by the preconditioner, rather than any change 
in the objective function. This indicates that the preconditioner does not 
change the final solution of the optimization problem, and the resulting mesh 
quality is essentially preserved.

From Table~\ref{table2}, the preconditioner clearly reduces the iteration count 
and the running time. Compared with small cases, the speedup is more obvious 
for large problems. This is because, for large mesh, the cost of energy and 
gradient evaluation dominates. Reducing the number of iterations then leads 
directly to a clear time saving. In summary, the preconditioner keeps the same 
optimization result while improving efficiency for large-scale problems.
\begin{figure}[htbp]
\centering
\subfloat[Model]{
\includegraphics[width=0.2\linewidth]{sphere.png}}
\subfloat[Init quality]{
\includegraphics[width=0.2\linewidth]{sphere_init.png}}
\hspace{0.01\linewidth}
\subfloat[RRE]{
\includegraphics[width=0.2\linewidth]{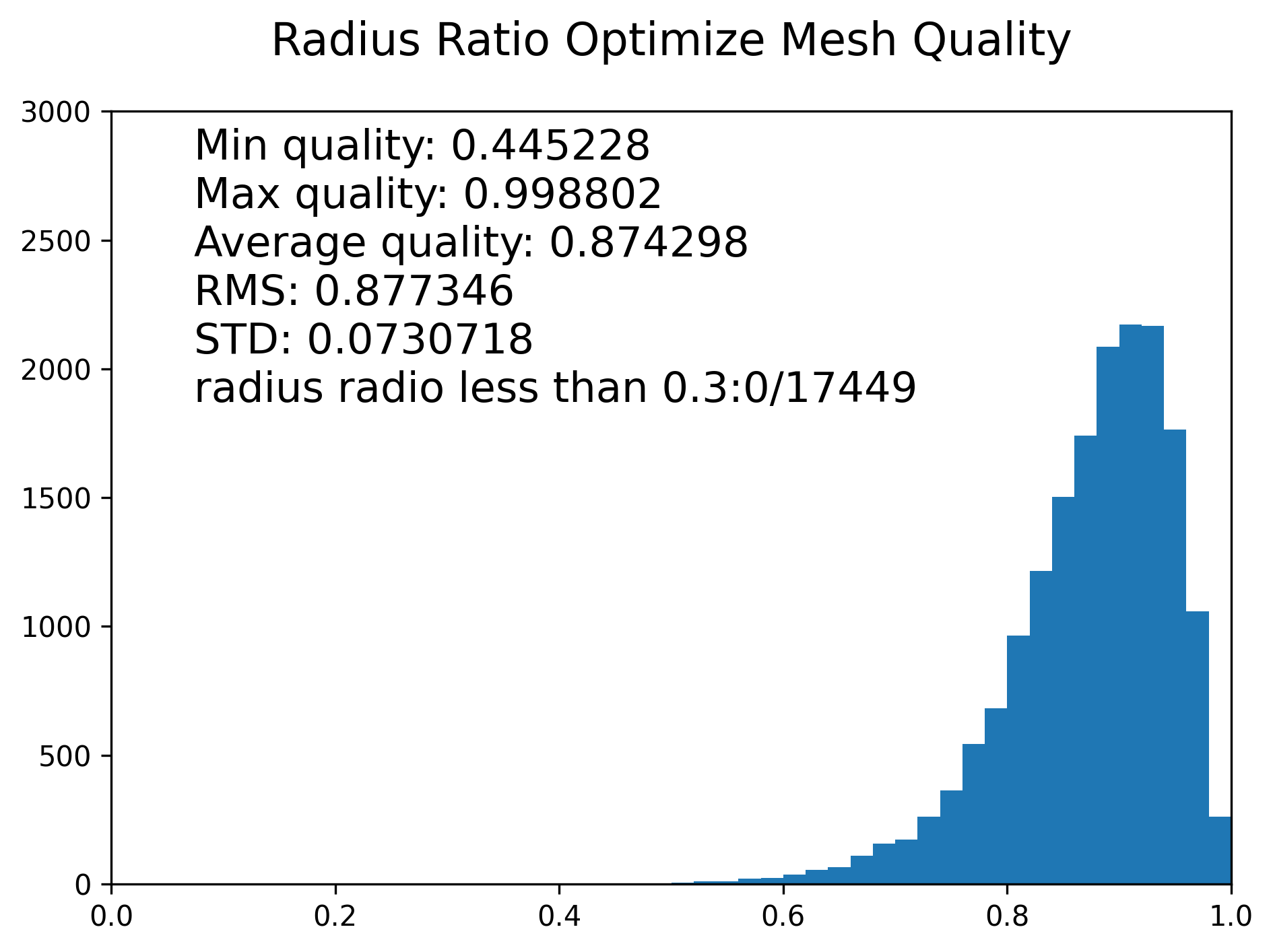}}
\hspace{0.01\linewidth}
\subfloat[Precondition RRE]{
\includegraphics[width=0.2\linewidth]{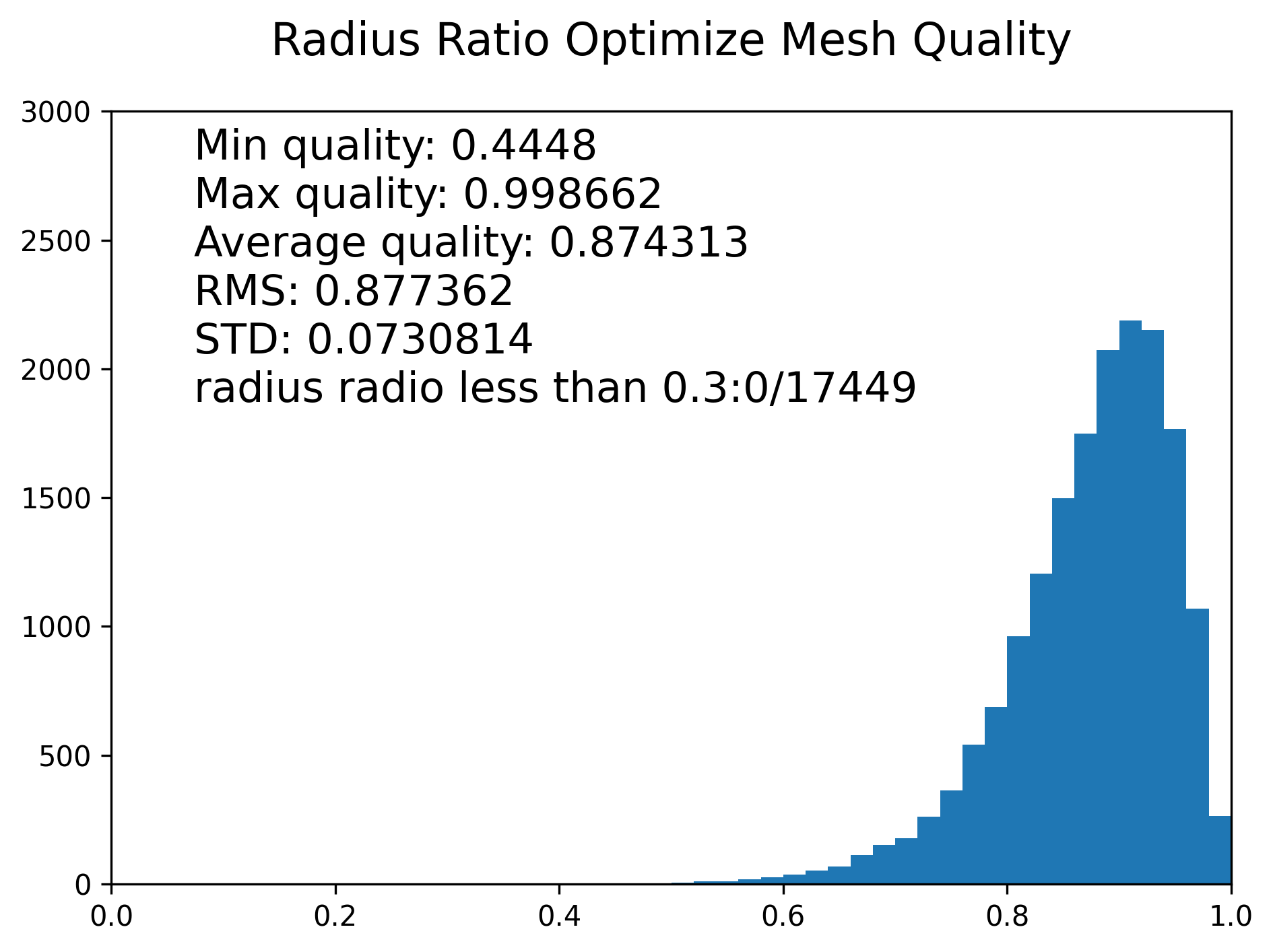}}
\vfill
\subfloat[Model]{
\includegraphics[width=0.2\linewidth]{lshape.png}}
\subfloat[Init quality]{
\includegraphics[width=0.2\linewidth]{lshape_init.png}}
\subfloat[RRE]{
\includegraphics[width=0.2\linewidth]{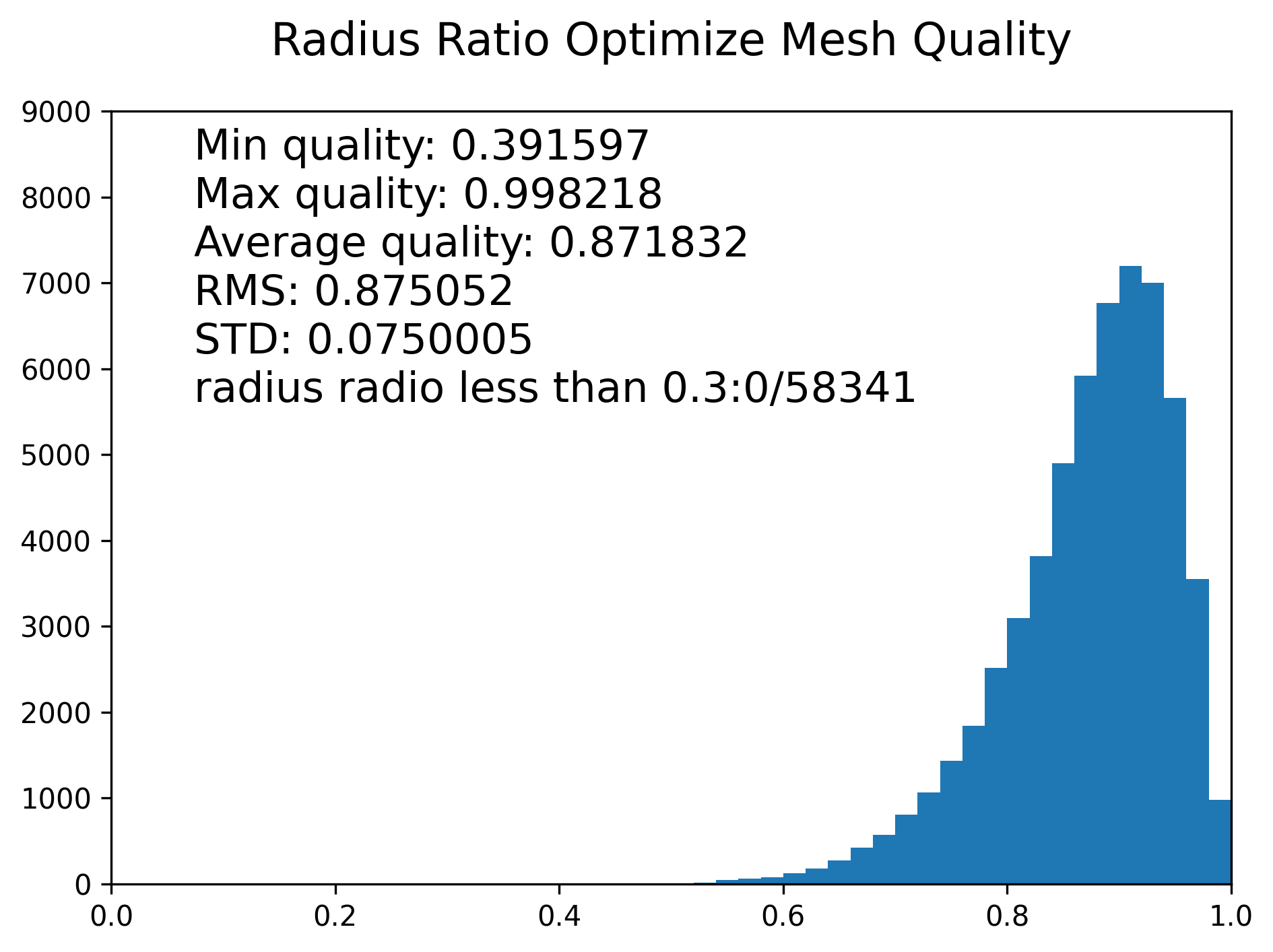}}
\subfloat[Precondition RRE]{
\includegraphics[width=0.2\linewidth]{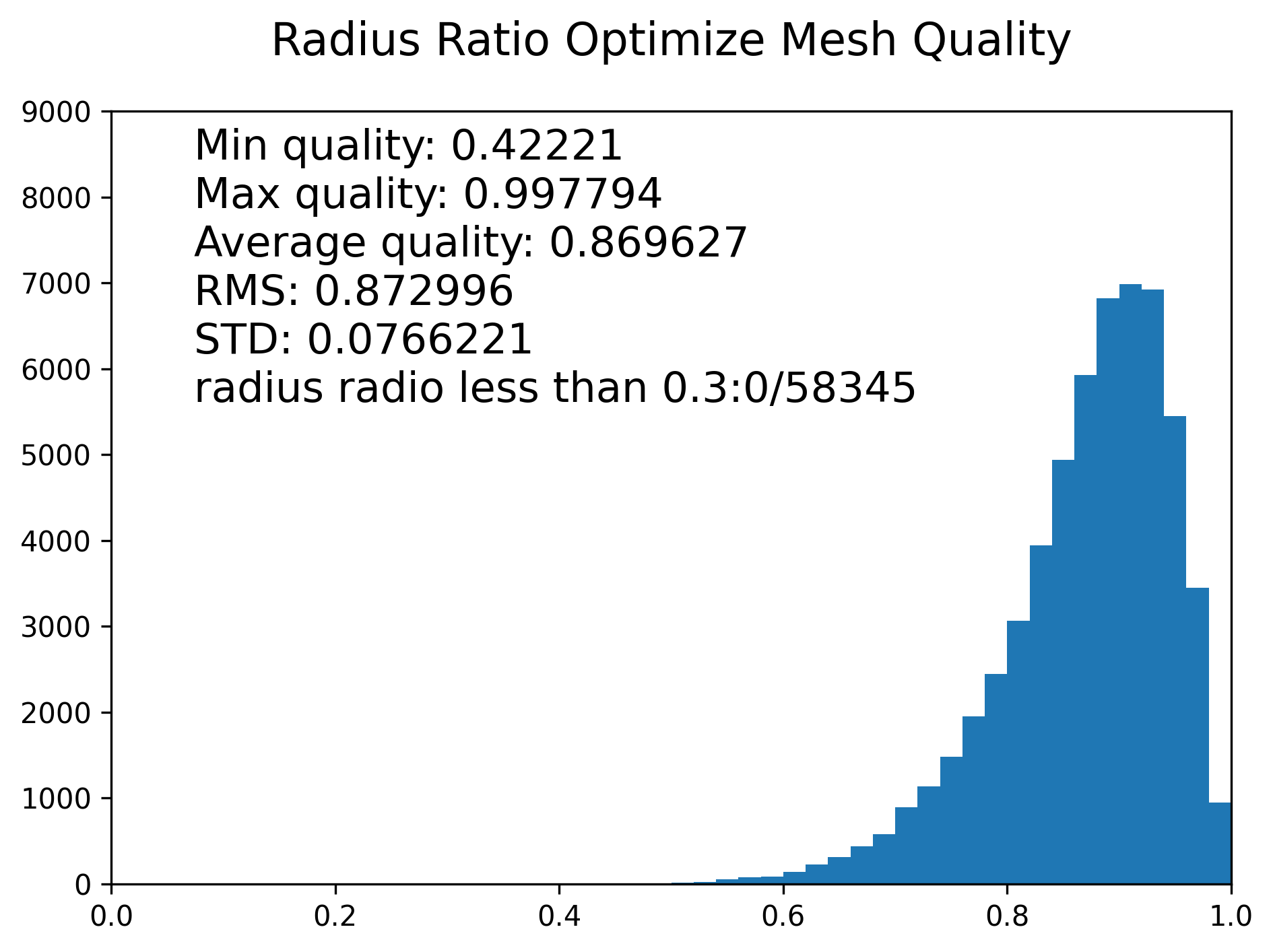}}
\vfill
\subfloat[Model]{
\includegraphics[width=0.2\linewidth]{12sphere.png}}
\subfloat[Init quality]{
\includegraphics[width=0.2\linewidth]{12sphere_init.png}}
\subfloat[RRE]{
\includegraphics[width=0.2\linewidth]{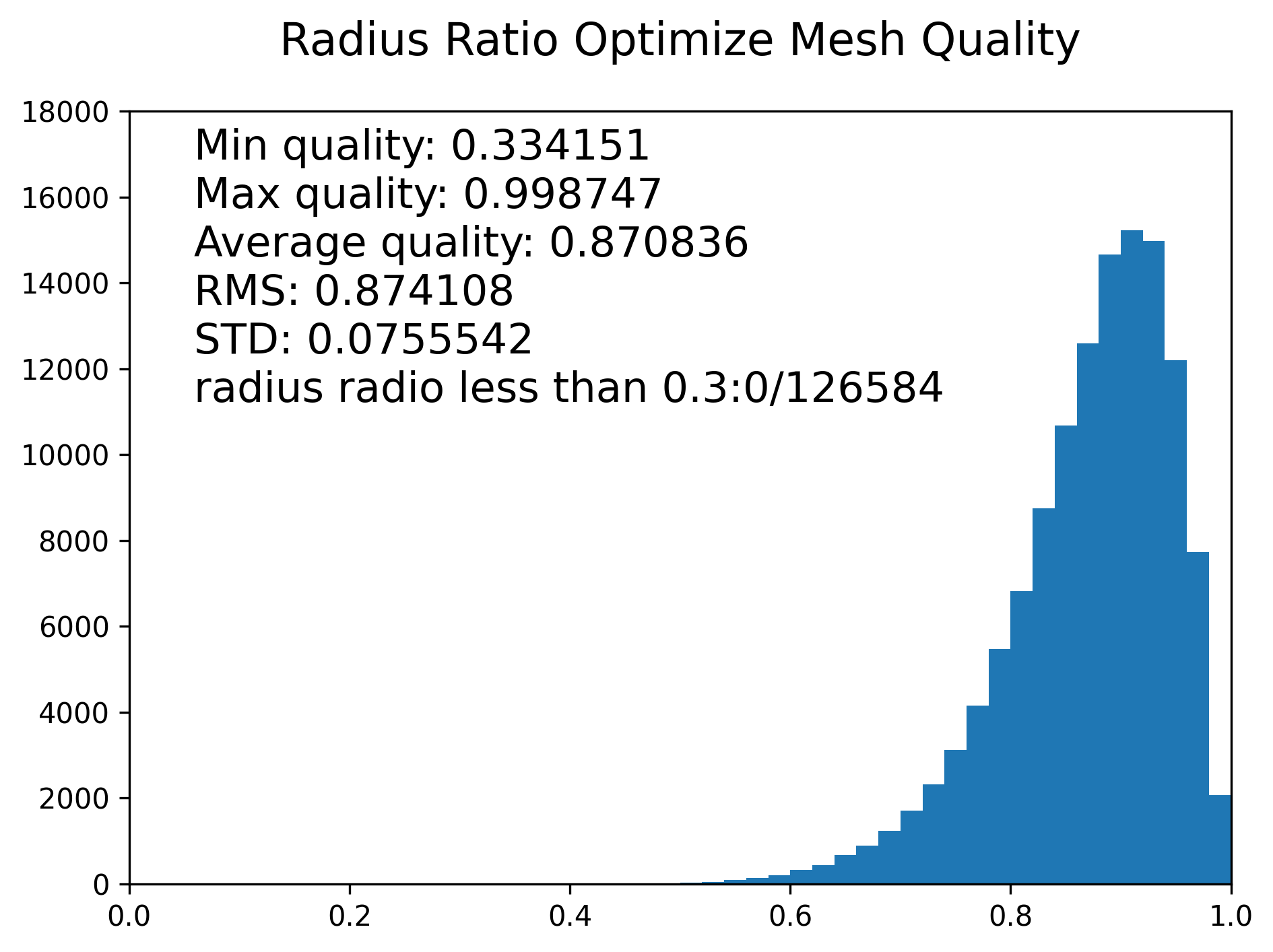}}
\subfloat[Precondition RRE]{
\includegraphics[width=0.2\linewidth]{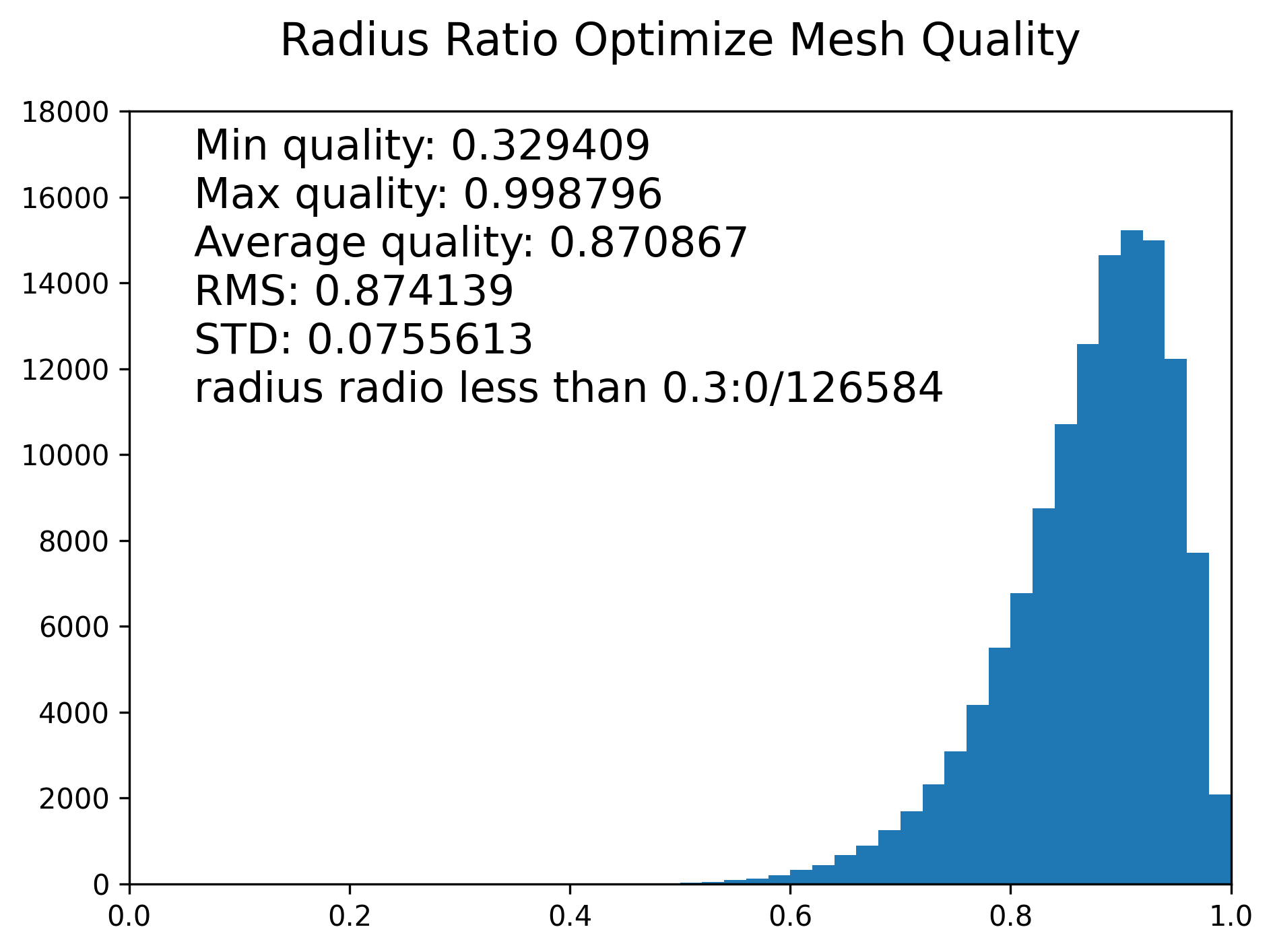}}
\caption{Use PNLCG algorithm instead of PLBFGS algorithm}
\label{fig:pnlcgexam}
\end{figure}

\begin{table}[htbp]
    \caption{Use PNLCG algorithm instead of PLBFGS algorithm}
    \label{table3}
    \centering
    %\begin{adjustwidth}{-1.5cm}{0cm}
    \begin{tabular}{|c|c|c|c|c|c|c|}
        \hline
        Model & Method & \makecell{Iteration count/\\Filp count} & Num of cell &
        \makecell{Min dihedral\\ angle(deg.)} & \makecell{Min radius \\ratio} & Time(sec.) \\
        \hline
        Sphere & Init & / & 18574 & 0.708 & 0.014 & /\\
        \hline
        & RRE & 113/2 & 17449 & 21.88 & 0.445 & 8.53 \\
        \hline
        & Precondition RRE & 59/2 & 17449 & 21.88 & 0.445 & 4.77 \\
        \hline
        Lshape & Init & / & 62002 & 0.516 & 0.011 & / \\
        \hline
        & RRE & 381/3 & 58341 & 18.67 & 0.392 & 83.09 \\
        \hline
        & Precondition RRE & 97/3 & 58345 & 19.24 & 0.422 & 25.29 \\
        \hline
        \makecell{Intersection of\\ 12 spheres} &  Init & / & 134710 & 0.230 & 9.63e-4 & /\\
        \hline
        & RRE & 363/3 & 126584 & 17.65 & 0.334 & 171.34 \\
        \hline
        & Precondition RRE & 124/3 & 126584 & 17.82 & 0.329 & 79.08 \\
        \hline
    \end{tabular}
%\end{adjustwidth}
\end{table}

The proposed preconditioner is designed to accelerate first-order optimization 
methods. In the previous examples, we used the PLBFGS method. However, the 
construction of the preconditioner depends on the structure of the radius ratio 
energy, not on a specific optimization method. In principle, it can be used 
with other methods that have a similar iteration structure. To verify this, 
Fig.~\ref{fig:pnlcgexam} and Table~\ref{table3} show results when the 
optimizer is replaced by the PNLCG method, for the sphere, L-shaped domain, 
and intersection of 12 spheres cases. The results show that RRE still improves 
mesh quality under PNLCG. With the preconditioner, both the iteration count and 
the running time are reduced, while the resulting mesh quality remains comparable 
to that obtained without preconditioning. This shows that the preconditioner also works 
well in the PNLCG framework. It improves the search direction and speeds up 
convergence. Therefore, the proposed preconditioner is general and can be used 
with different first-order methods.

\begin{figure}[htbp]
\centering
\subfloat[Model]{
\includegraphics[width=0.2\linewidth]{sphere.png}}
\subfloat[Init quality]{
\includegraphics[width=0.2\linewidth]{sphere_init.png}}
\hspace{0.01\linewidth}
\subfloat[RRE]{
\includegraphics[width=0.2\linewidth]{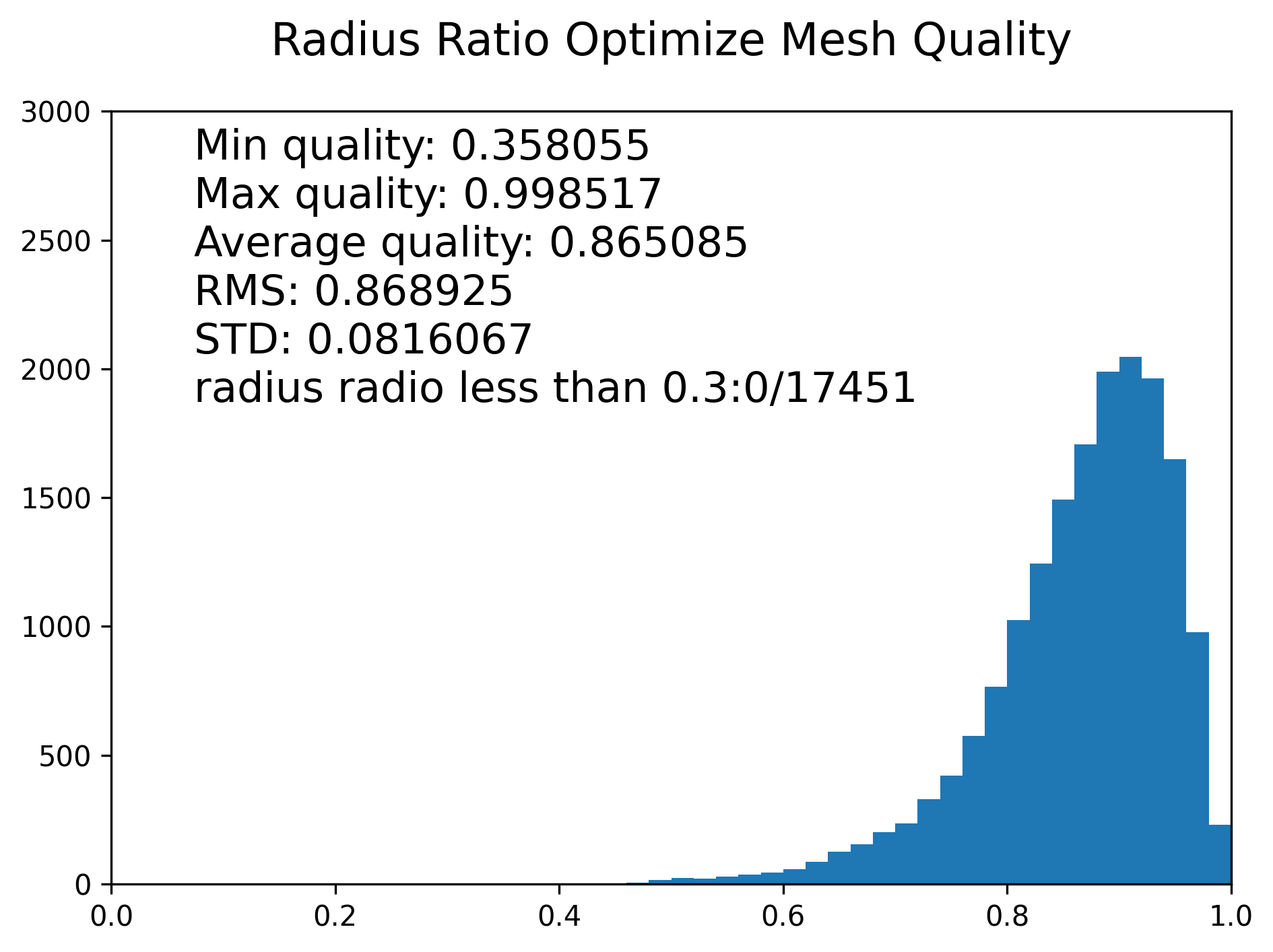}}
\hspace{0.01\linewidth}
\subfloat[Precondition RRE]{
\includegraphics[width=0.2\linewidth]{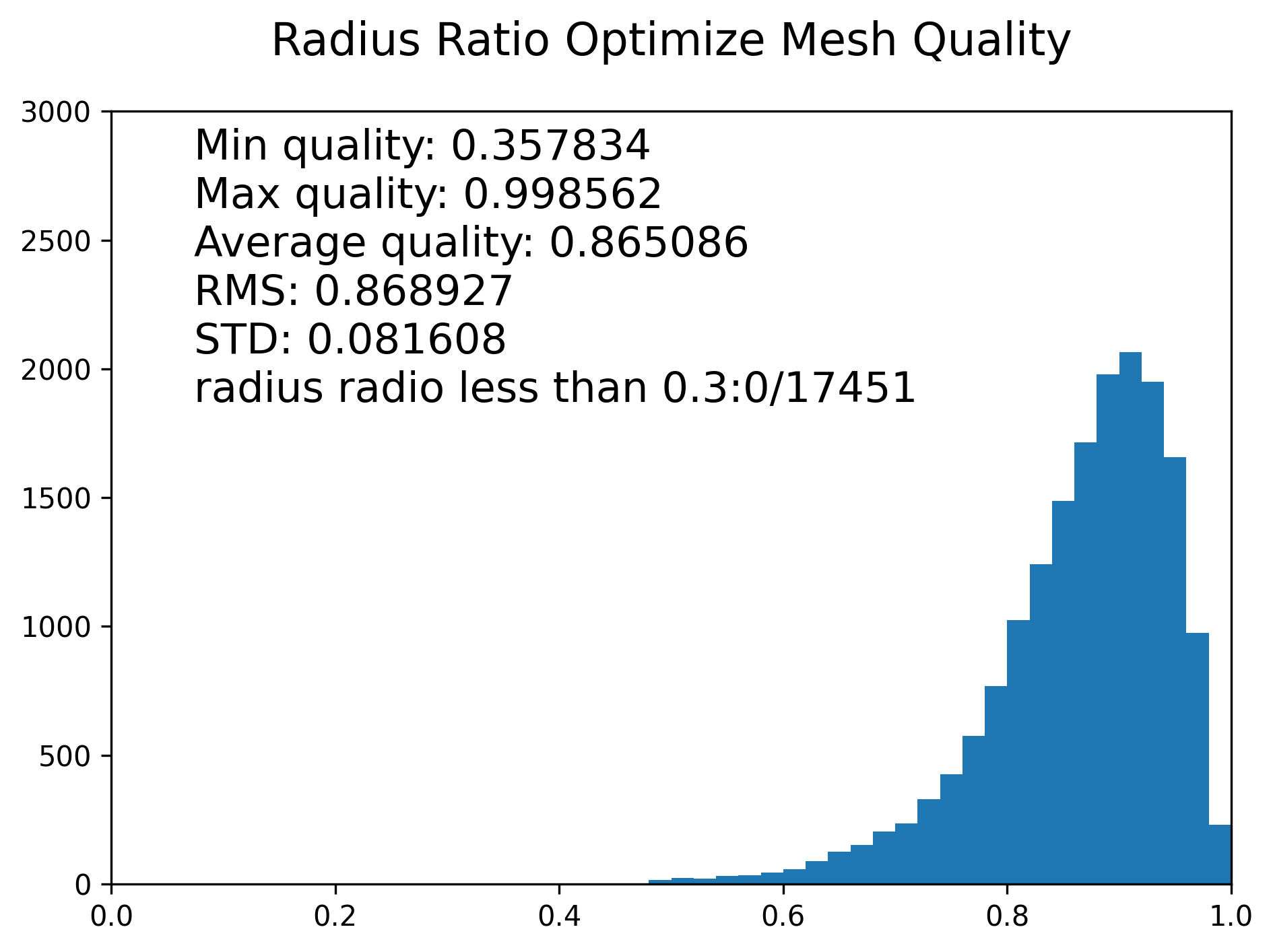}}
\vfill
\subfloat[Model]{
\includegraphics[width=0.2\linewidth]{lshape.png}}
\subfloat[Init quality]{
\includegraphics[width=0.2\linewidth]{lshape_init.png}}
\subfloat[RRE]{
\includegraphics[width=0.2\linewidth]{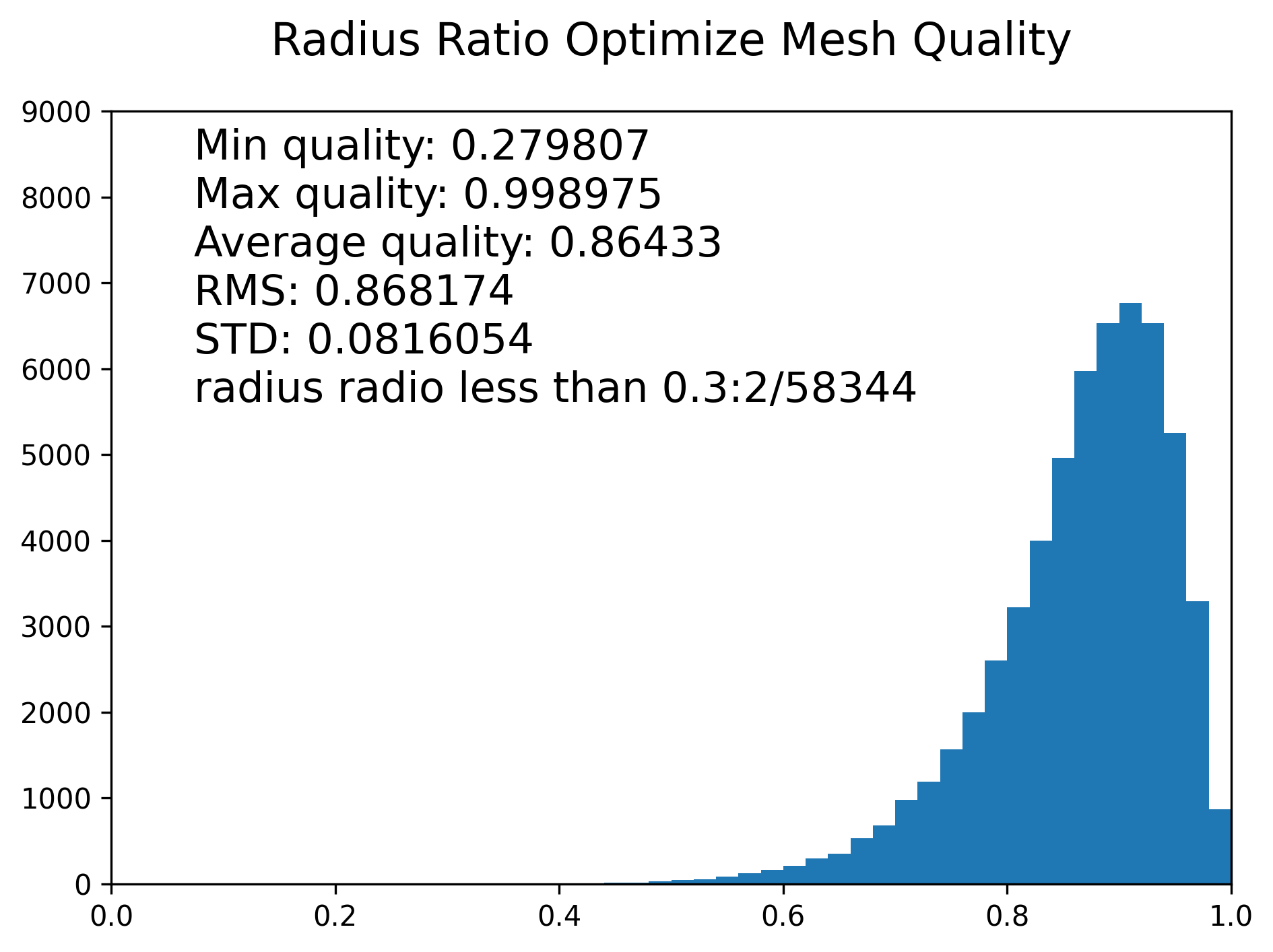}}
\subfloat[Precondition RRE]{
\includegraphics[width=0.2\linewidth]{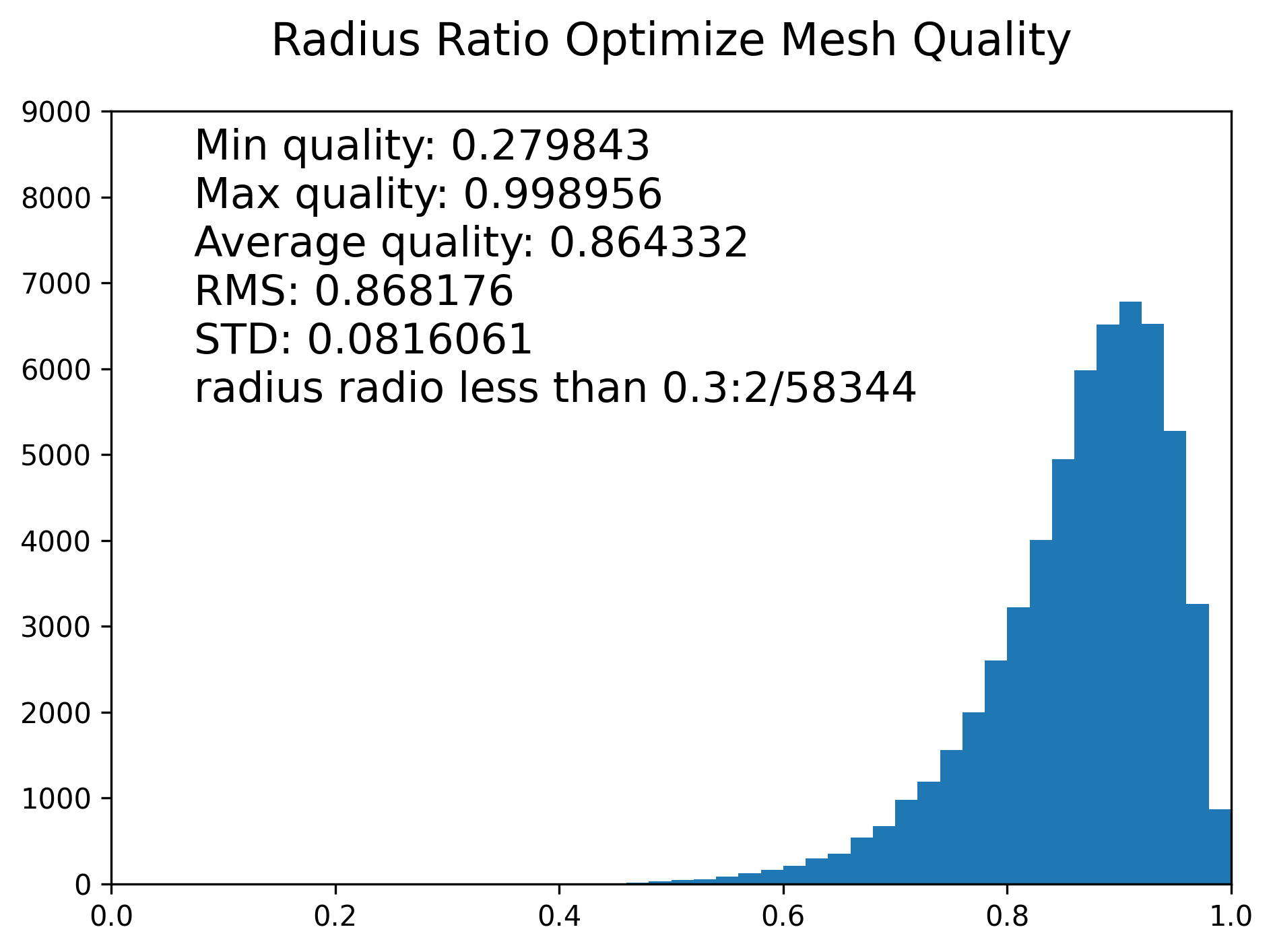}}
\vfill
\subfloat[Model]{
\includegraphics[width=0.2\linewidth]{12sphere.png}}
\subfloat[Init quality]{
\includegraphics[width=0.2\linewidth]{12sphere_init.png}}
\subfloat[RRE]{
\includegraphics[width=0.2\linewidth]{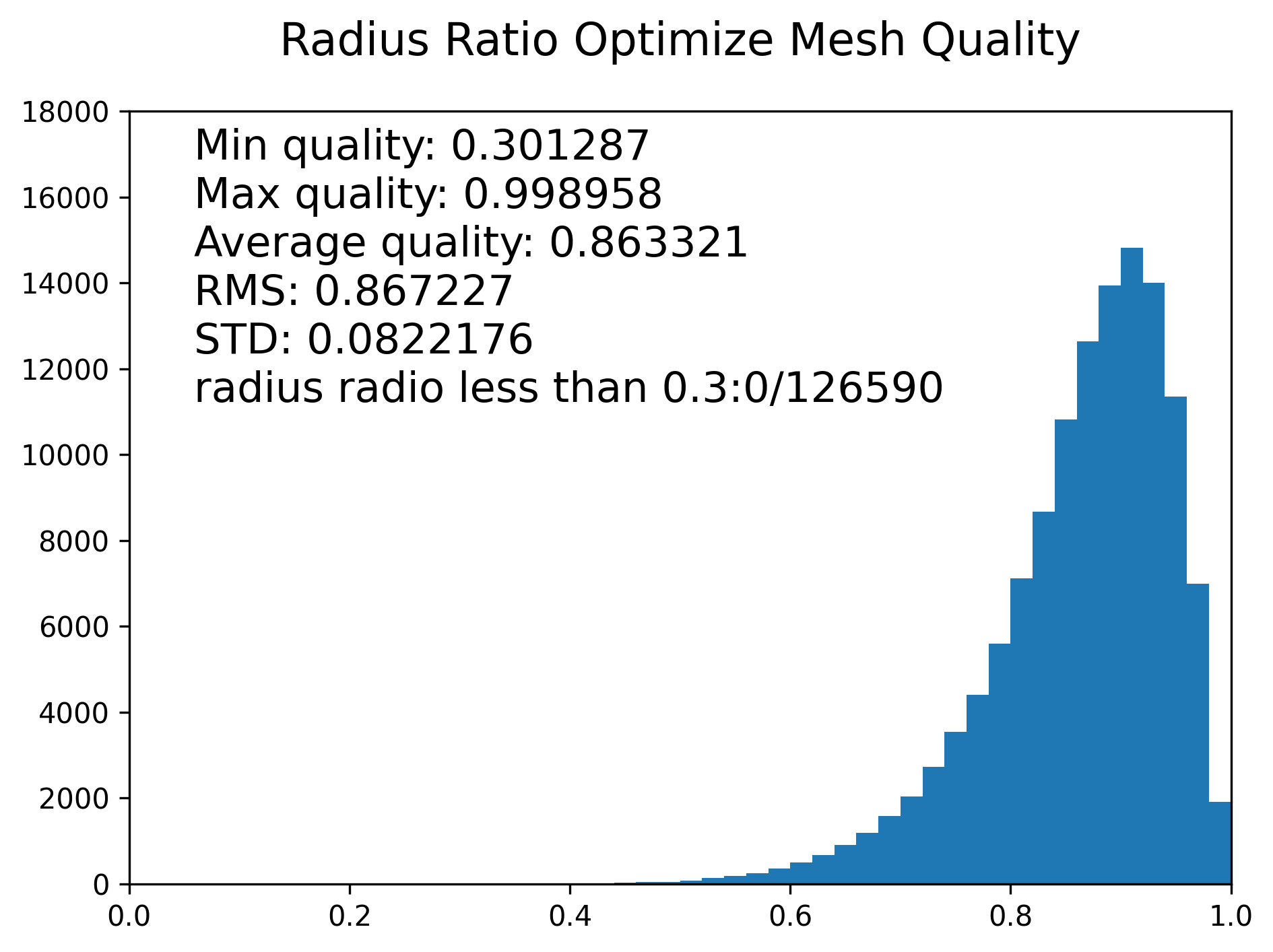}}
\subfloat[Precondition RRE]{
\includegraphics[width=0.2\linewidth]{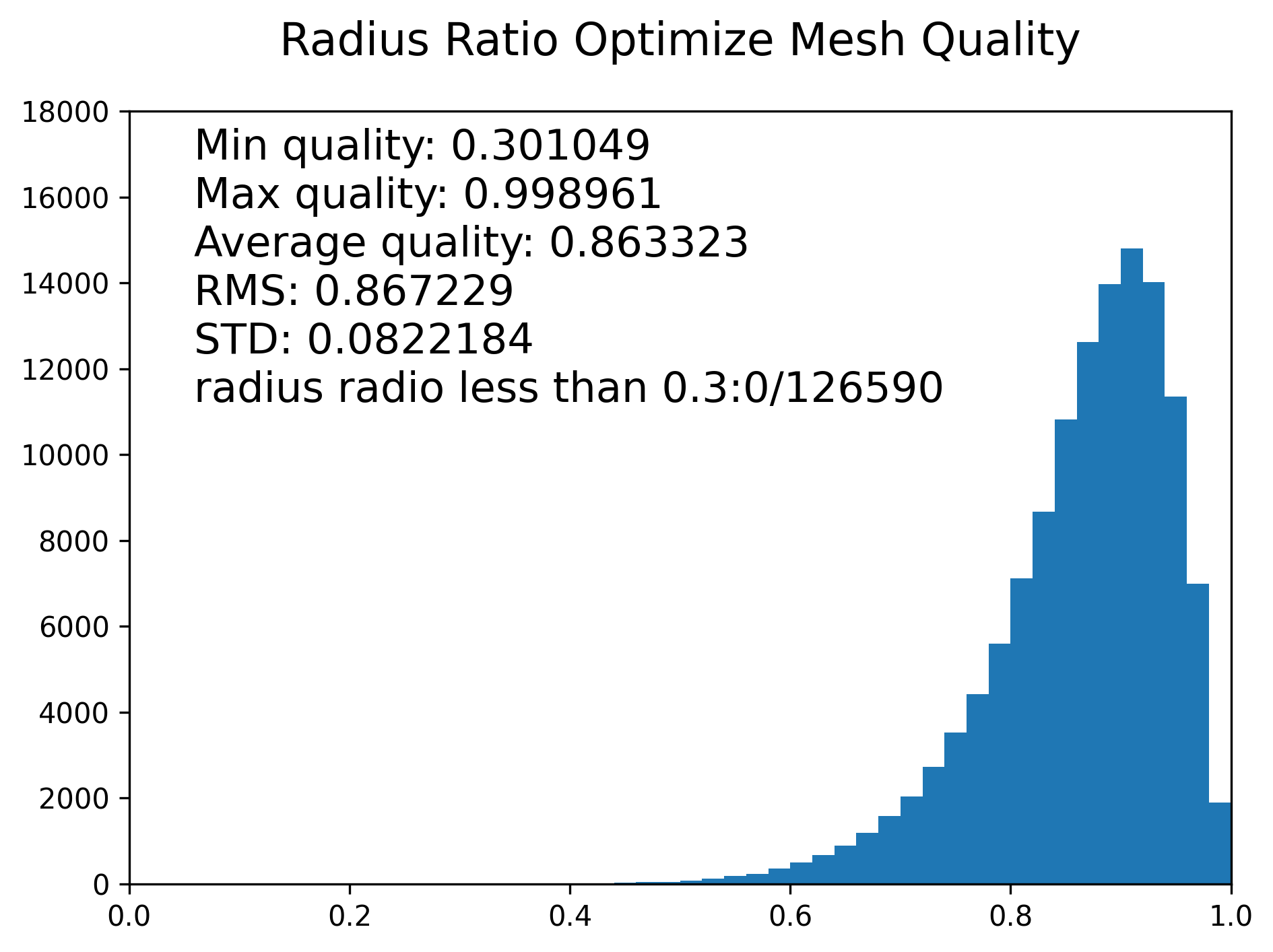}}
\caption{Fixed boundary point RRE mesh optimization}
\label{fig:bdfixexam}
  \end{figure}

\begin{table}[htbp]
    \caption{Fixed boundary point RRE mesh optimization}
    \label{table5}
    \centering
    %\begin{adjustwidth}{-1.5cm}{0cm}
    \begin{tabular}{|c|c|c|c|c|c|c|}
        \hline
        Model & Method & \makecell{Iteration count/\\Filp count} & Num of cell &
        \makecell{Min dihedral\\ angle(deg.)} & \makecell{Min radius\\ ratio} & Time(sec.) \\
        \hline
        Sphere & Init & / & 18574 & 0.708 & 0.014 & /\\
        \hline
        & RRE & 95/2 & 17451 & 17.18 & 0.358 & 5.60 \\
        \hline
        & Precondition RRE & 46/2 & 17451 & 17.18 & 0.358 & 3.34 \\
        \hline
        Lshape & Init & / & 62002 & 0.516 & 0.011 & / \\
        \hline
        & RRE & 154/3 & 58344 & 16.68 & 0.280 & 34.38 \\
        \hline
        & Precondition RRE & 73/3 & 58344 & 16.68 & 0.280 & 20.51 \\
        \hline
        \makecell{Intersection of\\ 12 spheres} & Init & / & 134710 & 0.230 & 9.63e-4 & /\\
        \hline
        & RRE & 152/5 & 126590 & 14.84 & 0.301 & 77.95 \\
        \hline
        & Precondiiton RRE & 103/5 & 126590 & 14.83 & 0.301 & 65.59 \\
        \hline
    \end{tabular}
%\end{adjustwidth}
\end{table}

Since this work focuses on mesh optimization rather than mesh generation, it is
reasonable to assume that the boundary node distribution is well-shaped and the
boundary mesh quality is satisfactory. Under this assumption, we fix the
boundary nodes and perform optimization only on the interior nodes. Fig.
~\ref{fig:bdfixexam} and Table \ref{table5} present the corresponding
numerical results. The results show that, even with fixed boundary nodes, RRE 
can still significantly improve the overall mesh quality and effectively 
enhance low-quality elements.  Compared with the case where boundary nodes can 
move, the improvement is slightly smaller, but still stable and effective. 
Moreover, since fixing boundary nodes eliminates the need for boundary movement 
and projection operations, and reduces the number of degrees of freedom involved 
in the optimization, the convergence speed is improved and the computational 
efficiency is higher. These results indicate that, in practical applications 
where the boundary mesh quality is already satisfactory, adopting a 
fixed-boundary optimization strategy can further improve computational 
efficiency while maintaining good optimization performance.
\subsection{Anisotropic mesh}
For anisotropic mesh, we first generate an initial isotropic mesh and then 
optimize it according to a given anisotropic metric. Below we present several 
two-dimensional and three-dimensional examples.

\begin{figure}[htbp]
\centering
\subfloat[Init mesh]{
\includegraphics[width=0.2\linewidth]{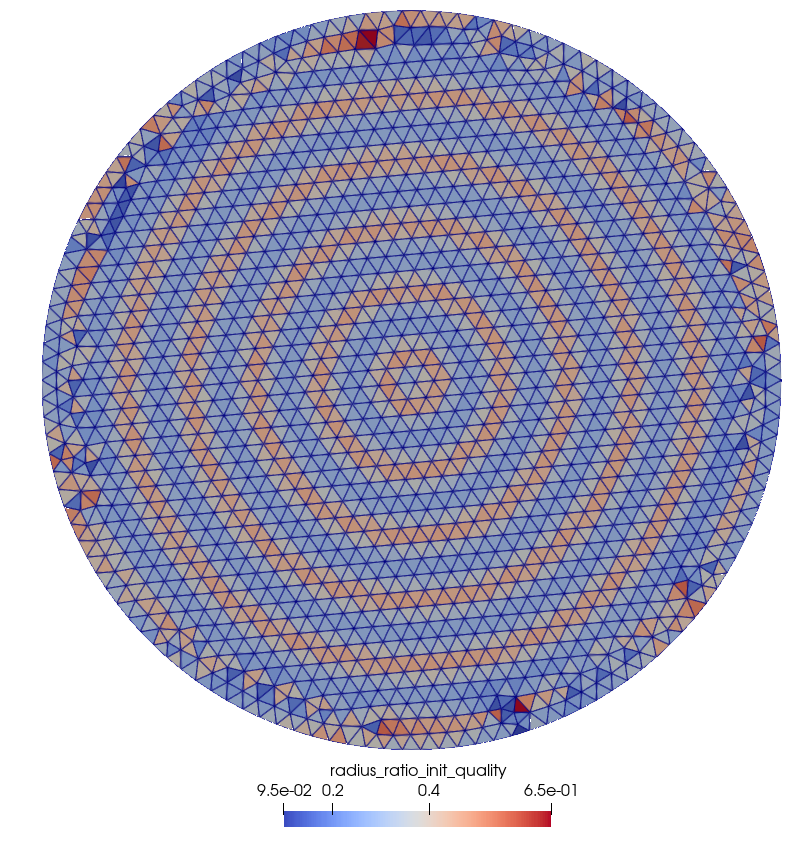}}
\hspace{0.01\linewidth}
\subfloat[Optimized mesh]{
\includegraphics[width=0.2\linewidth]{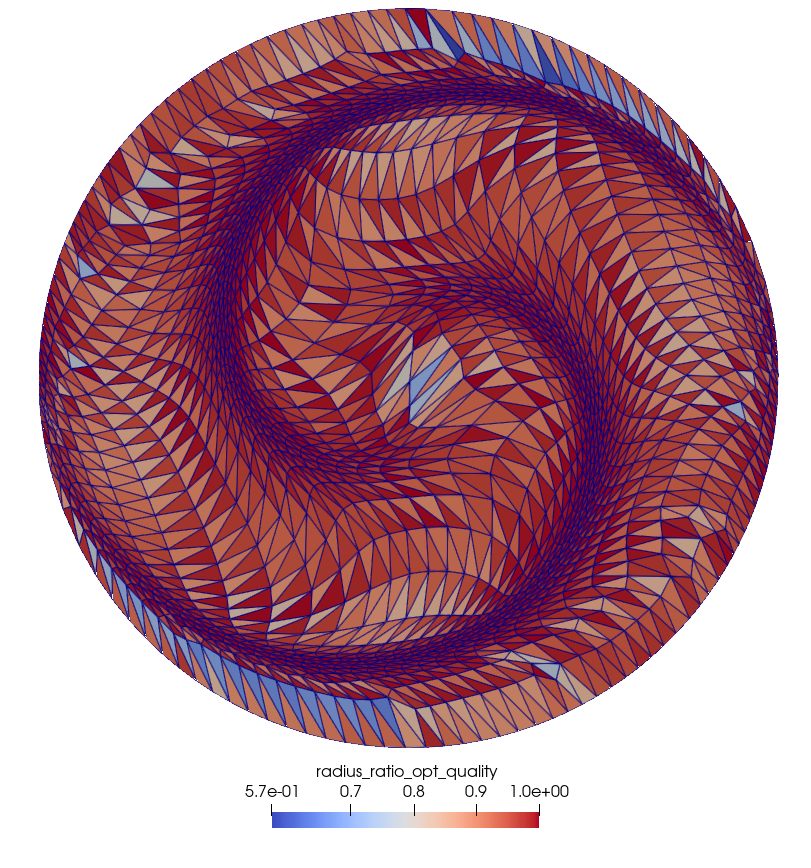}}
\vspace{0.01\linewidth}
\subfloat[Init quality]{
\includegraphics[width=0.2\linewidth]{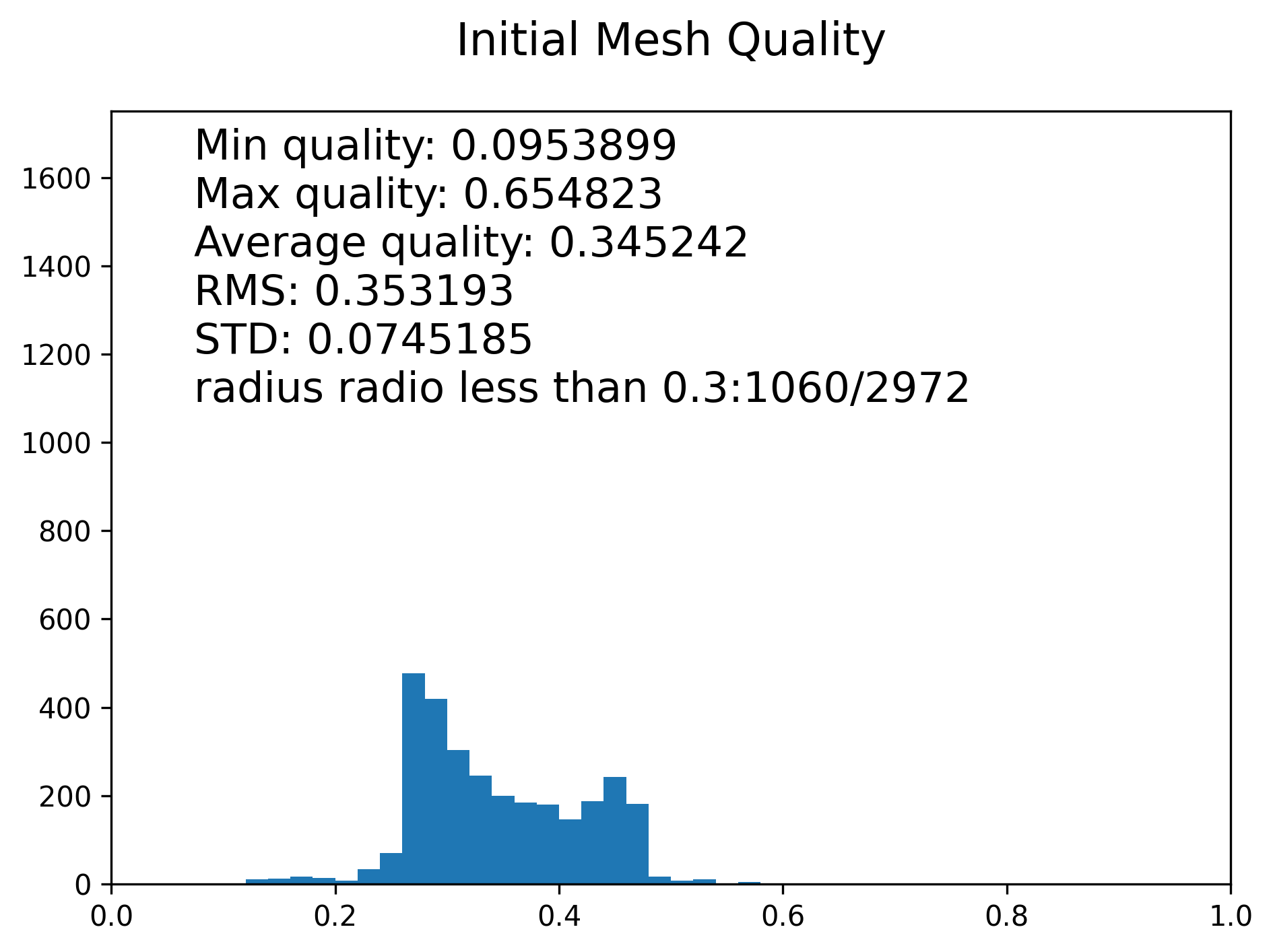}}
\hspace{0.01\linewidth}
\subfloat[RRE quality]{
\includegraphics[width=0.2\linewidth]{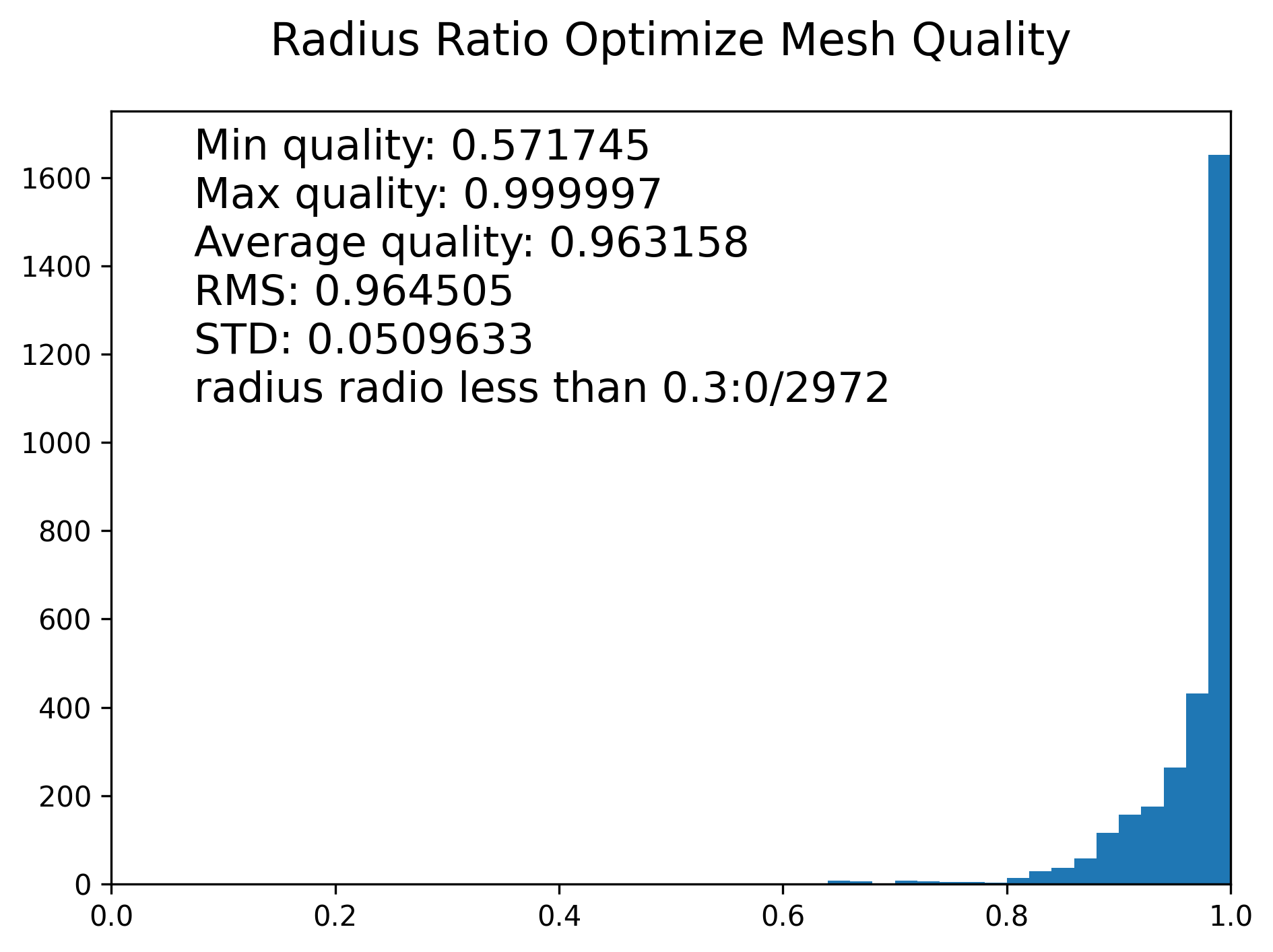}}
\caption{Spiral anisotropic mesh}
\label{fig:anicir1}
\end{figure}

Fig.~\ref{fig:anicir1} shows the unit disk with the anisotropic metric  
$\bm{M}(\bm{x})=\bm{Q}(\bm{x})^T\Lambda\bm{Q}(\bm{x})$ where
$$
\Lambda=\text{diag}(25,1),
\bm{Q}(\bm{x})=
\begin{bmatrix}
\cos(\phi(\bm{x})) & -\sin(\phi(\bm{x})) \\
\sin(\phi(\bm{x}) & \cos(\phi(\bm{x})))
\end{bmatrix}
$$
Let $r=\sqrt{x^2+y^2},\phi(\bm{x})=kr,k=6$. Under this metric, the mesh forms a 
spiral anisotropic pattern.

\begin{figure}[htbp]
\centering
\subfloat[Init mesh]{
\includegraphics[width=0.2\linewidth]{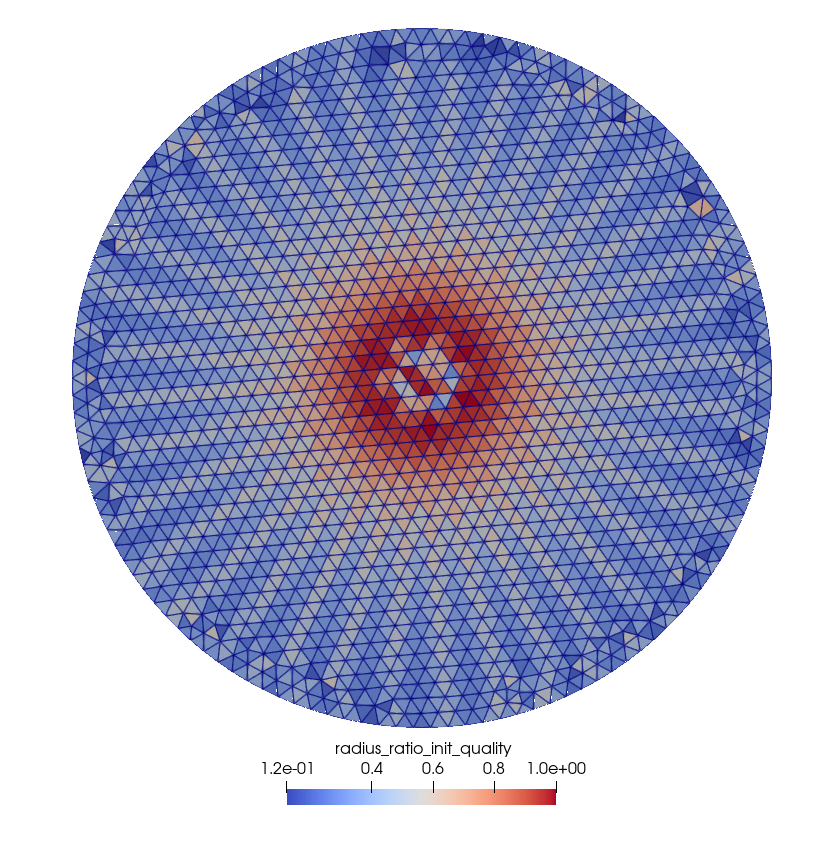}}
\hspace{0.01\linewidth}
\subfloat[Optimized mesh]{
\includegraphics[width=0.2\linewidth]{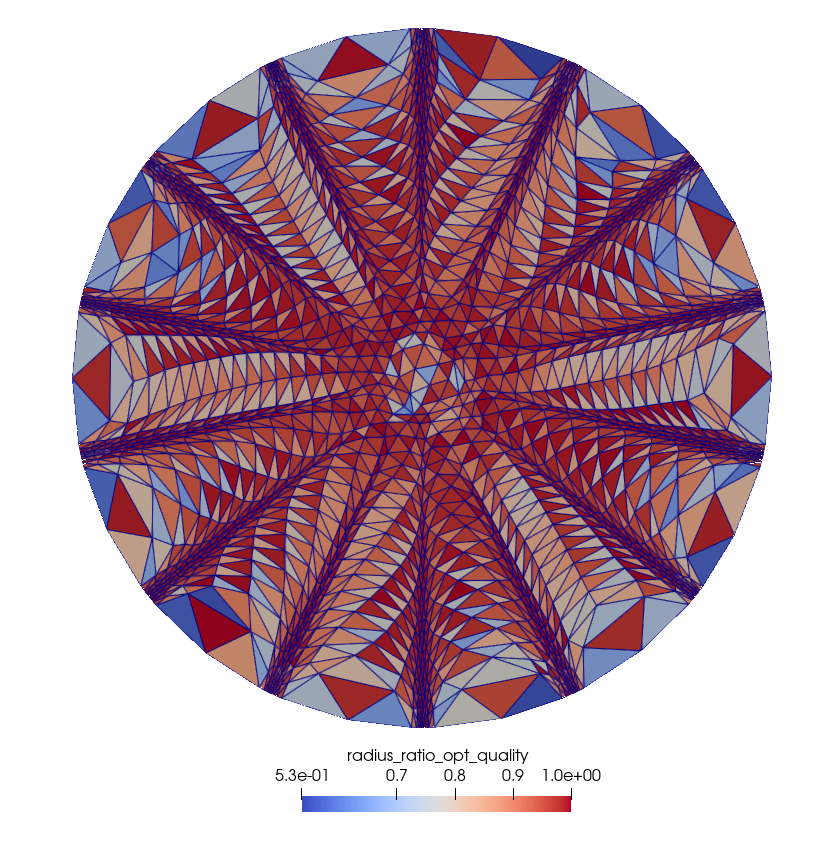}}
\vspace{0.01\linewidth}
\subfloat[Init quality]{
\includegraphics[width=0.2\linewidth]{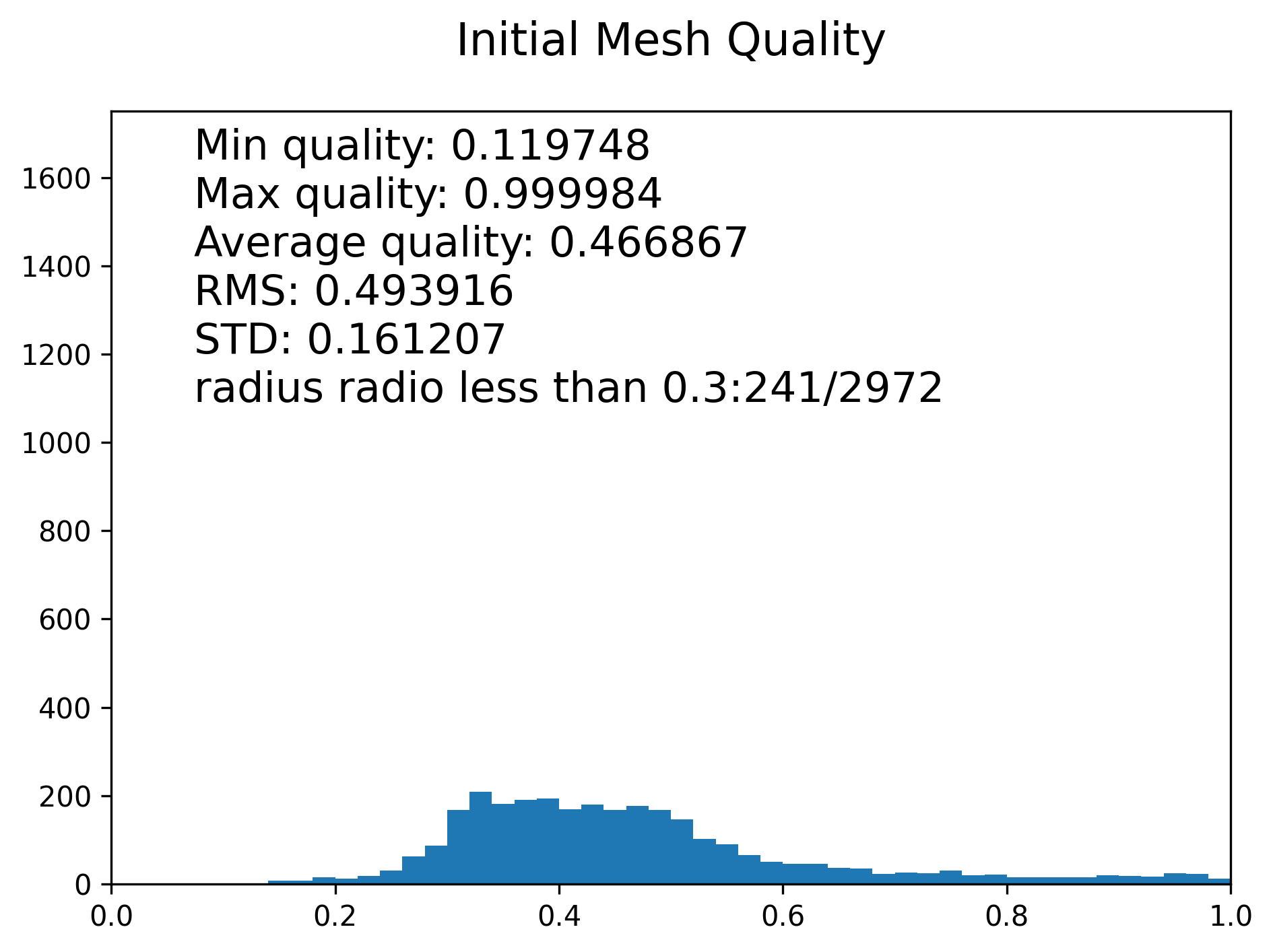}}
\hspace{0.01\linewidth}
\subfloat[RRE quality]{
\includegraphics[width=0.2\linewidth]{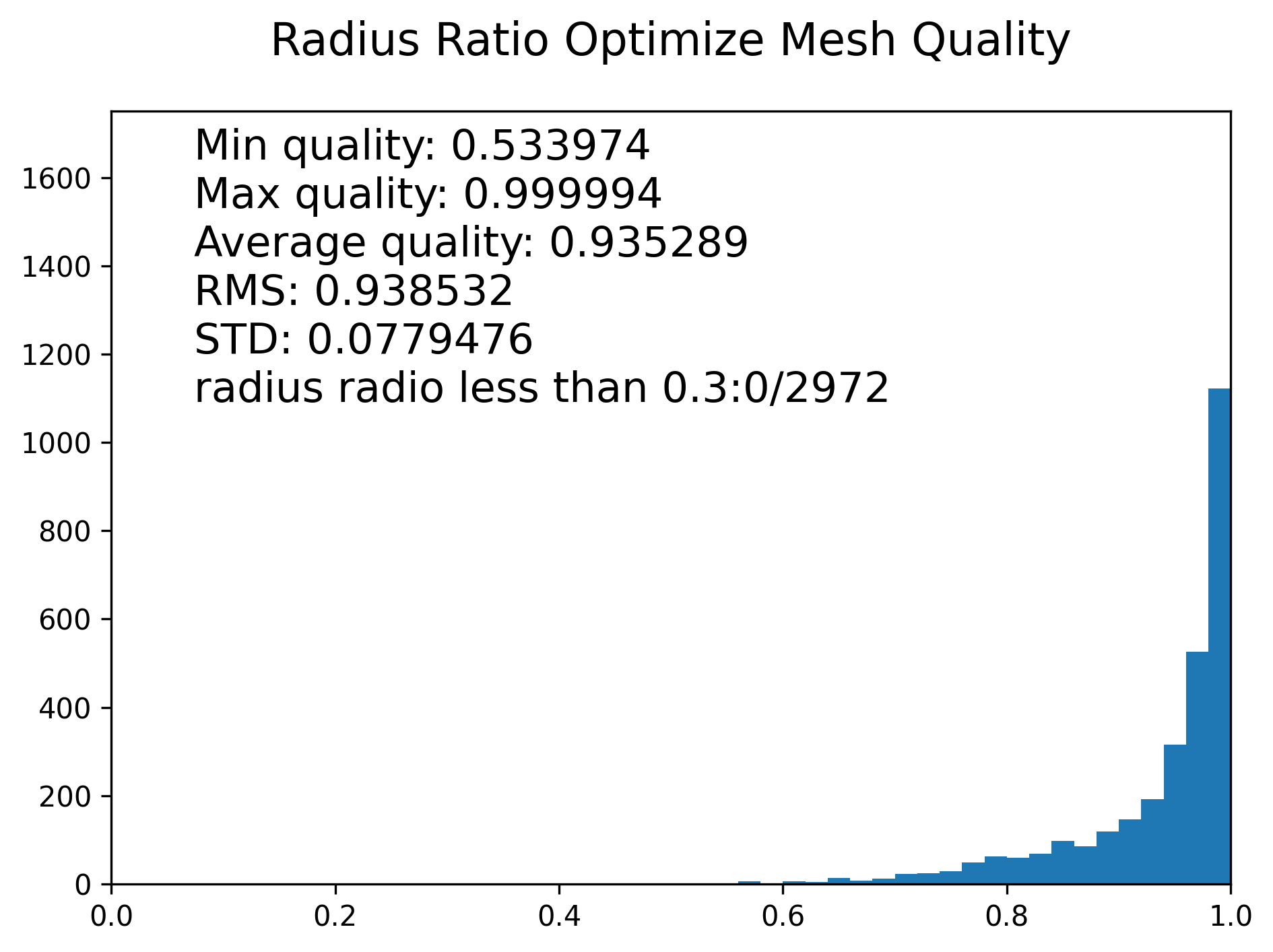}}
\caption{Radial anisotropic mesh}
\label{fig:anicir2}
\end{figure}

If we modify the metric by setting  $\theta=\arctan2(y,x),\phi(\bm{x})=k\theta,
k=6$, the mesh forms a radial anisotropic pattern from the center, as shown in 
Fig.~\ref{fig:anicir2}.

\begin{figure}[htbp]
\centering
\subfloat[Init mesh]{
\includegraphics[width=0.2\linewidth]{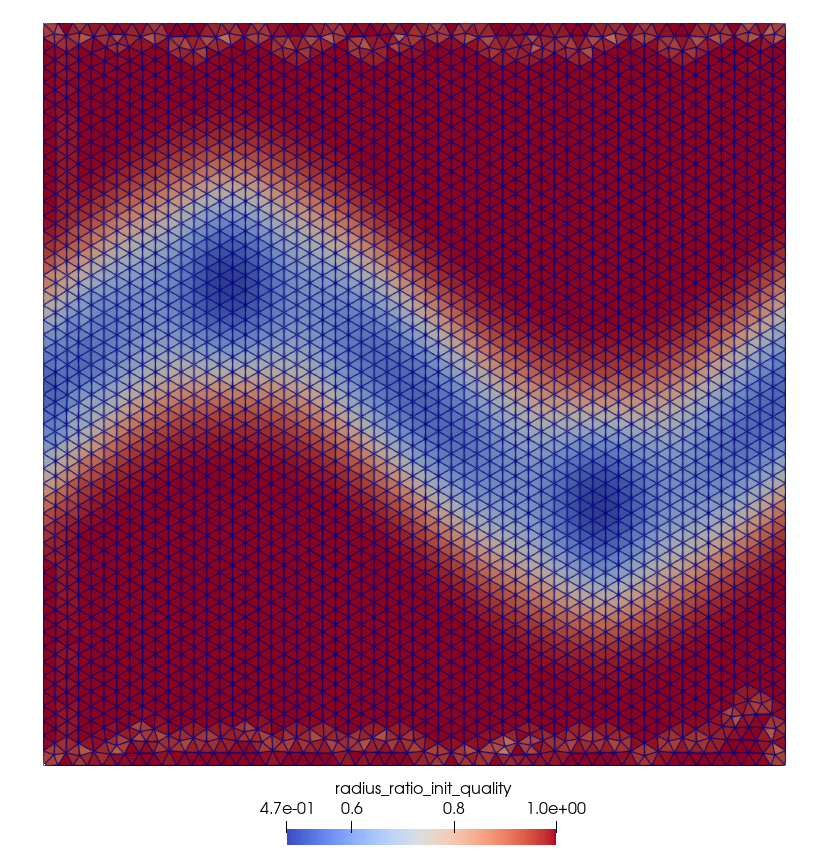}}
\hspace{0.01\linewidth}
\subfloat[Optimized mesh]{
\includegraphics[width=0.2\linewidth]{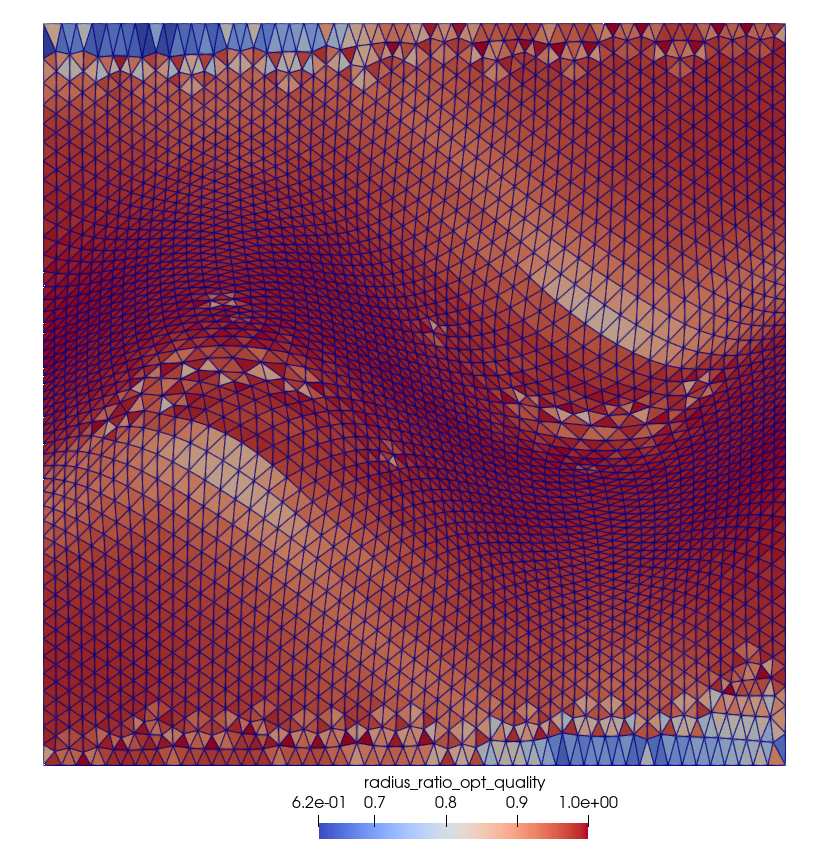}}
\vspace{0.01\linewidth}
\subfloat[Init quality]{
\includegraphics[width=0.2\linewidth]{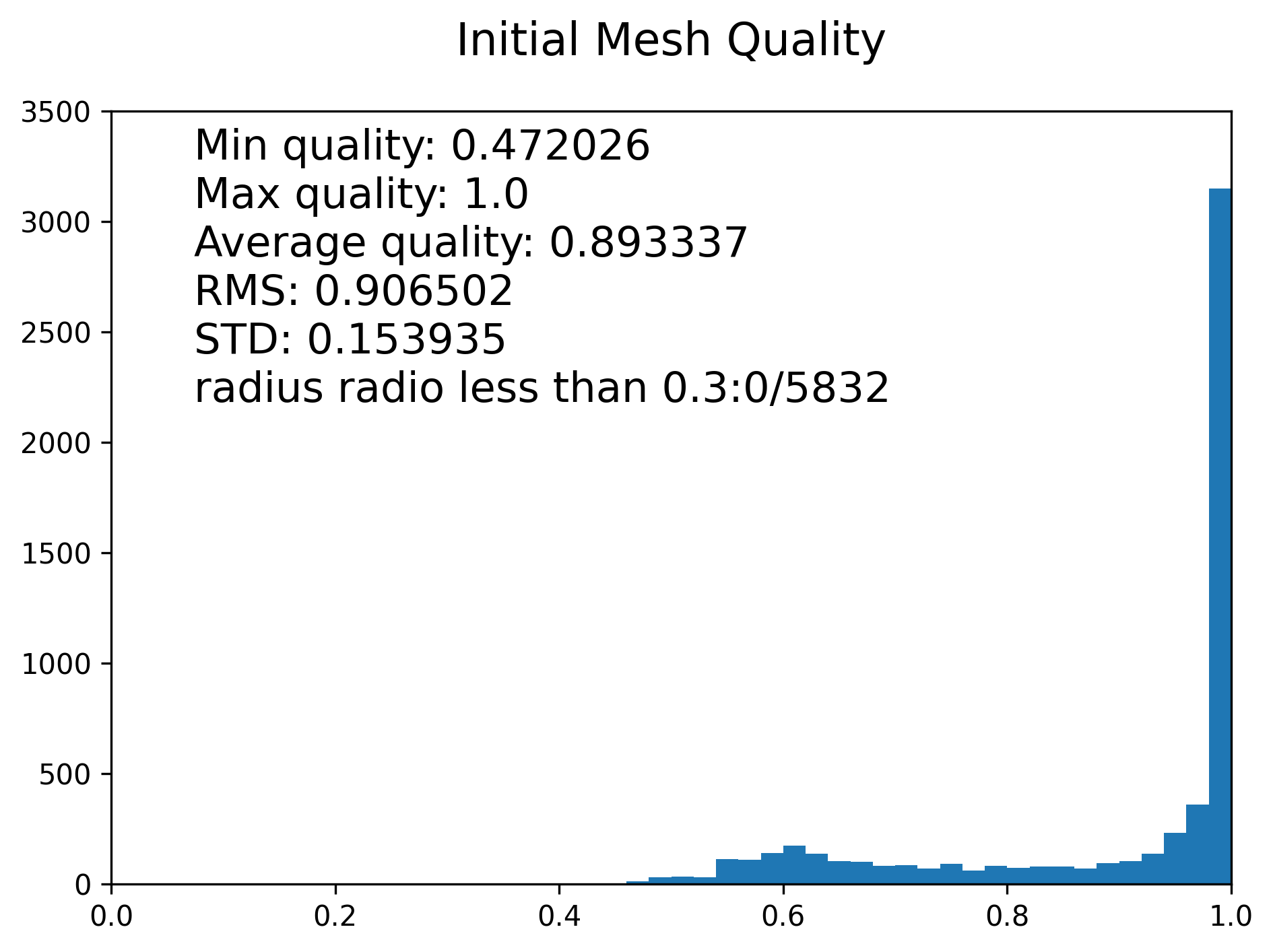}}
\hspace{0.01\linewidth}
\subfloat[RRE quality]{
\includegraphics[width=0.2\linewidth]{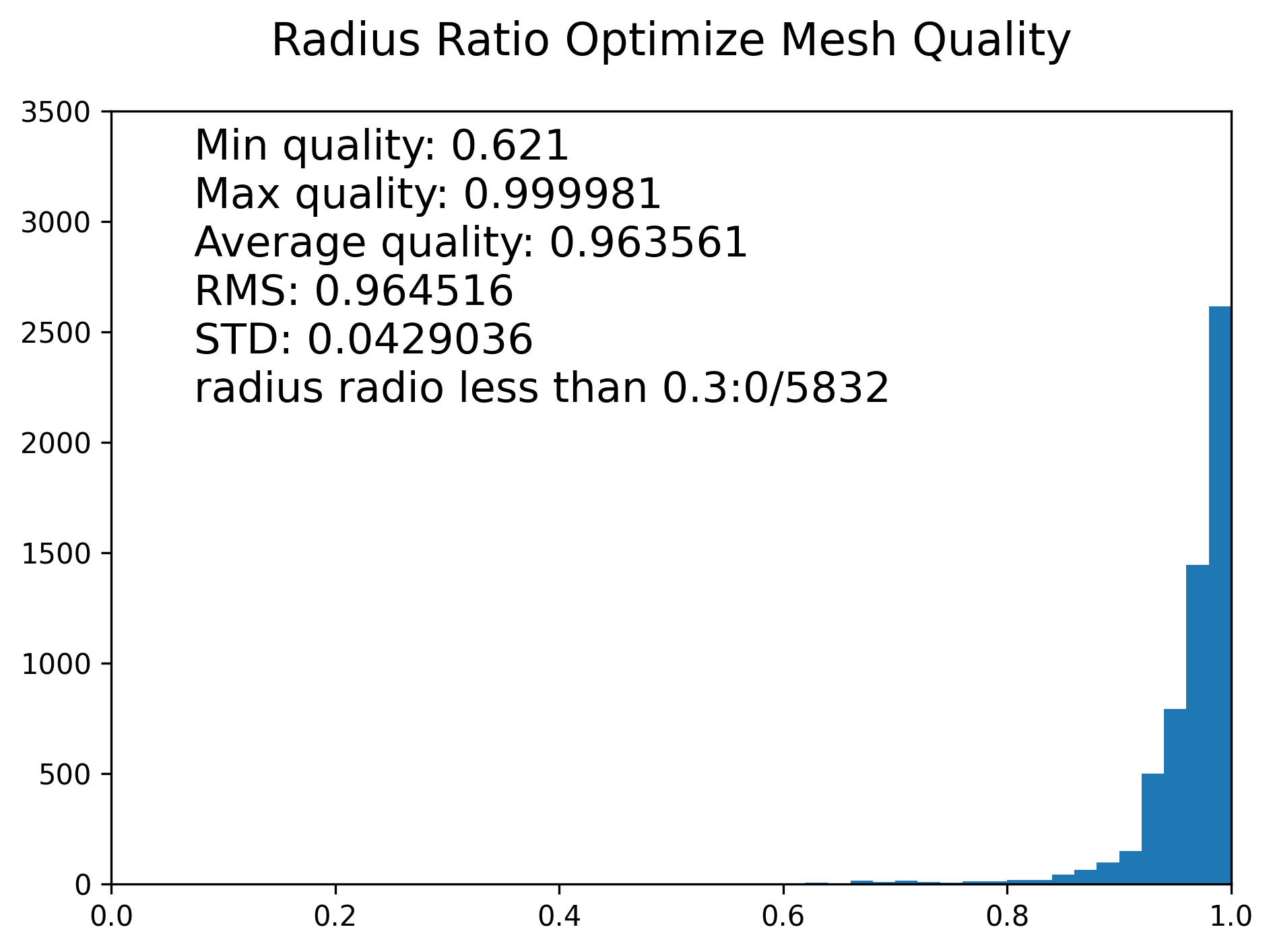}}
\caption{Curve-induced anisotropic mesh}
\label{fig:anisqu2d}
\end{figure}

Fig.~\ref{fig:anisqu2d} shows a square domain $[0,1]\times[0,1]$ with the
metric $\bm{M}(\bm{x})=\bm{Q}(\bm{x})^T\Lambda\bm{Q}(\bm{x})$,where
$$
\Lambda=\text{diag}(1,1+(t-1)\exp(-\frac{f(x,y)^2}{\sigma^2}),
f(x,y)=y-0.5-0.15\sin(2\pi x)
$$
with parameters $t=9,\sigma=0.12$. The first row of $\bm{Q}(\bm{x})$ is the unit 
tangent vector of $f(x,y)$, and the second row is the unit normal vector. This 
metric stretches the mesh along the curve $\Gamma: y=0.5+0.15\sin(2\pi x)$  and 
produces a smooth band structure.

\begin{figure}[htbp]
\centering
\subfloat[Init mesh]{
\includegraphics[width=0.2\linewidth]{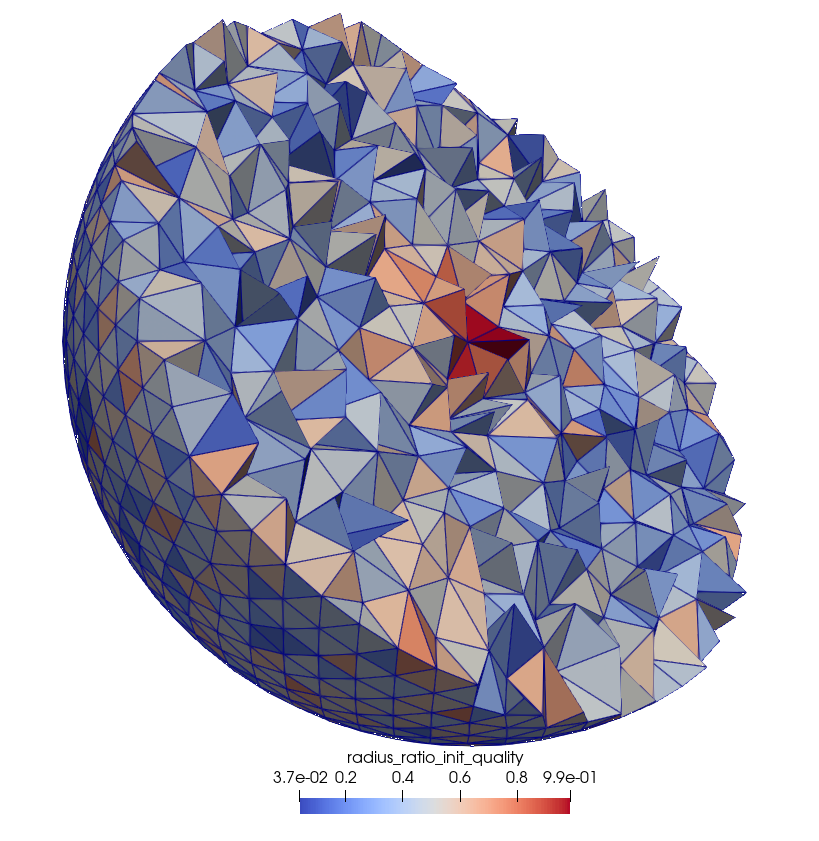}}
\hspace{0.01\linewidth}
\subfloat[Optimized mesh]{
\includegraphics[width=0.2\linewidth]{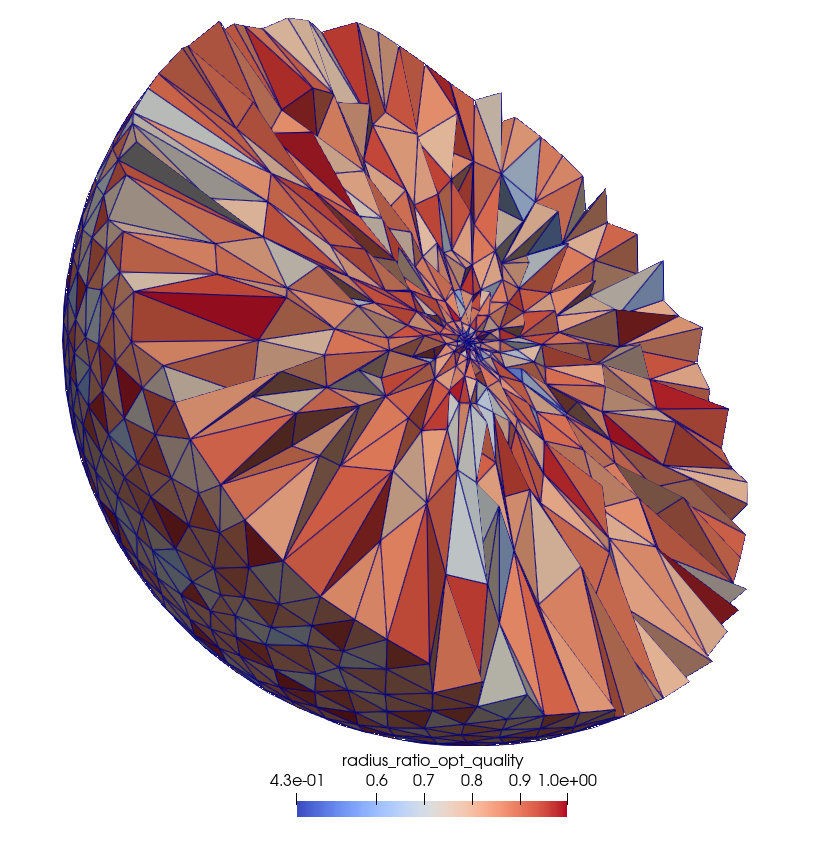}}
\vspace{0.01\linewidth}
\subfloat[Init quality]{
\includegraphics[width=0.2\linewidth]{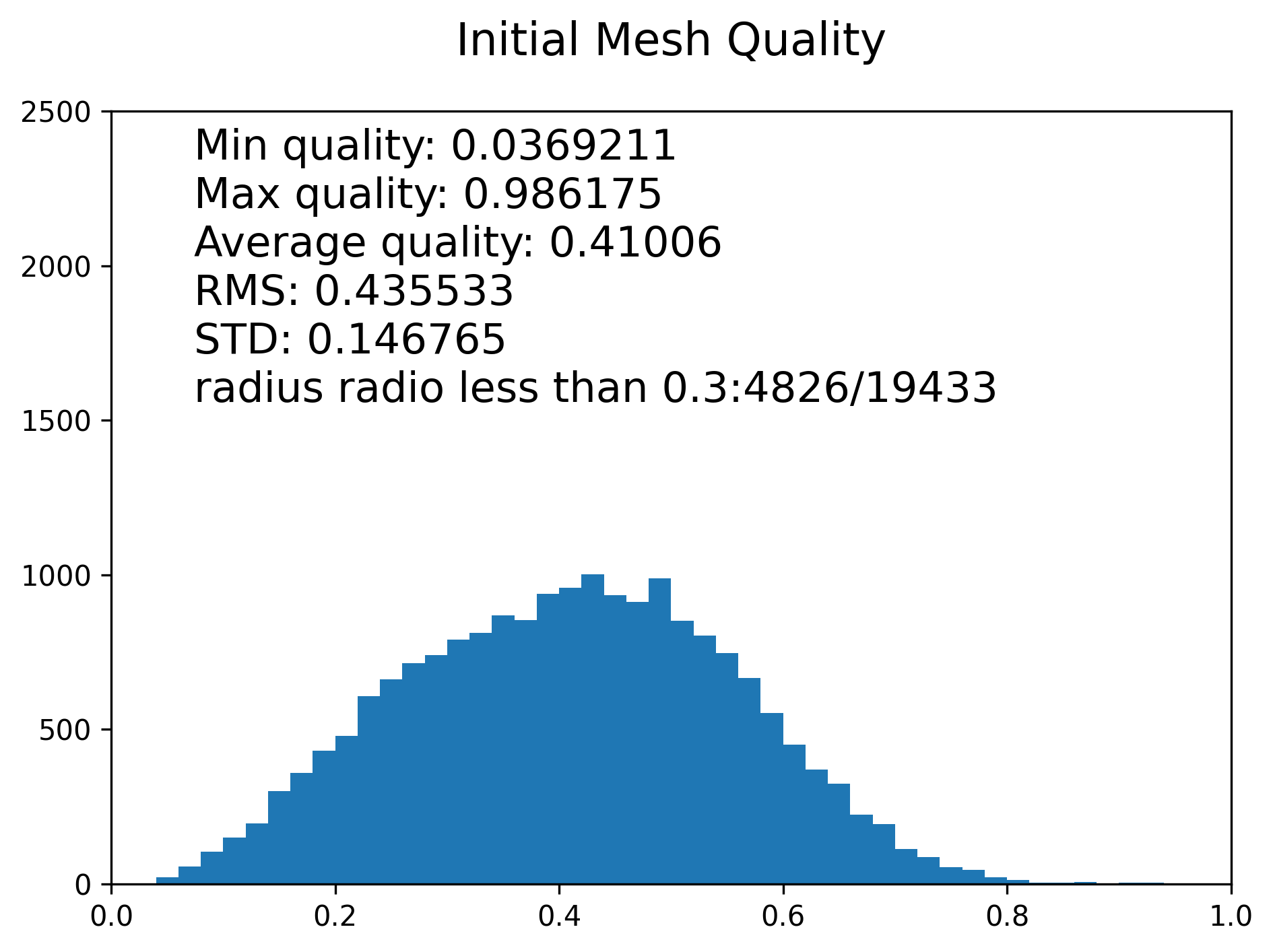}}
\hspace{0.01\linewidth}
\subfloat[RRE quality]{
\includegraphics[width=0.2\linewidth]{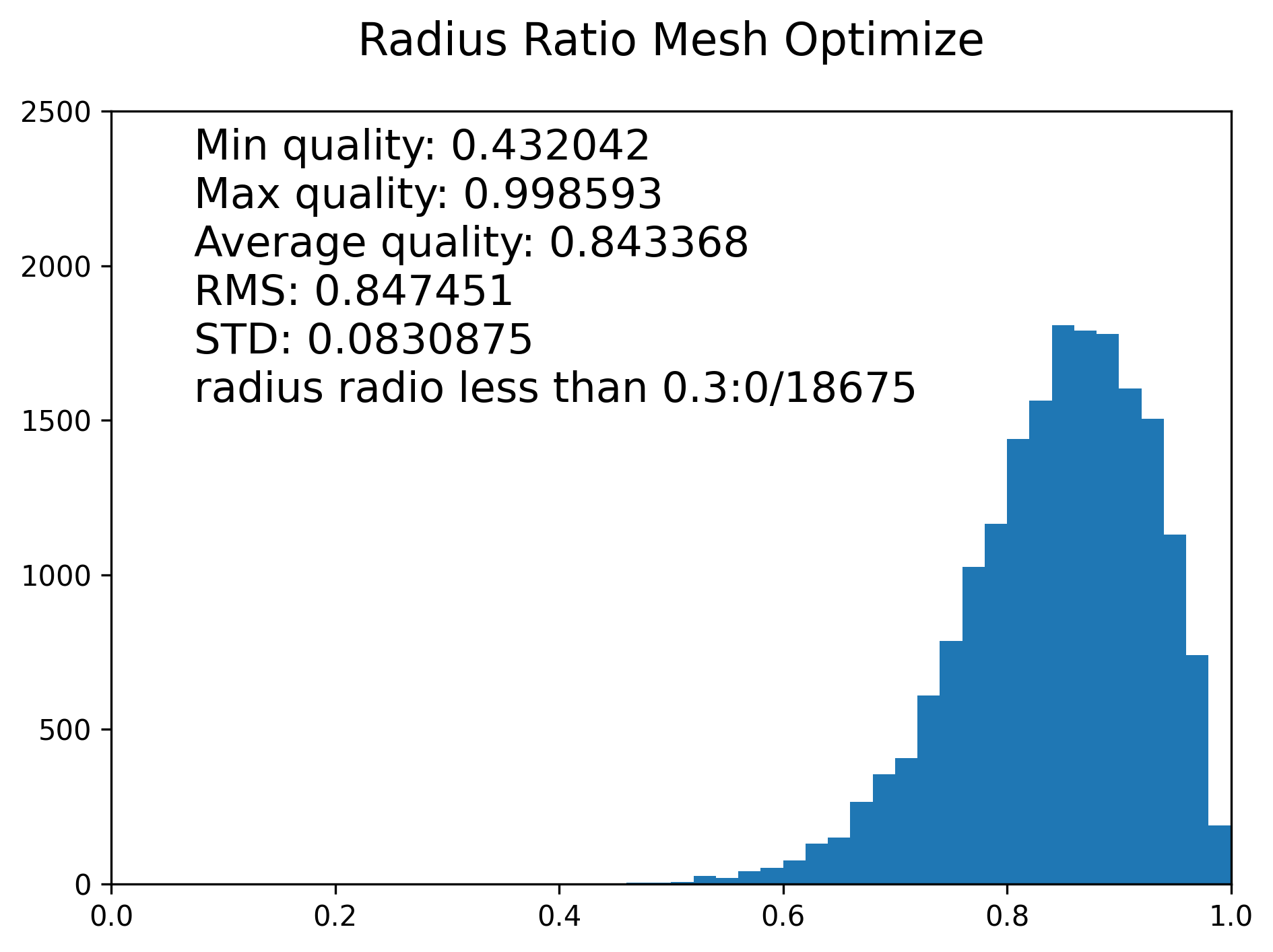}}
\caption{Unit sphere radial anisotropic mesh}
\label{fig:anispm1}
\end{figure}

Fig.~\ref{fig:anispm1} shows the unit sphere with the metric 
$\bm{M}(\bm{x}) = \lambda_r\bm{r}\bm{r}^T+\lambda_t(\bm{I}-\bm{r}\bm{r}^T)$,
where 
$$\bm{r}(\bm{x})=\frac{\bm{x}}{\|\bm{x}\|+\varepsilon},\lambda_r=1,\lambda_t=16$$ 
This metric stretches the mesh in the radial direction.

\begin{figure}[htbp] 
\centering
\subfloat[Init mesh]{
\includegraphics[width=0.2\linewidth]{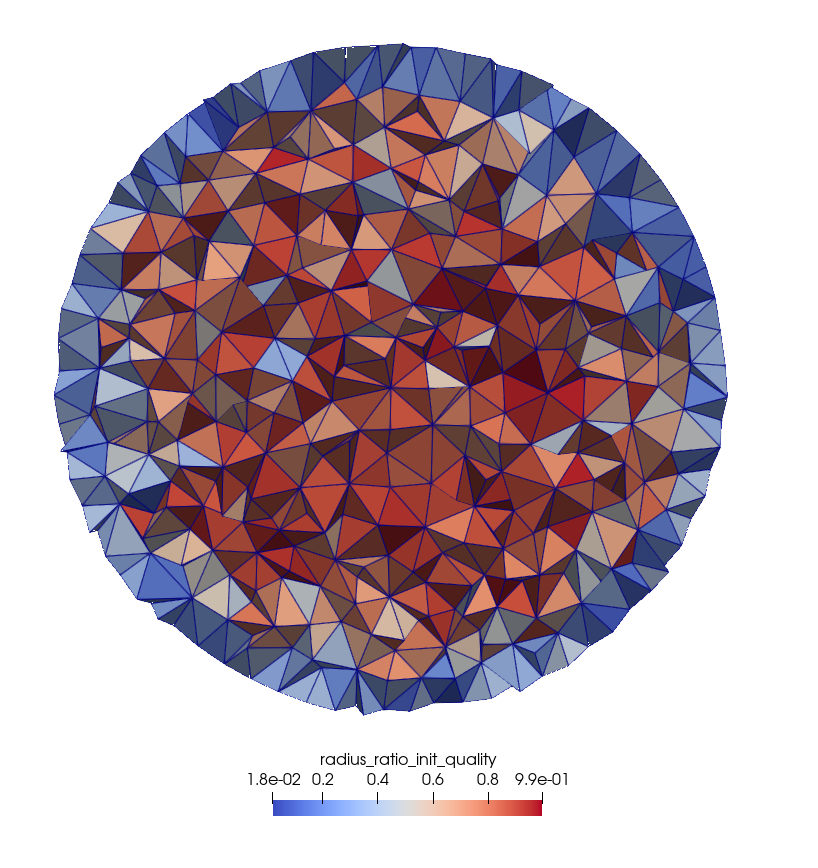}}
\hspace{0.01\linewidth}
\subfloat[Optimized mesh]{
\includegraphics[width=0.2\linewidth]{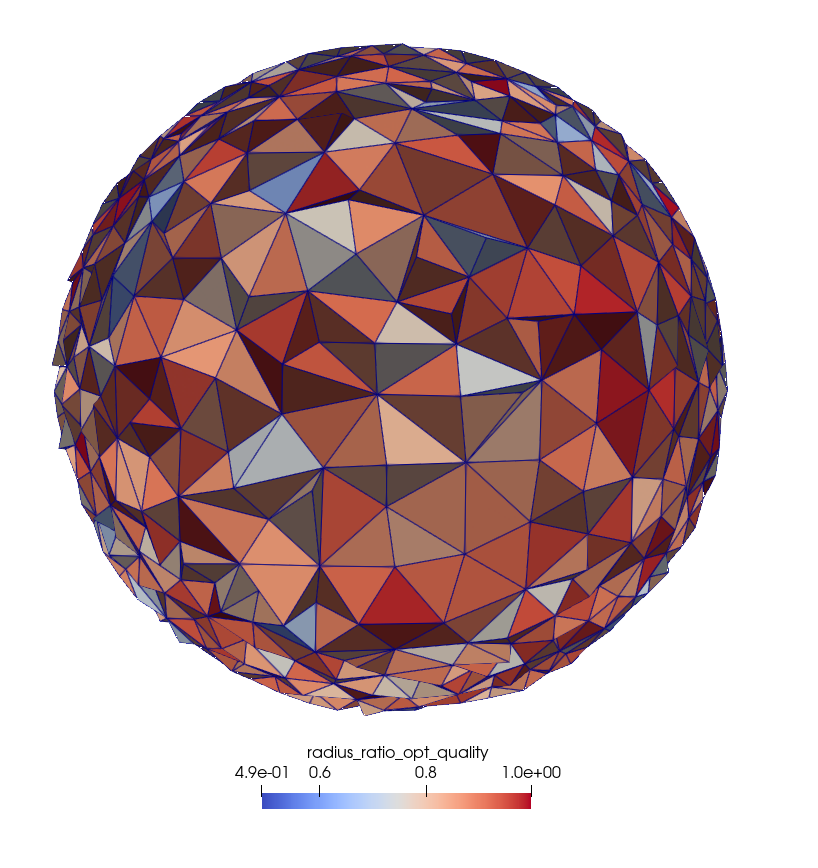}}
\vspace{0.01\linewidth}
\subfloat[Init quality]{
\includegraphics[width=0.2\linewidth]{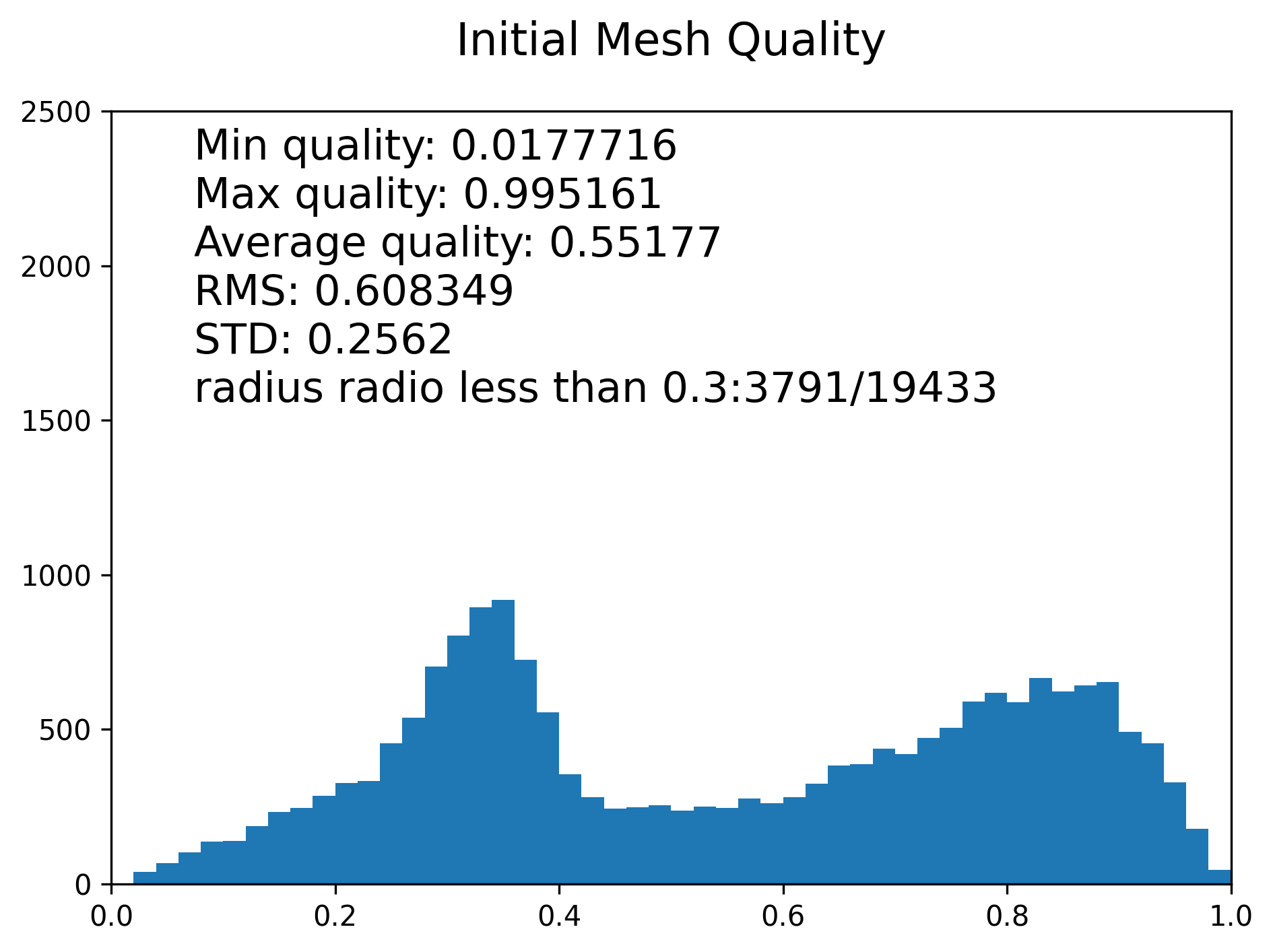}}
\hspace{0.01\linewidth}
\subfloat[RRE quality]{
\includegraphics[width=0.2\linewidth]{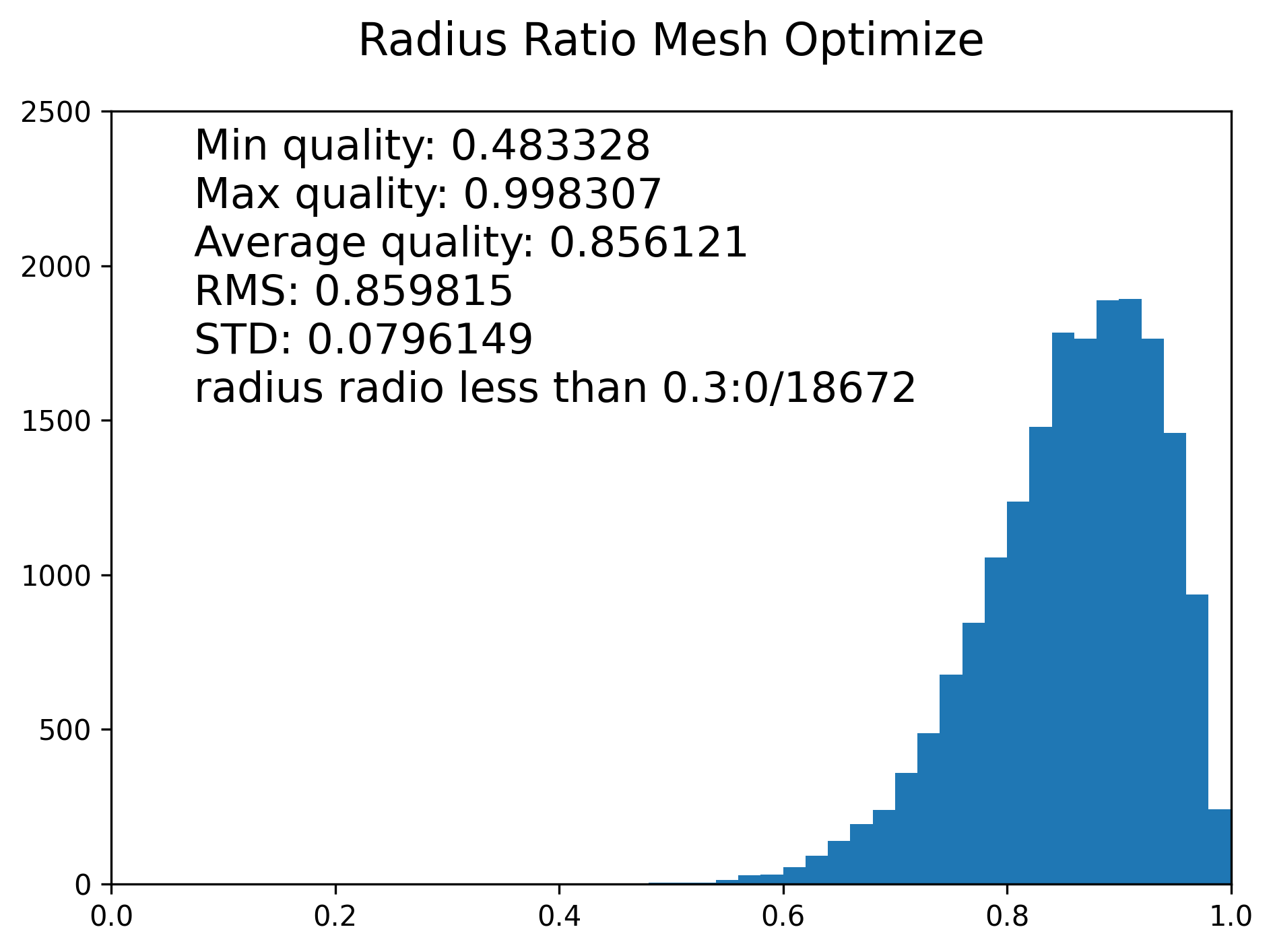}}
\caption{Unit sphere boundary layer anisotropic mesh}
\label{fig:anispm4}
\end{figure}

Fig.~\ref{fig:anispm4} shows another metric $\bm{M}(\bm{x}) 
= \lambda_n(d)\bm{r}\bm{r}^T+\lambda_t(\bm{I}-\bm{r}\bm{r}^T)$ on the unit 
sphere where
$$
d(\bm{x})=1-\|\bm{x}\|,\bm{r}(\bm{x})=\frac{\bm{x}}{\|\bm{x}\|+\varepsilon}, 
\lambda_n(d)=1+20{e}^{-10d},\lambda_t=1
$$
This metric produces a boundary layer structure with finite thickness.

\begin{figure}[htbp]
\centering
\subfloat[Init mesh]{
\includegraphics[width=0.2\linewidth]{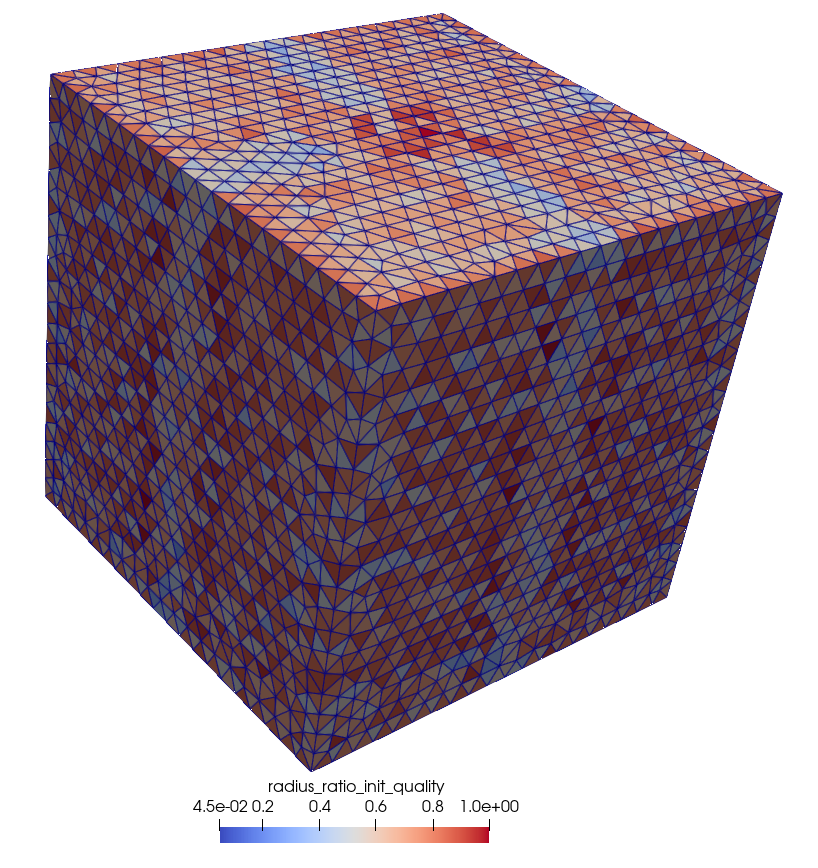}}
\hspace{0.01\linewidth}
\subfloat[Optimized mesh]{
\includegraphics[width=0.2\linewidth]{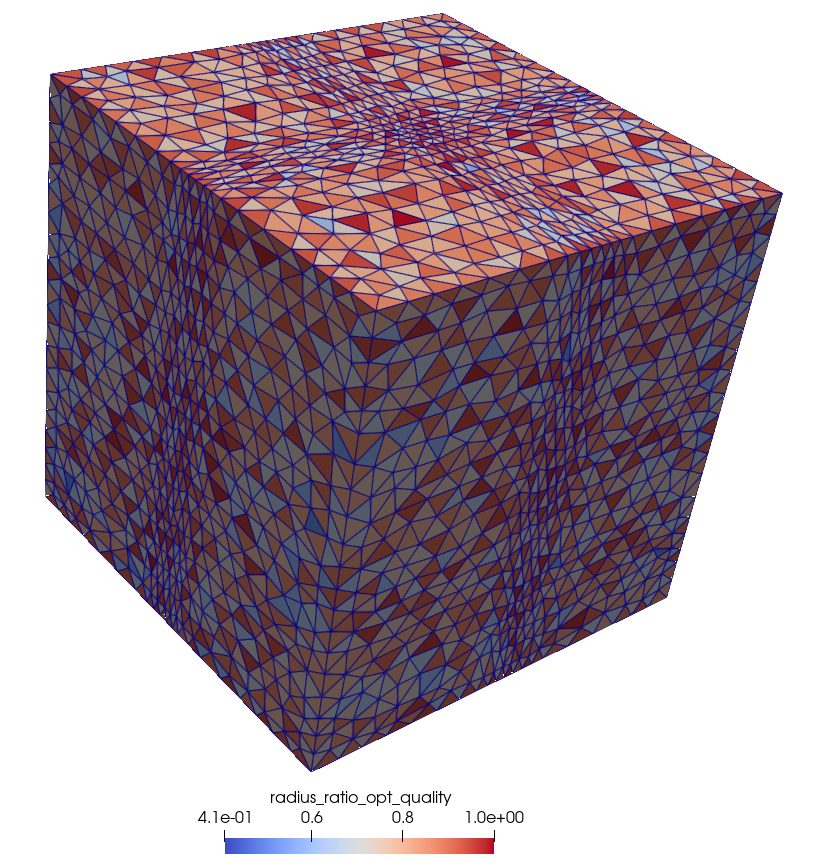}}
\vspace{0.01\linewidth}
\subfloat[Init quality]{
\includegraphics[width=0.2\linewidth]{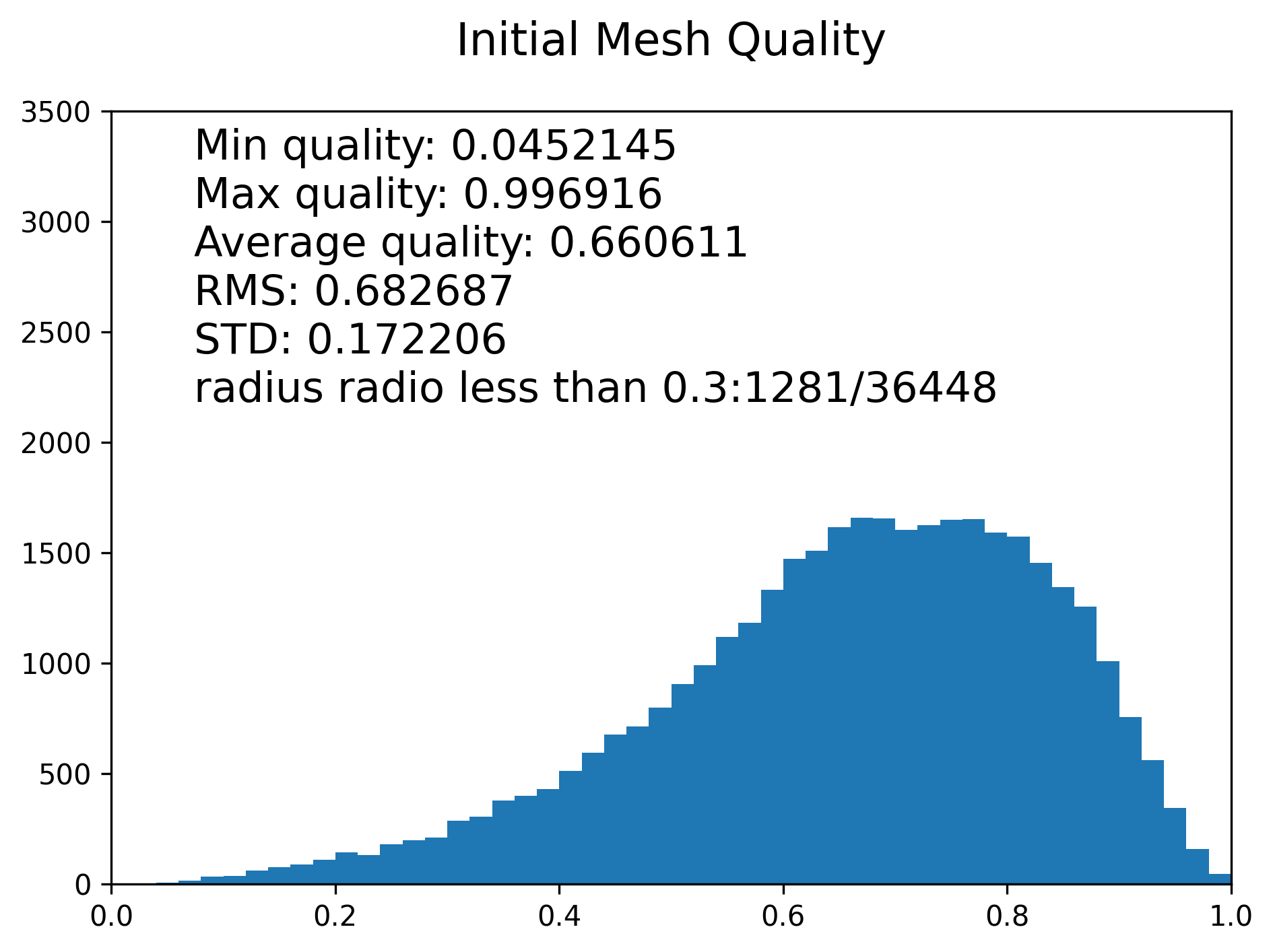}}
\hspace{0.01\linewidth}
\subfloat[RRE quality]{
\includegraphics[width=0.2\linewidth]{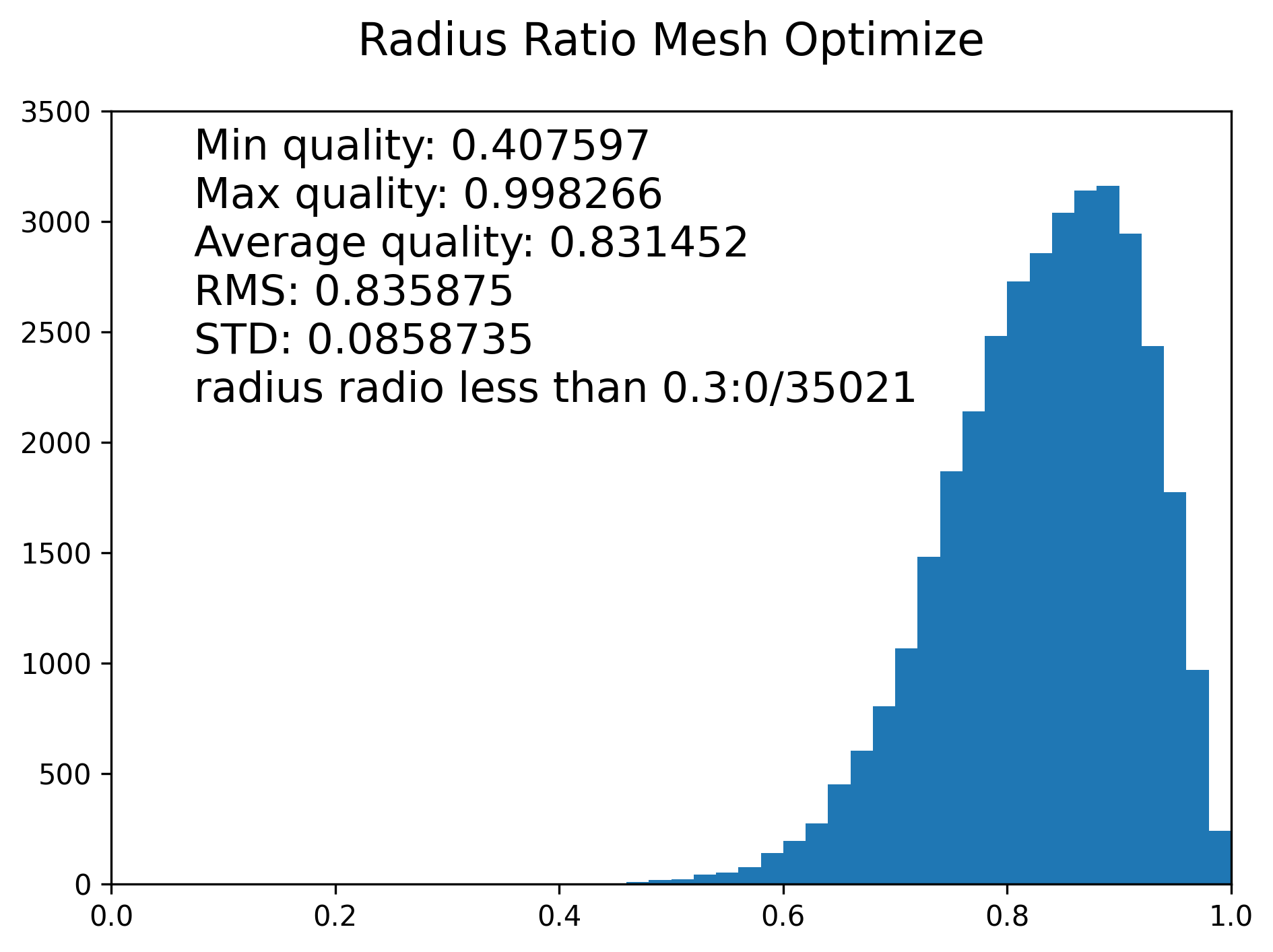}}
\caption{Cube anisotropic mesh}
\label{fig:anisqum1}
\end{figure}

Fig.~\ref{fig:anisqum1} shows the cube domain  $[0,1]^3$ with the metric
$\bm{M} = \text{diag}(m_0,m_1(\bm{x}),m_2(\bm{x}))$ where
$$
m_0=3,m_1(\bm{x})=1+(t-1)\exp(-\frac{(y-0.5)^2}{2\sigma^2}),
m_2(\bm{x})=1+(t-1)\exp(-\frac{(z-0.5)^2}{2\sigma^2})
$$
with $\sigma=0.05,t=10$. This metric introduces anisotropic regions near the 
planes $y=0.5$ and $z=0.5$, forming two intersecting anisotropic bands.

These examples show that the radius ratio based optimization method adapts well 
to anisotropic cases. For different metric tensors, the 
method captures both direction and scale changes, and aligns the mesh elements 
with the main anisotropic directions. In terms of element shape, the optimized 
tetrahedra tend to ideal shapes in the metric space. In strongly anisotropic 
regions, elements are stretched along the main directions, while in weak 
anisotropic regions, they return to near isotropic shapes. In all examples, 
the mesh quality is clearly improved. The minimum and average radius ratios 
increase, and low-quality elements are greatly reduced. This shows that the 
method has good shape optimization capabilities while maintaining anisotropic 
features.

\begin{table}[htbp]
    \caption{Anisotropic example optimization results}
    \label{tableani}
    \centering
    \begin{tabular}{|c|c|c|c|c|c|c|}
        \hline
        Model & Num of cell & Iteration count  & Time(sec.) \\
        \hline
        Unit disk: Spiral anisotropic mesh & 2972 & 203 & 0.70\\
        \hline
        Unit disk: Radial anisotropic mesh & 2972 & 182 & 0.65 \\
        \hline
        square: Curve-induced anisotropic mesh & 5832 & 148 & 1.10 \\
        \hline
        Unit sphere: radial anisotropic mesh& 18675 & 471 & 78.26 \\
        \hline
        Unit sphere: boundary layer anisotropic mesh& 18672 & 253 & 42.29 \\
        \hline
        Cube: anisotropic mesh & 35021 & 149 & 40.36 \\
        \hline
    \end{tabular}
\end{table}
Finally, Table~\ref{tableani} reports the number of cell, iteration counts, and 
running times for the two-dimensional and three-dimensional anisotropic cases. 
In general, the computation time increases with mesh size. The number of 
iterations depends on the metric. For strongly anisotropic cases (such as the 
radial metric on the unit sphere), the iteration count is higher. Even so, the 
method converges in a reasonable number of iterations in all cases. Overall, 
the proposed method shows good efficiency and stability for different scales 
and types of anisotropic problems.
\section{Conclusion and further work}
This paper proposes a simplex mesh optimization method based on the radius 
ratio energy. The method can effectively suppress slivers and improve 
the overall mesh quality. Based on the structure of the gradient of the radius 
ratio energy, we also construct a preconditioner. The preconditioned 
optimization reduces the number of iterations and speeds up convergence.We 
further extend the radius ratio energy to the anisotropic case and define an 
energy function suitable for anisotropic mesh. As a result, the RRE method 
can be applied to both isotropic and anisotropic mesh optimization.In future
work, we plan to extend the proposed method to high-order mesh and further 
improve computational efficiency by leveraging GPU-based parallel computing. 
\section*{Declaration of competing interest}
The authors declare that they have no known competing financial interests or 
personal relationships that could have appeared to influence the work reported 
in this paper.
\section*{Acknowledgement}
This work was supported in part by the National Key R\&D Program of China 
(2024YFA1012600), the National Natural Science Foundation of China (NSFC) 
(Grant No. 12371410, 12261131501), and the Graduate innovation Project of 
Xiangtan University(Grant No.XDCX2024Y183)
\bibliographystyle{abbrv}
\bibliography{references}  

@article{1ref,
  title={An apporach to automatic three-dimensional finite element mesh generation},
  author={Caendish, James C and Field, David A and Frey, William H},
  journal={International journal for numerical methods in engineering},
  volume={21},
  number={2},
  pages={329--347},
  year={1985},
  publisher={Wiley Online Library}
}

@incollection{2ref,
  title={Mesh generation and optimal triangulation},
  author={Bern, Marshall and Eppstein, David},
  booktitle={Computing in Euclidean geometry},
  pages={47--123},
  year={1995},
  publisher={World Scientific}
}

@article{5ref,
  title={Optimization of tetrahedral meshes based on element shape measures},
  author={Lo, SH},
  journal={Computers \& structures},
  volume={63},
  number={5},
  pages={951--961},
  year={1997},
  publisher={Elsevier}
}

@inproceedings{7ref,
  title={Aggressive tetrahedral mesh improvement},
  author={Klingner, Bryan Matthew and Shewchuk, Jonathan Richard},
  booktitle={Proceedings of the 16th international meshing roundtable},
  pages={3--23},
  year={2007},
  organization={Springer}
}

@article{10ref,
  title={Efficient mesh optimization schemes based on optimal Delaunay triangulations},
  author={Chen, Long and Holst, Michael},
  journal={Computer Methods in Applied Mechanics and Engineering},
  volume={200},
  number={9-12},
  pages={967--984},
  year={2011},
  publisher={Elsevier}
}

@article{11ref,
  title={Fast methods for computing centroidal Voronoi tessellations},
  author={Hateley, James C and Wei, Huayi and Chen, Long},
  journal={Journal of Scientific Computing},
  volume={63},
  number={1},
  pages={185--212},
  year={2015},
  publisher={Springer}
}

@article{12ref,
  title={Centroidal Voronoi tessellations: Applications and algorithms},
  author={Du, Qiang and Faber, Vance and Gunzburger, Max},
  journal={SIAM review},
  volume={41},
  number={4},
  pages={637--676},
  year={1999},
  publisher={SIAM}
}

@article{13ref,
  title={Finite volume methods},
  author={Eymard, Robert and Gallou{\"e}t, Thierry and Herbin, Rapha{\`e}le},
  journal={Handbook of numerical analysis},
  volume={7},
  pages={713--1018},
  year={2000},
  publisher={Elsevier}
}

@book{14ref,
  title={Finite element mesh generation},
  author={Lo, Daniel SH},
  year={2014},
  publisher={CRC press}
}

@article{15ref,
  title={A survey of unstructured mesh generation technology.},
  author={Owen, Steven J},
  journal={IMR},
  volume={239},
  number={267},
  pages={15},
  year={1998}
}

@inproceedings{16ref,
  title={Guaranteed-quality delaunay meshing in 3d (short version)},
  author={Chew, L Paul},
  booktitle={Proceedings of the thirteenth annual symposium on Computational geometry},
  pages={391--393},
  year={1997}
}

@article{19ref,
  title={Optimal delaunay triangulations},
  author={Chen, Long and Xu, Jin-chao},
  journal={Journal of Computational Mathematics},
  pages={299--308},
  year={2004},
  publisher={JSTOR}
}

@misc{20ref,
	title = {FEALPy: Finite Element Analysis Library in Python. https://github.com/weihuayi/
        fealpy},
	url = {https://github.com/weihuayi/fealpy},
	author = {Wei, Huayi and Huang, Yunqing},
    institution = {Xiangtan University},
	year = {Xiangtan University, 2017-2023},
}

@article{21ref,
  title={Gmsh: A 3-D finite element mesh generator with built-in pre-and post-processing facilities},
  author={Geuzaine, Christophe and Remacle, Jean-Fran{\c{c}}ois},
  journal={International journal for numerical methods in engineering},
  volume={79},
  number={11},
  pages={1309--1331},
  year={2009},
  publisher={Wiley Online Library}
}

@book{22ref,
  title={Numerical optimization},
  author={Nocedal, Jorge and Wright, Stephen J},
  year={1999},
  publisher={Springer}
}

@article{23ref,
  title={A survey of nonlinear conjugate gradient methods},
  author={Hager, William W and Zhang, Hongchao},
  journal={Pacific journal of Optimization},
  volume={2},
  number={1},
  pages={35--58},
  year={2006}
}

@article{24ref,
  title={On the limited memory BFGS method for large scale optimization},
  author={Liu, Dong C and Nocedal, Jorge},
  journal={Mathematical programming},
  volume={45},
  number={1},
  pages={503--528},
  year={1989},
  publisher={Springer}
}

@article{25ref,
  title={Note sur la convergence de m{\'e}thodes de directions conjugu{\'e}es},
  author={Polak, Elijah and Ribiere, Gerard},
  journal={Revue fran{\c{c}}aise d'informatique et de recherche op{\'e}rationnelle. S{\'e}rie rouge},
  volume={3},
  number={16},
  pages={35--43},
  year={1969},
  publisher={EDP Sciences}
}

@article{26ref,
  title={Methods of conjugate gradients for solving linear systems},
  author={Hestenes, Magnus R and Stiefel, Eduard and others},
  journal={Journal of research of the National Bureau of Standards},
  volume={49},
  number={6},
  pages={409--436},
  year={1952}
}

@article{27ref,
  title={Function minimization by conjugate gradients},
  author={Fletcher, Reeves and Reeves, Colin M},
  journal={The computer journal},
  volume={7},
  number={2},
  pages={149--154},
  year={1964},
  publisher={Oxford University Press}
}

@article{28ref,
  title={Lean algebraic multigrid (LAMG): Fast graph Laplacian linear solver},
  author={Livne, Oren E and Brandt, Achi},
  journal={SIAM Journal on Scientific Computing},
  volume={34},
  number={4},
  pages={B499--B522},
  year={2012},
  publisher={SIAM}
}

@incollection{29ref,
  title={Algebraic multigrid},
  author={Ruge, John W and St{\"u}ben, Klaus},
  booktitle={Multigrid methods},
  pages={73--130},
  year={1987},
  publisher={SIAM}
}

@article{30ref,
  title={Conjugate gradient method},
  author={Nazareth, John L},
  journal={Wiley Interdisciplinary Reviews: Computational Statistics},
  volume={1},
  number={3},
  pages={348--353},
  year={2009},
  publisher={Wiley Online Library}
}

@article{31ref,
  title={Generation of three-dimensional unstructured grids by the advancing-front method},
  author={L{\"o}hner, Rainald and Parikh, Paresh},
  journal={International Journal for Numerical Methods in Fluids},
  volume={8},
  number={10},
  pages={1135--1149},
  year={1988},
  publisher={Wiley Online Library}
}

@article{32ref,
  title={Automatic three-dimensional mesh generation by the finite octree technique},
  author={Shephard, Mark S and Georges, Marcel K},
  journal={International Journal for Numerical methods in engineering},
  volume={32},
  number={4},
  pages={709--749},
  year={1991},
  publisher={Wiley Online Library}
}

@article{33ref,
  title={A modified quadtree approach to finite element mesh generation},
  author={Yerry, Mark and Shephard, Mark},
  journal={IEEE Computer Graphics and Applications},
  volume={3},
  number={01},
  pages={39--46},
  year={1983},
  publisher={IEEE Computer Society}
}

@article{34ref,
  title={Unstructured mesh generation},
  author={Shewchuk, Jonathan Richard},
  journal={Combinatorial Scientific Computing},
  volume={12},
  number={257},
  pages={2},
  year={2012},
  publisher={CRC Press Boca Raton, FL}
}

@article{36ref,
  title={Introducing the target-matrix paradigm for mesh optimization via node-movement},
  author={Knupp, Patrick},
  journal={Engineering with Computers},
  volume={28},
  number={4},
  pages={419--429},
  year={2012},
  publisher={Springer}
}

@techreport{37ref,
  title={Metric type in the target-matrix mesh optimization paradigm},
  author={Knupp, Patrick},
  year={2020},
  institution={Lawrence Livermore National Lab.(LLNL), Livermore, CA (United States)}
}

@article{38ref,
  title={Sliver exudation},
  author={Cheng, Siu-Wing and Dey, Tamal K and Edelsbrunner, Herbert and Facello, Michael A and Teng, Shang-Hua},
  journal={Journal of the ACM (JACM)},
  volume={47},
  number={5},
  pages={883--904},
  year={2000},
  publisher={ACM New York, NY, USA}
}

@article{39ref,
  title={A nonlinear conjugate gradient method with a strong global convergence property},
  author={Dai, Yu-Hong and Yuan, Yaxiang},
  journal={SIAM Journal on optimization},
  volume={10},
  number={1},
  pages={177--182},
  year={1999},
  publisher={SIAM}
}

@inproceedings{40ref,
  title={Mesh Smoothing Schemes Based on Optimal Delaunay Triangulations.},
  author={Chen, Long},
  booktitle={IMR},
  pages={109--120},
  year={2004},
  organization={Citeseer}
}

@book{42sliver,
  title={Sliver-free three dimensional delaunay mesh generation},
  author={Li, Xiangyang},
  year={2001},
  publisher={University of Illinois at Urbana-Champaign}
}

@article{43ref,
  title={Two discrete optimization algorithms for the topological improvement of tetrahedral meshes},
  author={Shewchuk, Jonathan Richard},
  journal={Unpublished manuscript},
  volume={65},
  pages={2--7},
  year={2002}
}

@book{44ref,
  title={Delaunay mesh generation},
  author={Cheng, Siu-Wing and Dey, Tamal Krishna and Shewchuk, Jonathan and Sahni, Sartaj},
  year={2013},
  publisher={CRC Press Boca Raton}
}

@article{45ref,
  title={On the design of CGAL a computational geometry algorithms library},
  author={Fabri, Andreas and Giezeman, Geert-Jan and Kettner, Lutz and Schirra, Stefan and Sch{\"o}nherr, Sven},
  journal={Software: Practice and Experience},
  volume={30},
  number={11},
  pages={1167--1202},
  year={2000},
  publisher={Wiley Online Library}
}
\end{document}